\documentclass[12pt]{elsarticle}
\usepackage{graphicx} 

\usepackage[margin=1in]{geometry}
\usepackage{amsmath}  
\usepackage{amsthm}
\usepackage{amsfonts}
\usepackage{amssymb}
\usepackage{graphicx}  
\usepackage{float}
\usepackage{natbib}
\usepackage{hyperref}
\usepackage{bm}

\usepackage{caption}
\usepackage{subcaption}
\usepackage{tikz}

\theoremstyle{definition}

\newtheorem{remark}{Remark}

\newcommand{\fref}[1]{Figure~\ref{#1}}
\newcommand{\tref}[1]{Table~\ref{#1}}
\newcommand{\sref}[1]{Section~\ref{#1}}
\newcommand{\vm}[1]{\bm{#1}}
\newcommand{\vx}{\vm{x}}

\newcommand{\proj}[1]{\Pi^{\varepsilon} #1}

\newcommand{\diverge}[1]{\nabla \cdot #1}

\newcommand{\symgrad}[1]{\nabla_s #1}

\newcommand{\Nmatrix}{\bm{N}^{\partial E}}
\newcommand{\Projbeta}{\Pi_{\beta}\vm{\sigma}}
\newcommand{\Projbetavect}{\overline{\Projbeta}}
\newcommand{\symP}{\mathbb{P}}
\newcommand{\symPspace}{[\mathbb{P}_\ell(E)]^{2\times 2}_{\text{sym}}}

\newcommand{\veps}{\vm{\varepsilon}}
\newcommand{\vsigma}{\vm{\sigma}}
\newcommand{\vd}{\vm{d}}

\newcommand{\wtilde}[1]{\widetilde{#1}}

\newcommand{\ac}{}

\hypersetup{
    colorlinks,
    citecolor=black,
    filecolor=black,
    linkcolor=black,
    urlcolor=black
}

\graphicspath{ 
  {./figs/}}
  
\begin{document}
\begin{frontmatter}

\title{Stress-hybrid virtual element method on six-noded triangular meshes for compressible and nearly-incompressible linear elasticity}

\cortext[cor1]{Corresponding authors}

\author[inst1]{Alvin Chen}

\affiliation[inst1]{organization={Department of Mathematics},
            addressline={University of California}, 
            city={Davis},
            postcode={95616}, 
            state={CA},
            country={USA}}

\author[inst2]{Joseph E. Bishop}
\affiliation[inst2]{organization={Engineering Sciences Center},
            addressline={Sandia National Laboratories}, 
            city={Albuquerque},
            postcode={87185},
            state={NM},
            country={USA}}

\author[inst3]{N. Sukumar\corref{cor1}}
\ead{nsukumar@ucdavis.edu}
\affiliation[inst3]{organization={Department of Civil and Environmental Engineering},
            addressline={University of California}, 
            city={Davis},
            postcode={95616}, 
            state={CA},
            country={USA}}

\begin{abstract}
In this paper, we present a first-order Stress-Hybrid Virtual Element Method (SH-VEM) on six-noded triangular meshes for linear plane elasticity. We adopt the Hellinger--Reissner variational principle to construct a weak equilibrium condition and a stress based projection operator. \ac{In each element, the stress projection operator is expressed in terms of the nodal displacements, which leads to a displacement based formulation.}
This stress-hybrid approach assumes a globally continuous displacement field while the stress field is discontinuous across each element. The stress field is initially represented by divergence-free tensor polynomials based on Airy stress functions, 
\ac{but we also present a formulation that uses a penalty term to enforce the element equilibrium conditions, referred to as the Penalty Stress-Hybrid Virtual Element Method (PSH-VEM)}. 
Numerical results are presented for PSH-VEM and SH-VEM, and we compare their convergence to the composite triangle FEM and B-bar VEM on benchmark problems in linear elasticity. The SH-VEM converges optimally in the $L^2$ norm of the displacement, energy seminorm, and the $L^2$ norm of hydrostatic stress. Furthermore, the results reveal that PSH-VEM converges in most cases at a faster rate than the expected optimal rate, but it requires the selection of a suitably chosen penalty parameter.
\end{abstract}


\end{frontmatter}

\section{Introduction}\label{sec:intro}
\ac{Quadrilaterals and hexahedra are common elements used in many applications of the Finite Element Method (FEM).} 
However, high-quality quadrilateral and hexahedral meshes cannot be automatically generated for complex geometries. Triangular and tetrahedral meshes have well-established automatic mesh generators, but these elements tend to be overly stiff and suffer from volumetric locking in the nearly-incompressible limit. Many techniques such as the B-bar and selective integration formulations~\cite{Malkus:1978:MFE,Piltner:1999:ASC}, method of incompatible modes~\cite{Wilson:1973:IDM}, assumed enhanced strain~\cite{Simo:1986:jam, Simo:1990:ijnme}, stabilized elements~\cite{Belytschko:1986:EIO, Belytschko:1991:cmame}, mixed elements~\cite{Veubeke:1965:DAE}, and hybrid stress methods~\cite{Jog:2005:fead,Jog:2010:mom,Jog:2006:ijnme,Pian:1964:DOE, Pian:1984:ijnme, Pian:1986:RBI} have been developed to alleviate locking in the nearly-incompressible limit. However, many of these methods are only applicable to quadrilateral and hexahedral elements. Progress has been made to construct modified triangular elements that are robust and accurate: the variational multiscale approach~\cite{Scovazzi:2016:SSA}, hybrid-stress elements using Airy stress functions~\cite{KarimiPour:2023:ASF,Li:2016:HHS,Ma:2013:R18,Rezaiee-Pajand:2020:TSB,Wang:2017:ASQ}, B\'ezier elements~\cite{Kadapa:2019:NQB}, F-bar and reduced integration~\cite{Danielson:2014:FNT,Neto:2005:FBL} to name a few. Another promising approach is the use of composite elements~\cite{Camacho:1996:CMI,Foulk:2021:ECT,Guo:2000:TCF,Leonetti:2015:ACM,Nguyen:2019:PCF,Thoutireddy:2002:TCF}. In this approach, an element is split into sub-triangles (sub-tetrahedra), and then on each sub-triangle the strain is approximated and combined to recover the strain field over the entire element. Similar constructions are also used in~\cite{Boerner:2008:AME} to construct a three-dimensional brick element for nearly incompressible nonlinear elasticity problems and in~\cite{Cook:1975:APS,Cook:1987:APH,Cook:1990:SOF,Eom:2009:AMP} for triangular bending elements. 

The recent development of the \ac{Virtual Element Method} (\ac{VEM})~\cite{basicprinciple,elasticdaveiga,Veiga2014TheHG} has provided an alternative way to extend these methods to alleviate volumetric locking for general polygons and polyhedra. Methods such as mixed and hybrid formulations~\cite{Artioli:2018:cmame,Artioli:2020:M2AN,Caceres:2019:anm,Dassi:2021:m3as}, B-bar~\cite{Park:2020:meccanica}, assumed strain~\cite{berrone:2023:los,Chen:2023:SFS,Chen:2023:SFV,enhanced:VEM,Lamperti:2023:cmech}, and recently hybrid stress for general quadrilaterals~\cite{chen:2023:shv} have been adopted in the virtual element method. The approach of using triangular sub-elements has also been applied in~\cite{Wriggers:2017:EVE} to construct a stabilization term for the virtual element method.

In this paper, we focus on nearly-incompressible two-dimensional linear elasticity problems; 
however, extensions to three-dimensional problems appears to be feasible. We first propose a virtual element method for a six-noded element using \ac{equilibrated} stress fields originally motivated for finite elements 
in~\cite{Cen:2011:cmame,Pian:1964:DOE,Spilker:1981:PIH}. This method is applicable to both convex and nonconvex elements; however, in our application we focus mainly on triangular meshes with an additional node placed along each edge of the triangle. Following~\cite{chen:2023:shv}, we refer to this approach as the Stress-Hybrid Virtual Element Method (SH-VEM). This approach
draws on the stress-hybrid finite element method of Pian and Sumihara~\cite{Pian:1984:ijnme} and the resulting method resembles the stress recovery procedure of Artioli et al.~\cite{Artioli:2019:ESR}. Unlike the standard mixed or hybrid methods based on the Hellinger--Reissner principle, which require the use of more specialized approximation spaces and additional traction degrees of freedom, the stress-hybrid approach uses a discontinuous stress field that is only dependent on the nodal displacement degrees of freedom.
When using a fifteen-parameter expansion, the resulting method does not require a stabilization term and shows immunity to shear and volumetric locking. 

The remainder of this paper is structured as follows. In~\sref{sec:formulation}, we use the Hellinger--Reissner variational principle to construct a weak formulation of the linear elasticity problem. In~\sref{sec:VEM}, we construct the necessary projection operators and present the virtual element space. In~\sref{sec:implementation}, we present the numerical implementation of the projections, discuss possible choices of stress basis functions, and construct the element stiffness matrix.
In~\sref{sec:penalty}, we present a virtual element method, \ac{referred to as the Penalty Stress-Hybrid Virtual Element Method (PSH-VEM)}, that combines the stress-hybrid VEM with an equilibrium penalty term in order to construct an element that does not rely on the Airy stress functions, alleviates shear locking, and requires no stabilization term. 
In~\sref{sec:results}, we present numerical results for Composite Triangle FEM (CT FEM)~\cite{Guo:2000:TCF}, B-bar VEM~\cite{Park:2020:meccanica}, SH-VEM, and PSH-VEM. 
We present results for the four methods on a series of benchmark problems in plane strain 
linear elasticity in the nearly-incompressible limit: bending of a thin cantilever beam, Cook's membrane, infinite plate with a circular hole, hollow pressurized cylinder, a manufactured problem with exact solution, and the punch problem. We close with our main findings and some final remarks in~\sref{sec:conclusion}. 

\section{Variational formulation}\label{sec:formulation}
Consider an elastic body occupying a bounded region $\Omega \subset \mathbb{R}^2$. Assume that the boundary $\partial \Omega$ can be decomposed into a disjoint union $\partial \Omega = \Gamma_u \cup \Gamma_t$ and $\Gamma_u \cap \Gamma_t =\emptyset$. On $\Gamma_u$, displacement boundary conditions $\vm{u}=\vm{u}_0$ are imposed, while tractions $\vm{t}=\bar{\vm{t}}$ are applied on $\Gamma_t$. The displacement field is denoted by $\vm{u} \in [H^1(\Omega)]^2$, the small-strain tensor $\vm{\varepsilon} = \symgrad{\vm{u}}$, the Cauchy stress tensor $\vm{\sigma} \in [L^2(\Omega)]^{2\times2}_{\textrm{sym}}$, the material moduli tensor is $\mathbb{C}$, and $\vm{b}\in [L^2(\Omega)]^2$ represents the body force per unit volume.

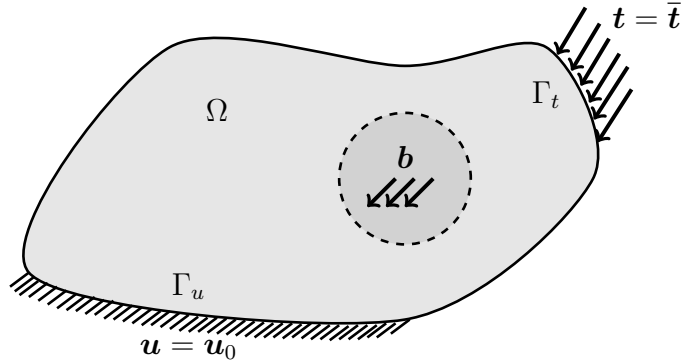
\begin{figure}[!htb]
\centering
\begin{tikzpicture}[scale=2.5]

            \filldraw[line width=1pt, black,fill=gray!20] plot [smooth cycle] coordinates {(-2,.5) (-1.25,1.7) (0,1.6) (.75,1.7) (1,1) (0,.25) };

            \draw[line width=1pt] (-2,.5) -- ++(-.05-.05,-.025-.05);
            
            \draw[line width=1pt] (-2+.025,.5-.025) -- ++(-.05-.05,-.025-.05);
            \draw[line width=1pt] (-2+.05,.5-.045) -- ++(-.05-.05,-.025-.05);
            \draw[line width=1pt] (-2+.1,.5-.055) -- ++(-.05-.05,-.025-.05);
            \draw[line width=1pt] (-2+.15,.5-.075) -- ++(-.05-.05,-.025-.05);
            \draw[line width=1pt] (-2+.2,.5-.095) -- ++(-.05-.05,-.025-.05);

            \draw[line width=1pt] (-2+.25,.5-.115) -- ++(-.05-.05,-.025-.05);
            \draw[line width=1pt] (-2+.3,.5-.125) -- ++(-.05-.05,-.025-.05);
            \draw[line width=1pt] (-2+.35,.5-.13) -- ++(-.05-.05,-.025-.05);
            \draw[line width=1pt] (-2+.4,.5-.145) -- ++(-.05-.05,-.025-.05);
            \draw[line width=1pt] (-2+.45,.5-.15) -- ++(-.05-.05,-.025-.05);

            \draw[line width=1pt] (-2+.5,.5-.17) -- ++(-.05-.05,-.025-.05);
            \draw[line width=1pt] (-2+.55,.5-.175) -- ++(-.05-.05,-.025-.05);
            \draw[line width=1pt] (-2+.6,.5-.185) -- ++(-.05-.05,-.025-.05);
            \draw[line width=1pt] (-2+.65,.5-.195) -- ++(-.05-.05,-.025-.05);
            \draw[line width=1pt] (-2+.7,.5-.205) -- ++(-.05-.05,-.025-.05);

            \draw[line width=1pt] (-2+.75,.5-.205) -- ++(-.05-.05,-.025-.05);
            \draw[line width=1pt] (-2+.8,.5-.215) -- ++(-.05-.05,-.025-.05);
            \draw[line width=1pt] (-2+.85,.5-.225) -- ++(-.05-.05,-.025-.05);
            \draw[line width=1pt] (-2+.9,.5-.235) -- ++(-.05-.05,-.025-.05);
            \draw[line width=1pt] (-2+.95,.5-.24) -- ++(-.05-.05,-.025-.05);

            \draw[line width=1pt] (-2+1,.5-.24) -- ++(-.05-.05,-.025-.05);
            \draw[line width=1pt] (-2+1.05,.5-.245) -- ++(-.05-.05,-.025-.05);
            \draw[line width=1pt] (-2+1.1,.5-.25) -- ++(-.05-.05,-.025-.05);
            \draw[line width=1pt] (-2+1.15,.5-.255) -- ++(-.05-.05,-.025-.05);
            \draw[line width=1pt] (-2+1.2,.5-.26) -- ++(-.05-.05,-.025-.05);

            \draw[line width=1pt] (-2+1.25,.5-.265) -- ++(-.05-.05,-.025-.05);
            \draw[line width=1pt] (-2+1.3,.5-.27) -- ++(-.05-.05,-.025-.05);
            \draw[line width=1pt] (-2+1.35,.5-.27) -- ++(-.05-.05,-.025-.05);
            \draw[line width=1pt] (-2+1.4,.5-.27) -- ++(-.05-.05,-.025-.05);
            \draw[line width=1pt] (-2+1.45,.5-.275) -- ++(-.05-.05,-.025-.05);

            \draw[line width=1pt] (-2+1.5,.5-.27) -- ++(-.05-.05,-.025-.05);
            \draw[line width=1pt] (-2+1.55,.5-.275) -- ++(-.05-.05,-.025-.05);
            \draw[line width=1pt] (-2+1.6,.5-.28) -- ++(-.05-.05,-.025-.05);
            \draw[line width=1pt] (-2+1.65,.5-.28) -- ++(-.05-.05,-.025-.05);
            \draw[line width=1pt] (-2+1.7,.5-.285) -- ++(-.05-.05,-.025-.05);

            \draw[line width=1pt] (-2+1.75,.5-.28) -- ++(-.05-.05,-.025-.05);
            \draw[line width=1pt] (-2+1.8,.5-.28) -- ++(-.05-.05,-.025-.05);
            \draw[line width=1pt] (-2+1.85,.5-.28) -- ++(-.05-.05,-.025-.05);
            \draw[line width=1pt] (-2+1.9,.5-.27) -- ++(-.05-.05,-.025-.05);
            \draw[line width=1pt] (-2+1.95,.5-.27) -- ++(-.05-.05,-.025-.05);

            \draw[line width=1pt] (0+.025,.25) -- ++(-.05-.05,-.025-.05);;

            \draw[line width=1.5pt , ->] (.75+.21,1.7+.21) -- ++(-.15,-.25);
            \draw[line width=1.5pt , ->] (.75+.21+.075,1.7+.21-.085) -- ++(-.15,-.25);
            \draw[line width=1.5pt , ->] (.75+.21+.125,1.7+.21-.165) -- ++(-.15,-.25);
            \draw[line width=1.5pt , ->] (.75+.21+.18,1.7+.21-.245) -- ++(-.16,-.26);
            \draw[line width=1.5pt , ->] (.75+.21+.225,1.7+.21-.315) -- ++(-.18,-.28);
            \draw[line width=1.5pt , ->] (.75+.21+.245,1.7+.21-.435) -- ++(-.18,-.28);

            \filldraw[line width=1pt,fill=gray!35, dashed] (0,1) circle (.35);
            \draw[line width=1.5pt , ->] (0-1/10+1/20,1) -- ++(-0.15,-.15);
            \draw[line width=1.5pt , ->] (0+1/20,1) -- ++(-0.15,-.15); 
            \draw[line width=1.5pt , ->] (0+1/10+1/20,1) -- ++(-0.15,-.15);

            
            \draw[] (-1.15,0.2) node[anchor = north]{$\bm{u} = \bm{u}_0$ };
            \draw[] (-1.15,0.35-.05) node[anchor = south]{$\Gamma_u$ };

            \draw[] (-1,1.25) node[anchor = south]{$\Omega$ };

            \draw[] (0,1) node[anchor = south]{$\bm{b}$ };

            \draw[] (.75+.21+.125-.2,1.7+.21-.165-.3) node[anchor = east]{$\Gamma_t$ };
            \draw[] (.75+.21+.125+.2,1.7+.21-.165) node[anchor = south]{$\bm{t} = \overline{\bm{t}}$ };

\end{tikzpicture}
	\caption{\ac{Two-dimensional solid that occupies the region $\Omega$ with body force $\vm{b}$, and is subjected to displacement and traction boundary conditions.}}\label{fig:domain}
\end{figure}

Following~\cite{chen:2023:shv}, we introduce the Hellinger--Reissner functional for linear elasticity:
\begin{align*}
    {\Pi}_{\textrm{HR}}[\vm{u},\vsigma] = -\frac{1}{2}\int_{\Omega}{\vsigma: \mathbb{C}^{-1}:\vsigma \, d\vm{x}}  + \int_{\Omega}{\vsigma:\symgrad{\vm{u}}\, d\vm{x}} - \int_{\Omega}{\vm{b}\cdot \vm{u} \, d\vm{x}} - \int_{\Gamma_t}{\Bar{\vm{t}}\cdot \vm{u} \, ds} .
\end{align*}
After taking the first variation of $\Pi_{\textrm{HR}}(\cdot,\cdot)$ and 
requiring it to be stationary, we obtain
the weak statement of the equilibrium equations and strain-displacement relations:
\begin{subequations}\label{eq:weak_equations}
\begin{align}
    &\int_{\Omega}{\vsigma : \symgrad{(\delta\vm{u})} \, d\vm{x}} - \int_{\Omega}{\vm{b} \cdot \delta \vm{u} \, d\vm{x}} - \int_{\Gamma_t}{\Bar{\vm{t}}\cdot \delta\vm{u} \, ds} =0 \ \ \forall 
    \delta \vm{u} \in {\cal V}_u, \label{eq:weak_equilibrium} \\ 
     &\int_{\Omega}{\delta\vsigma : \left(\nabla_s\vm{u}-\mathbb{C}^{-1}:\vsigma\right) \, d\vm{x}} = 0 \ \ \forall \delta \vsigma \in {\cal V}_\sigma, \label{eq:weak_compatibility}
\end{align}
\end{subequations}
where ${\cal V}_u\subset [H^1(\Omega)]^2$ contains vector-valued functions that vanish on $\Gamma_u$ and ${\cal V}_\sigma$ contains symmetric tensor-valued functions in $[L^2(\Omega)]^{2\times 2}_{\textrm{sym}}$. 
\section{ Stress-hybrid virtual element discretization}\label{sec:VEM}
In this section, we introduce the projection operations used in the virtual element method. In the virtual element method, the polygonal basis functions are not explicitly constructed and only known on the boundary of an element; therefore, projection operators are needed to 
\ac{give a polynomial approximation of the displacement and stress fields in the interior of elements}. We first define the energy projection of the displacement field using an energy orthogonality condition \ac{in~\sref{sec:3_1_energy_proj}}. \ac{In~\sref{sec:3_2_stress_proj}}, we define the stress-hybrid projection using the orthogonality condition based on the weak strain-displacement relation obtained from the Hellinger--Reissner principle. Finally, we define the virtual element space and the virtual basis functions over a single element \ac{in~\sref{sec:3_3_vem_space}}.

\begin{figure}[!bht]
    \centering
    \begin{subfigure}{0.32\textwidth}
    \centering
    \includegraphics[width=\textwidth]{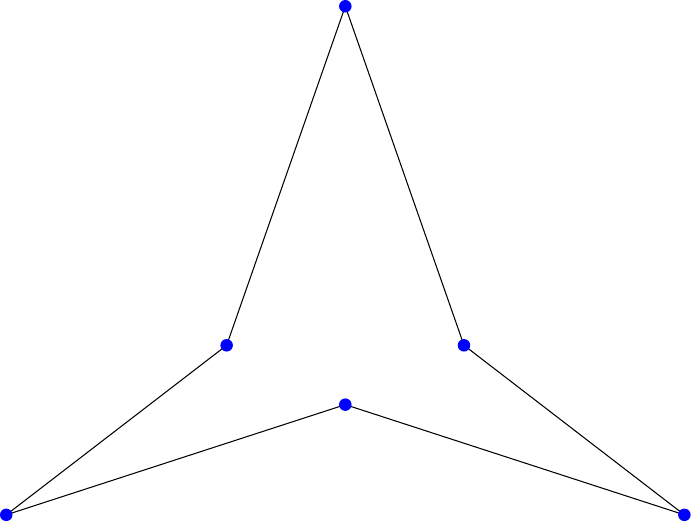}
    \caption{}\label{fig:example_single_elements_a}
    \end{subfigure}
     \hfill
     \begin{subfigure}{0.25\textwidth}
         \centering
         \includegraphics[width=\textwidth]{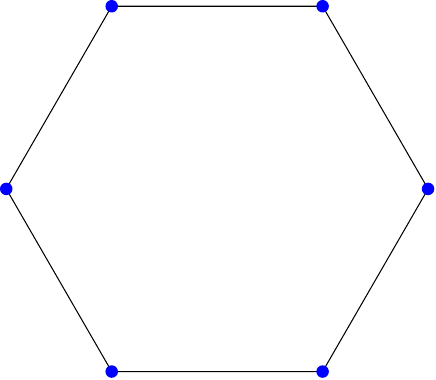}
         \caption{}\label{fig:example_single_elements_b}
     \end{subfigure}
     \hfill
     \begin{subfigure}{0.32\textwidth}
         \centering
         \includegraphics[width=\textwidth]{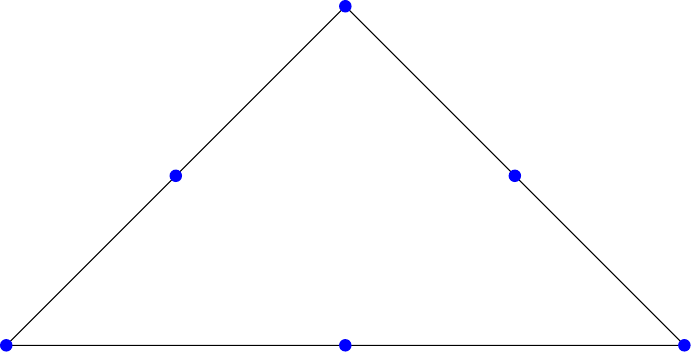}
         \caption{}\label{fig:example_single_elements_c}
     \end{subfigure}
    \caption{Examples of admissible six-noded elements (a) nonconvex element, (b) hexagonal element and (c) six-noded triangular element. }
    \label{fig:example_single_elements}
\end{figure}

Let $\mathcal{T}^h$ be a decomposition of $\Omega$ into six-noded elements (see~\fref{fig:example_single_elements}). For each element $E$, denote the diameter by $h_E$, its centroid by $\vx_E$ and the coordinate of the $i$-th vertex by $\vx_i = (x_i, y_i)$. 
\subsection{Energy projection operator of the displacement field}\label{sec:3_1_energy_proj}
For the displacement field approximation in each element, we choose a basis of first order scaled vector monomials given by:
\begin{subequations}
    \begin{align}
        \widehat{\vm{M}}(E) &= \begin{bmatrix}
    \begin{Bmatrix}
    1 \\ 0 
    \end{Bmatrix}, 
    \begin{Bmatrix}
    0 \\ 1
    \end{Bmatrix}, 
    \begin{Bmatrix}
    -\eta \\ \xi
    \end{Bmatrix}, 
    \begin{Bmatrix}
    \eta \\ \xi
    \end{Bmatrix}, 
    \begin{Bmatrix}
    \xi \\ 0
    \end{Bmatrix}, 
    \begin{Bmatrix}
    0 \\ \eta
    \end{Bmatrix}
        \end{bmatrix}\label{eq:disp_poly_basis},\\
        \intertext{where}
        \xi &= \frac{x-x_E}{h_E}, \quad \eta = \frac{y-y_E}{h_E}. \\
        \intertext{Let $\vm{m}_\alpha$ denote the $\alpha$th column of $\widehat{\vm{M}}(E)$, and define $\wtilde{\vm{M}}$ as the matrix representation given by} 
            \wtilde{\vm{M}} &= \begin{bmatrix}
        1 & 0 & -\eta & \eta & \xi & 0 \\ 0 & 1 & \xi & \xi & 0 & \eta
    \end{bmatrix}.\label{eq:M_basis_matrix}
    \end{align}
\end{subequations}
\ac{Then, for each element $E$}, define the energy projection of the displacement field $\proj: [H^1(E)]^2 \to [\symP_1(E)]^2$ \ac{such that for any $\vm{v} \in [H^1(E)]^2$, $\proj{\vm{v}}$ is a vector polynomial function that satisfies the energy orthogonality condition:}
\begin{equation}
    \int_{E}{\veps(\vm{m}_\alpha):\mathbb{C}:\veps(\vm{v}-\proj{\vm{v}}) \, d\vx} =0 \quad \forall \vm{m}_\alpha \in \widehat{\vm{M}}(E).
\end{equation}
For $\alpha=1,2,3$, $\vm{m}_\alpha$ correspond to rigid-body modes, so we get $\veps(\vm{m}_\alpha)=0$. This results in three trivial equations $0=0$, and therefore we require three additional conditions to define a unique projection. Let $Q_0: [H^1(E)]^2 \times [H^1(E)]^2 \to \mathbb{R}$ be the discrete $L^2$ inner product:
\begin{subequations}
  \begin{align}
    Q_0(\vm{u},\vm{v}) &= \frac{1}{6} \sum_{k=1}^{6}{\vm{u}(\vx_k)\cdot \vm{v}(\vx_k)}, \\
    \intertext{and then impose the additional three constraints on $\proj{\vm{v}}$:}
    Q_0(\vm{m}_{\alpha},\vm{v}-\proj{\vm{v}}) &=\frac{1}{6}\sum_{k=1}^{6}{\vm{m}_{\alpha}}(\vx_k)\cdot (\vm{v}-\proj{\vm{v}})(\vx_k) = 0 \quad (\alpha=1,2,3).\label{eq:energy_proj_constraints}
\end{align}  
\end{subequations}
The condition~\eqref{eq:energy_proj_constraints} enforces the nodal average of $\vm{v}$ and its projection $\proj{\vm{v}}$ for a rigid-body mode to be equal. Now we define the energy projection $\proj{\vm{v}}$ as the unique function that satisfies the following two relations:
\begin{subequations}
    \begin{align}
    &\frac{1}{6}\sum_{k=1}^{6}{\vm{m}_{\alpha}}(\vx_k)\cdot (\vm{v}-\proj{\vm{v}})(\vx_k) =0 \quad (\alpha = 1,2,3), \label{eq:proj_P0_condtion} \\
    &\int_{E}{\veps(\vm{m}_\alpha):\mathbb{C}:\veps(\vm{v}-\proj{\vm{v}}) \, d\vx} =0 \quad (\alpha = 4,5,6) \label{eq:proj_integral_condition}.
    \end{align}
\end{subequations}
For implementation, it is more convenient to write~\eqref{eq:proj_integral_condition} in matrix-vector form \ac{using Voigt notation}. From~\cite{Chen:2023:SFV}, this relation can be written as
\begin{subequations}
    \begin{align}
    \int_{E}{(\vm{S}\proj\vm{v})^T(\vm{C}\vm{S}\vm{m}_\alpha)\, d\vx} &= \int_{E}{(\vm{S}\vm{v})^T(\vm{C\vm{S}\vm{m}_\alpha})\, d\vx},\label{eq:proj_integral_condition_matrix} \\
    \intertext{where $\vm{S}$ is the matrix symmetric gradient operator that is
    given by}
    \vm{S} &=\begin{bmatrix}
        \frac{\partial}{\partial x} & 0 \\ 0 & \frac{\partial}{\partial y} \\ \frac{\partial}{\partial y} & \frac{\partial}{\partial x} \end{bmatrix},
        \\
    \intertext{and $\vm{C}$ is the matrix representation of the material moduli tensor.
    For our study, we assume that the material is homogeneous, isotropic and linearly elastic. We also use the plane strain condition since we are interested in applications when the material is nearly incompressible. Under these assumptions, the matrix $\vm{C}$ is given by} 
    \vm{C} = \frac{E_Y}{(1+\nu)(1-2\nu)} & \begin{bmatrix}
        1-\nu & \nu & 0 \\ \nu & 1-\nu & 0 \\ 0 & 0 & \frac{1-2\nu}{2}
    \end{bmatrix},
    \end{align}
\end{subequations}
where $E_Y$ is the Young's modulus and $\nu$ is the Poisson's ratio of the material.

\subsection{Stress-hybrid projection operator}\label{sec:3_2_stress_proj}
Let $\ell\in \mathbb{Z}_{\geq 0}$ be the \ac{largest degree of the polynomials} used for the stress approximation on an element. Then following~\cite{chen:2023:shv}, define the projection operator $\Projbeta$ for the stress-hybrid formulation \ac{over the space $\symPspace$ of $2\times 2$ polynomial symmetric tensors of degree less than or equal to $\ell$ } by the condition
\begin{subequations}
    \begin{align*}
    &\int_{E}{\symP : \left(\nabla_s\vm{u}-\mathbb{C}^{-1}:\Projbeta\right) \, d\vm{x}} = 0 \quad \forall \symP \in \symPspace, \\
    \intertext{which is rewritten as}
    &\int_{E}{\symP:\mathbb{C}^{-1}:\Projbeta \, d\vx} = \int_{E}{ \symP : \nabla_s\vm{u}\, d\vx} \quad \forall \symP \in \symPspace.
\end{align*}
\end{subequations} 
Now applying the divergence theorem and simplifying, we obtain
\begin{align}
    \int_{E}{\symP:\mathbb{C}^{-1}:\Projbeta \, d\vx }= \int_{\partial E}{\left(\symP\cdot \vm{n}\right)\cdot\vm{u} \,ds} - \int_{E}{\left(\diverge{\symP}\right)\cdot \vm{u} \, d\vx}, \label{eq:proj_tensor_equation}
\end{align}
Let $\overline{\symP}$, $\overline{\Projbeta}$ be the Voigt representation of $\symP$ and 
$\Projbeta$, respectively. Then, \eqref{eq:proj_tensor_equation} can be written as 

\begin{subequations}
\begin{align}
    \int_{E}{\overline{\symP}^T\vm{C}^{-1}\overline{\Projbeta} \, d\vx} &= \int_{\partial E}{\overline{\symP}^T\Nmatrix\vm{u} \, ds} - \int_{E}{ \left(\vm{\partial}\overline{\symP}\right)^T\vm{u} \, d\vx}\label{eq:matrix_vect_proj} \quad \forall \symP \in \symPspace,\\
    \intertext{where $\Nmatrix$ is the matrix representation of the outward normals along the boundary of $E$, and $\vm{\partial}$ is the matrix divergence operator that are given by}
    &\Nmatrix := \begin{bmatrix}
        n_x & 0 \\ 0 & n_y \\ n_y & n_x
    \end{bmatrix}, \quad \vm{\partial} := \begin{bmatrix}
        \frac{\partial}{\partial x} & 0 & \frac{\partial}{\partial y} \\ 0 & \frac{\partial}{\partial y} & \frac{\partial}{\partial x}
    \end{bmatrix}.
\end{align}
\end{subequations}
\subsection{Virtual element space}\label{sec:3_3_vem_space}
In the standard first order virtual element method~\cite{basicprinciple}, the virtual element space consists of harmonic functions that are piecewise affine and continuous along the boundary of an element. On each vertex of the element, there are two displacement degrees of freedom, which allows for the construction of a linear projection operator of the displacement field and a constant projection operator of the stress field. However, this space is not sufficient to compute the projection $\Projbetavect$ for arbitrary $\ell$ without adding additional degrees of freedom; therefore, we follow the approach of~\cite{berrone:2023:los,Chen:2023:SFV} to define the \ac{enhanced virtual element space for each element $E$:}
\begin{equation}\label{eq:vem_space}
    \begin{split}
    \ac{\vm{V}_h(E)} = \biggl\{  & \vm{v}_h \in [H^1(E)]^2 : \Delta \vm{v}_h \in [\mathbb{P}_{\ell-1}(E)]^2, \ \vm{v}_h |_{e} \in [\mathbb{P}_1(e)]^2 \ \forall e \in \partial E, \\ 
    & \vm{v}_h |_{\partial E} \in [C^0(\partial E)]^2, \
    \int_{E}{\vm{v}_h \cdot \vm{p} \ d\vx} = \int_{E}{ \proj{\vm{v}_h} \cdot \vm{p} \ d\vx} \quad \forall \vm{p} \in [\mathbb{P}_{\ell-1}(E)]^2   \Biggr\}, 
    \end{split}
\end{equation}
where $\Delta$ is the vector Laplacian operator and $e$ is an edge of the element. \ac{We note that for a stress basis that uses degree $\ell$ polynomials, the space requires the enhancing property $\int_{E}{\vm{v}_h \cdot \vm{p} \ d\vx} = \int_{E}{ \proj{\vm{v}_h} \cdot \vm{p} \ d\vx}$ to hold for all vectorial polynomials up to degree $\ell-1$.} For each element $E$ we also assign a basis for the local space $\ac{\vm{V}_h(E)}$. Let $\{\phi_i\}$ be the scalar polygonal basis functions in standard VEM~\cite{basicprinciple} that satisfy the 
Kronecker-delta 
property $\phi_i(\vx_j) = \delta_{ij}$. Using the scalar basis, we define the matrix of vector-valued basis functions by
\begin{subequations}
   \begin{align}
        \vm{\varphi} &= \begin{bmatrix}
    \phi_1 & \phi_2 & \phi_3 &\cdots & \phi_6 & 0 & 0 & 0 & \cdots & 0 \\ 
    0 & 0 & 0 &  \cdots &  0& \phi_1 & \phi_2 & \phi_3 & \cdots & \phi_6
    \end{bmatrix} := \begin{bmatrix}
        \vm{\varphi}_1 & \vm{\varphi}_2 & \dots & \vm{\varphi}_{12}
    \end{bmatrix},\\
    \intertext{then any function $\vm{v}_h \in \ac{\vm{V}_h(E)}$ can be represented as:}
    \vm{v}_h(\vx) &= \sum_{i=1}^{12}{\vm{\varphi}_{i}(\vx)v_i} = \vm{\varphi}\vd,  
\end{align}
\end{subequations}
where $v_i$ is the $i$-th degree of freedom of $\vm{v}_h$ \ac{and $\vd$ is the displacement vector with components $v_i$}.
\section{Stress-hybrid implementation}\label{sec:implementation}
With the definition of the projection operators and the virtual element space, we now present the implementation of the stress-hybrid virtual element method. We first construct the energy projection \ac{in~\sref{sec:4_1_energy_implementation}} and discuss different choices of stress basis when constructing the stress-hybrid projection \ac{in~\sref{sec:4_2_stress_basis} and~\sref{sec:4_3_stress_implementation}}. We then build the element stiffness matrix and the element force vector \ac{in~\sref{sec:4_4_stiffness_force}}.  
\subsection{Implementation of the energy projection}\label{sec:4_1_energy_implementation}
We follow the construction of the energy projection from~\cite{Chen:2023:SFV}. Let $\vm{v}_h = \vm{\varphi}_i$ be the $i$-th basis function $(i=1,2\dots 12)$ of $\ac{\vm{V}_h(E)}$, then after substituting into~\eqref{eq:proj_P0_condtion} and~\eqref{eq:proj_integral_condition_matrix}, the energy projection $\proj{\vm{\varphi}_i}$ is given by the two relations:   
 \begin{subequations}\label{eq:EProj}
    \begin{align}
         \frac{1}{6}\sum_{k=1}^{6}{\vm{m}_{\alpha}(\vx_k)\cdot\proj{\vm{\varphi}_i}(\vx_k)} &= \frac{1}{6}\sum_{k=1}^{6}{\vm{m}_{\alpha}}(\vx_k)\cdot \vm{\varphi}_i(\vx_k) \quad (\alpha=1,2,3), \\
         \int_{E}{\left(\bm{S}\proj{\vm{\varphi}_i}\right)^T\left(\vm{C}\vm{S}\vm{m}_{\alpha}\right)\,d\vx} &= \int_{E}{\left(\vm{S}\vm{\varphi}_i\right)^T\left(\vm{C}\vm{S}\vm{m}_{\alpha}\right)\,d \vx} \quad (\alpha=4,5,6).
 \end{align}  
 \end{subequations}
The projection $\proj{\vm{\varphi}_i}$ is an affine vector polynomial, therefore we expand it in terms of the basis in~\eqref{eq:disp_poly_basis}. That is:
\begin{align*}
    \proj{\vm{\varphi}_i} = \sum_{\mu=1}^{6}{(\vm{\Pi}_{*}^{\veps})_{\mu i}\vm{m}_{\mu}} \quad (i=1,2,\dots,12).
\end{align*}
After substituting in~\eqref{eq:EProj} and simplifying, we obtain the energy projection coefficients $\vm{\Pi}_{*}^{\veps}$ as~\cite{Chen:2023:SFV} 
\begin{subequations}
    \begin{align}
    \vm{\Pi}_{*}^{\veps} &= \vm{G}^{-1}\vm{B},\label{eq:Pi_eps_matrix}
    \intertext{where for $\mu=1,2,\dots,6$}
    \vm{G}_{\alpha \mu} &=\begin{cases}
    \frac{1}{6}\sum_{j=1}^{6}{\vm{m}_\mu(\vx_j)\cdot \vm{m}_\alpha(\vx_j)} \quad (\alpha =1,2,3) \\
    \int_{E}{\left(\vm{S}\vm{m}_\mu\right)^T(\vm{C}\vm{S}\vm{m}_{\alpha}) \, d\vx} \quad 
    (\alpha= 4,5,6),
    \end{cases} \\
    \intertext{and for $i=1,2,\dots,12$}
    \vm{B}_{\alpha i} &=
    \begin{cases}
        \frac{1}{6}\sum_{j=1}^{6}{\vm{\vm{\varphi}_i}(\vx_j)\cdot \vm{m}_\alpha(\vx_j)} \quad (\alpha =1,2,3) \\ 
        \int_{E}{\left(\vm{S}\vm{\varphi}_i\right)^T(\vm{C}\vm{S}\vm{m}_{\alpha}) \, d\vx} \quad (\alpha=4,5,6).
    \end{cases}
    \end{align}
\end{subequations}
\ac{By utilizing properties of the virtual element basis, the coefficients of $\vm{B}$ can be exactly computed. See~\cite{Chen:2023:SFV} for a more detailed construction of the matrix $\vm{B}$.}

\subsection{Choice of stress basis for a hybrid formulation}\label{sec:4_2_stress_basis}
It is known that for a stress-based method, using a stress approximation that satisfies the element equilibrium condition results in more accurate stress distributions for homogeneous problems~\cite{Jirousek:1989:GOA,Spilker:1981:ijnme,Spilker:1981:PIH,Wu:1995:OOA}. For plane isotropic elasticity, a convenient set of suitable stresses that satisfy the equilibrium equations with zero body force are derived from the Airy stress functions~\cite{Cen:2011:cmame,Ghosh:1994:VCF,Rezaiee-Pajand:2020:TSB}. A collection of fifteen potential stress fields is given in Voigt representation by~\cite{Cen:2011:cmame}:
\begin{equation}
   \vm{P}=\bigg[ \begin{smallmatrix}
        1 & 0 & 0 & \eta & 0 & \xi & 0 & 0 & 2\xi\eta & -\eta^2 & \xi^2-\eta^2 & \xi(\xi^2-6\eta^2) & \xi^3 & 3\xi^2\eta & \eta(3\xi^2-2\eta^2)\\ 0 & 1 & 0 & 0 & \xi & 0 & \eta & 2\xi\eta & 0 & \xi^2 & \eta^2-\xi^2 & 3\xi\eta^2 & -\xi(2\xi^2-3\eta^2) & -\eta(6\xi^2-\eta^2) & \eta^3 \\ 0 & 0 & 1 & 0 & 0& -\eta & -\xi & -\xi^2 & -\eta^2 & 0 & -2\xi\eta & -\eta(3\xi^2-2\eta^2) & -3\xi^2\eta & \xi(2\xi^2-3\eta^2) & -3\xi\eta^2 
    \end{smallmatrix}\bigg]\label{eq:airy_basis}.
\end{equation}
For a six-noded element, there are twelve displacement degrees of freedom and three rigid-body modes, so a minimum of nine terms are needed for the stress approximation~\cite{Pian:1988:ijnme}. However, it is also known that using an incomplete stress approximation will result in a stiffness matrix that is not rotationally invariant~\cite{Cook:1974:jsd,Pian:1984:ijnme}, or result in elements with inaccurate stress distributions~\cite{Wu:1995:OOA}.
Therefore, we only consider the complete bases with the first eleven terms of~\eqref{eq:airy_basis} and all fifteen terms. We denote any stress basis with $k$ independent terms by $k\beta$.  
\begin{remark}
The set of stress basis in~\eqref{eq:airy_basis} is an extension of the $5\beta$ hybrid-stress basis introduced in~\cite{Pian:1964:DOE} and the $15\beta$ basis introduced in~\cite{Spilker:1981:PIH}. In~\cite{enhanced:VEM}, a similar collection of divergence-free polynomials is used as a basis for the enhanced strain VEM. 
\end{remark}

We also examine a basis similar to the hybrid basis given in~\cite{enhanced:VEM}, which uses a mix of nine uncoupled constant and linear polynomials and four divergence-free quadratic polynomials:
\begin{equation}
\vm{P}^{*}=\begin{bmatrix}
    1 & 0 & 0 & \xi & 0 & 0 & \eta & 0 & 0 & 0 & 2\xi\eta & -\eta^2 & \xi^2-\eta^2 \\
    0 & 1 & 0 & 0& \xi & 0 & 0& \eta & 0 & 2\xi\eta & 0 & \xi^2 & \eta^2-\xi^2 \\ 
    0 & 0 & 1 & 0 & 0 & \xi & 0 & 0 & \eta  & -\xi^2 & -\eta^2 & 0 & -2\xi\eta
\end{bmatrix}\label{eq:hybrid_basis}.    
\end{equation}

\subsection{Implementation of the stress-hybrid projection}\label{sec:4_3_stress_implementation}
We have the discrete relation of~\eqref{eq:matrix_vect_proj} given by:
 \begin{subequations}
      \begin{align}
     \int_{E}{\overline{\symP}^T\vm{C}^{-1}\overline{\Projbeta} \, d\vx} &= \int_{\partial E}{\overline{\symP}^T\Nmatrix\vm{u}_h \, ds} - \int_{E}{ \left(\vm{\partial}\overline{\symP}\right)^T\vm{u}_h \, d\vx}.\label{eq:matrix_vect_proj_discrete_a} \\
     \intertext{The last term in~\eqref{eq:matrix_vect_proj_discrete_a} is not computable for general $\overline{\symP}$; however since the elements of $\vm{\partial}\overline{\symP}$ are polynomials of degree at most $\ell-1$, \ac{we apply the enhancing property of the local virtual element space $\vm{V}_h(E)$ to rewrite as}}
     \int_{E}{\overline{\symP}^T\vm{C}^{-1}\overline{\Projbeta} \, d\vx} &= \int_{\partial E}{\overline{\symP}^T\Nmatrix\vm{u}_h \, ds} - \int_{E}{ \left(\vm{\partial}\overline{\symP}\right)^T\proj{\vm{u}_h} \, d\vx}.\label{eq:matrix_vect_proj_discrete_b}
 \end{align}
 \end{subequations}
Expanding $\vm{u}_h$ in terms of the basis in $\ac{\vm{V}_h(E)}$, we have $\vm{u}_h = \vm{\varphi}\vd$, where $\vd$ is the displacement vector. We also expand $\overline{\Projbeta}$ in terms of $\vm{P}$: $\overline{\Projbeta} = \vm{P}\vm{\beta}$,\ac{ where $\vm{\beta}$ is the vector of stress coefficients}, 
 and since $\symP$ is arbitrary we take $\overline{\symP} = \vm{P}_i \ \ (i=1,2,\dots,15)$. After substituting in~\eqref{eq:matrix_vect_proj_discrete_b} for each $i=1,2,\dots,15$ and simplifying, we obtain the system:
\begin{equation}\label{eq:matrix_vect_proj_system}
     \left(\int_{E}{\vm{P}^T \vm{C}^{-1}\vm{P} \, d\vx}\right)\vm{\beta}  = \left(\int_{\partial E}{\vm{P}^T\Nmatrix\vm{\varphi} \, ds} - \int_{E}{\left(\vm{\partial} \vm{P}\right)^T\proj{\vm{\varphi}} \, d\vx}\right)\vd.
\end{equation}
For any choice of $\vm{P}$ in~\eqref{eq:airy_basis}, we have the divergence-free condition $\vm{\partial} \vm{P}=\vm{0}$, so we obtain
\begin{equation}
      \left(\int_{E}{\vm{P}^T \vm{C}^{-1}\vm{P} \, d\vx}\right)\vm{\beta}  =\left(\int_{\partial E}{\vm{P}^T\Nmatrix\vm{\varphi} \, ds} \right)\vd.
 \end{equation}
 Now define the corresponding matrices $\vm{H}$ and $\vm{L}$ by
 \begin{subequations}\label{eq:H_L_beta}
      \begin{align}
      \vm{H} &= \int_{E}{\vm{P}^T \vm{C}^{-1}\vm{P} \, d\vx}, \quad
     \vm{L} = \int_{\partial E}{\vm{P}^T\Nmatrix\vm{\varphi} \, ds}. \label{eq:HandL_matrix}\\
     \intertext{Since $\vm{P}$ has linearly independent columns and $\vm{C}^{-1}$ is symmetric positive-definite, the matrix $\vm{H}$ is invertible. Then the stress coefficients are given by }
     \vm{\beta} &= \vm{H}^{-1}\vm{L}\vd := \vm{\Pi}_{\beta}\vd , \label{eq:beta_displacement}
 \end{align}
 \end{subequations}
 where $\vm{\Pi}_{\beta}$ is the matrix representation
 of the stress-hybrid projection operator with respect to
 the symmetric tensor polynomial basis $\vm{P}$.

\begin{remark}
    For a nondivergence-free basis, the use of~\eqref{eq:matrix_vect_proj_discrete_b} can lead to an overly stiff element. An alternate approach used in~\cite{enhanced:VEM,Lamperti:2023:cmech}, is to introduce additional internal moment degrees of freedom to compute the last integral in~\eqref{eq:matrix_vect_proj_discrete_a}. These additional DOFs result in better performing elements but require an additional static condensation on the element stiffness matrix. Another approach is to follow the idea of using composite elements~\cite{Guo:2000:TCF} to compute the stress projection operator. For each element, we construct a sub-triangulation and assume the displacement field (displacement projection) is affine on each subtriangle. This piecewise displacement field can then be used to compute the integral in~\eqref{eq:matrix_vect_proj_discrete_a}. We tested a virtual element formulation based on composite elements and found that the
    resulting elements were more flexible but still suffered from volumetric locking.
\end{remark} 

\subsection{Element stiffness and forcing}\label{sec:4_4_stiffness_force}
Following~\cite{chen:2023:shv}, we define the discrete system:
\begin{subequations}
    \begin{align*}
        &a_h^E(\vm{u}_h,\delta\vm{u}_h) = \ell^E_h(\delta \vm{u}_h),\\
        \intertext{where }
        a_h^E(\vm{u}_h,\delta\vm{u}_h)&:=\int_{E}{\overline{\Projbeta(\delta \vm{u}_h)}^T \vm{C}^{-1} \overline{\Projbeta( \vm{u}_h)} \, d\vx}, \\ 
        \ell^E_h(\delta \vm{u}_h) &:= \int_{E}{(\delta \vm{u}_h)^T\vm{b} \, d\vx} + \int_{\Gamma_t\cap \partial E}{(\delta\vm{u}_h)^T\ac{\Bar{\vm{t}}} \, ds}.
    \end{align*}
\end{subequations}
After expanding $\overline{\Projbeta}$ and simplifying, we construct the element stiffness matrix
\begin{align}\label{eq:SH-VEM_stiffness}
    \vm{K}_E = (\vm{\Pi}_\beta)^T\left(\int_{E}{ \vm{P}^T\vm{C}^{-1}\vm{P}  \, d\vx}\right)\vm{\Pi}_\beta = \vm{\Pi}_\beta^T \vm{H} \vm{\Pi}_\beta.
\end{align}
\ac{We note that unlike in the standard virtual element method, the element stiffness matrix in the stress-hybrid formulation does not have an additional stabilization term. Instead, the stability of the element stiffness is enforced by a suitable choice of stress basis functions.}

Similarly, for every element $E$, the element force vector is given by
 \begin{align}\label{eq:force_vect}
    \vm{f}_{E} :=  \int_{E}{\vm{\varphi}^T\vm{b}  \, d\vm{x}} + \int_{\Gamma_t \cap \partial E}{\vm{\varphi}^T\Bar{\vm{t}} \, ds}.
 \end{align}

\section{Equilibrium penalty stress-hybrid method}\label{sec:penalty}
\ac{A difficulty with the stress-hybrid approach is the need to construct a suitable basis for the stress field. In the previous section, we utilized stress expansions derived from the Airy stress functions; however, it is difficult to find these stress functions for anistropic and nonlinear materials or for three-dimensional problems.}
From our tests, we also found that using only \ac{equilibrated} basis functions requires many higher order functions and more terms than the optimal basis to retain stability, which results in stiffer solutions and greater integration costs. We examine an approach by~\cite{Wu:1995:OOA}, which uses a penalty term to weakly enforce the element equilibrium equations. The penalty formulation was originally shown to mitigate shear locking in the Pian-Sumihara element and to improve performance in axisymmetric problems~\cite{Wu:1995:OOA,Wu:1996:PEA}. A similar method was developed in~\cite{Sze:2000:OIF} to improve the bending solution of the Pian-Sumihara element on distorted meshes. Later, the penalty equilibrium formulation was applied to fracture mechanics~\cite{Xiao:1999:APE} and extended to the Hu--Washizu variational principle~\cite{Cao:2002:3BE,Cao:2003:HAB}. The addition of an equilibrium penalty term is also used in~\cite{Blacker:1994:SPR,Wiberg:1993:PRB} to construct superconvergent stress recovery methods.  

We start with a modified Hellinger--Reissner functional with an equilibrium penalty term:
\begin{equation}
    \begin{split}
        {\Pi}_{\textrm{HR*}}[\vm{u},\vsigma] = -\frac{1}{2}\int_{\Omega}{\vsigma: \mathbb{C}^{-1}:\vsigma \, d\vm{x}}  + \int_{\Omega}{\vsigma:\symgrad{\vm{u}}\, d\vm{x}} - \int_{\Omega}{\vm{b}\cdot \vm{u} \, d\vm{x}} - \int_{\Gamma_t}{\Bar{\vm{t}}\cdot \vm{u} \, ds} \\ - \frac{\alpha}{2} \int_{\Omega}{({\diverge{\vsigma}}+\vm{b}) \cdot (\diverge{\vsigma}+\vm{b}) \, d\vx},
    \end{split}
\end{equation}
where $\alpha > 0$ is a penalty parameter. After taking the first variation, the stationary condition results in the weak equilibrium equations and modified strain-displacement relations:
\begin{subequations}
    \begin{align}
    &\int_{\Omega}{\vsigma : \symgrad{(\delta\vm{u})} \, d\vm{x}} - \int_{\Omega}{\vm{b} \cdot \delta \vm{u} \, d\vm{x}} - \int_{\Gamma_t}{\Bar{\vm{t}}\cdot \delta\vm{u} \, ds} =0 \ \ \forall 
    \delta \vm{u} \in {\cal V}_u,\label{eq:mod_weak_equilibrium} \\ 
     &\int_{\Omega}{\delta\vsigma : \left(\nabla_s\vm{u}-\mathbb{C}^{-1}:\vsigma\right) \, d\vm{x}}-\alpha\int_{\Omega}{(\diverge{\delta \vsigma})\cdot (\diverge{\vsigma}+b)} = 0 \ \ \forall \delta \vsigma \in {\cal V}_\sigma.\label{eq:mod_strain_disp_relation}
\end{align}
\end{subequations}
\ac{Note that the penalty parameter $\alpha$ only appears explicitly in the modified strain-displacement relation~\eqref{eq:mod_strain_disp_relation}, which is used to construct the stress projection operator. }
\subsection{Penalized stress-hybrid projection}
Using~\eqref{eq:mod_strain_disp_relation} over each element $E$, we define a stress projection operator by the condition:
\begin{subequations}
\begin{align*}
     &\int_{E}{\symP : \left(\nabla_s\vm{u}-\mathbb{C}^{-1}:\Projbeta\right) \, d\vm{x}}-\alpha\int_{E}{(\diverge{\symP}) \cdot (\diverge{\Projbeta}+\vm{b}) \, d\vx} = 0, \\
     \intertext{which can be rewritten as}
     &\int_{E}{\symP :\mathbb{C}^{-1} :\Projbeta \, d\vx} +\alpha\int_{E}{(\diverge{\symP})\cdot (\diverge{\Projbeta}) \,d\vx} = \int_{E}{ \symP : \nabla_s\vm{u}\, d\vx} - \alpha\int_{E}{(\diverge{\symP})\cdot \vm{b} \, d\vx}.
\end{align*}
\end{subequations}
After applying the divergence theorem, we get
\begin{equation}\label{eq:penalty_proj_tensor_equation}
     \begin{split}
      \int_{E}{\symP :\mathbb{C}^{-1} :\Projbeta \, d\vx} \ + \ &  \alpha\int_{E}{(\diverge{\symP})\cdot (\diverge{\Projbeta}) \,d\vx} \\
      &= \int_{\partial E}{(\symP \cdot \vm{n}) \cdot \vm{u} \, ds } 
       - \int_{E}{(\diverge{\symP})\cdot \vm{u} \, d\vx} 
      - \alpha \int_{E}{(\diverge{\symP})\cdot \vm{b} \, d\vx}.
      \end{split}
\end{equation}
Using Voigt notation, we rewrite~\eqref{eq:penalty_proj_tensor_equation} in terms of matrix-vector operations
\begin{equation}\label{eq:penalty_proj_matrix_equation}
      \begin{split}
      \int_{E}{\overline{\symP}^T\vm{C}^{-1}\overline{\Projbeta} \, d\vx} +\alpha\int_{E}{(\vm{\partial}{\overline{\symP}})^T \vm{\partial}{\overline{\Projbeta}} \,d\vx} &= \int_{\partial E}{\overline{\symP}^T \Nmatrix\vm{u} \, ds } - \int_{E}{(\vm{\partial}{\overline{\symP}})^T \vm{u} \, d\vx} \\ &-\alpha \int_{E}{(\vm{\partial}{\overline{\symP}})^T \vm{b} \, d\vx}.
      \end{split}
\end{equation}

\subsection{Implementation of penalized stress-hybrid projection}
From~\eqref{eq:penalty_proj_matrix_equation}, we have the discrete expression:
\begin{subequations}
    \begin{align}
      \begin{split}
      \int_{E}{\overline{\symP}^T\vm{C}^{-1}\overline{\Projbeta} \, d\vx} +\alpha\int_{E}{(\vm{\partial}{\overline{\symP}})^T \vm{\partial}{\overline{\Projbeta}} \,d\vx} 
      = & \int_{\partial E}{\overline{\symP}^T \Nmatrix\vm{u}_h \, ds } - \int_{E}{(\vm{\partial}{\overline{\symP}})^T \vm{u}_h \, d\vx} \\ &-\alpha \int_{E}{(\vm{\partial}{\overline{\symP}})^T \vm{b} \, d\vx}.
      \end{split} \\
      \intertext{The second term on the right-hand side is not computable from the element DOF's; however, by applying the \ac{enhancing property} of the virtual element space~\eqref{eq:vem_space}, we rewrite the relation in a computable form}
      \begin{split}
          \int_{E}{\overline{\symP}^T\vm{C}^{-1}\overline{\Projbeta} \, d\vx} +\alpha\int_{E}{(\vm{\partial}{\overline{\symP}})^T \vm{\partial}{\overline{\Projbeta}} \,d\vx} = & \int_{\partial E}{\overline{\symP}^T \Nmatrix\vm{u}_h \, ds } - \int_{E}{(\vm{\partial}{\overline{\symP}})^T \proj\vm{u}_h \, d\vx} \\ &-\alpha \int_{E}{(\vm{\partial}{\overline{\symP}})^T \vm{b} \, d\vx}.
      \end{split}
\end{align}
\end{subequations}
Letting $\overline{\symP} = \vm{P}, \overline{\Projbeta} = \vm{P}\vm{\beta}$ and $\vm{u}_h = \vm{\varphi}\vm{d}$, we obtain the system of equations:
\begin{equation}
    \begin{split}
      \biggl(\int_{E}{\vm{P}^T\vm{C}^{-1}\vm{P} \, d\vx}  &+ \alpha\int_{E}{(\vm{\partial}{\vm{P})^T \vm{\partial}\vm{P} \,d\vx}}\biggr)\vm{\beta} \\
       &= \biggl(\int_{\partial E}{\vm{P}^T \Nmatrix\vm{\varphi} \, ds } - \int_{E}{(\vm{\partial}{\vm{P})^T \proj\vm{\varphi}} \, d\vx}\biggr)\vm{d} - \alpha \int_{E}{(\vm{\partial}\vm{P})^T \vm{b} \, d\vx}.
      \end{split}
\end{equation}
Define the corresponding matrices by  
\begin{subequations}\label{eq:penalty_H_L_matrices}
    \begin{align}
        \vm{H} &= \int_{E}{\vm{P}^T\vm{C}^{-1}\vm{P} \, d\vx} , \quad \vm{H}_p = \int_{E}{(\vm{\partial}{\vm{P})^T \vm{\partial}\vm{P} \,d\vx}}, \\ 
        \vm{L} &= \int_{\partial E}{\vm{P}^T \Nmatrix\vm{\varphi} \, ds }- \int_{E}{(\vm{\partial}{\vm{P})^T \proj\vm{\varphi}} \, d\vx} , \quad \vm{L}_p =\int_{E}{(\vm{\partial}\vm{P})^T \vm{b} \, d\vx}.\\
    \intertext{Then the stress coefficients are given as }
    \vm{\beta} &= (\vm{H}+\alpha\vm{H}_p)^{-1}(\vm{L}\vm{d}-\alpha\vm{L}_p).
    \end{align}
\end{subequations}

\subsection{Element stiffness matrix and element force vector}
To construct the element stiffness matrix, we first consider the discrete equilibrium equations based on~\eqref{eq:mod_weak_equilibrium}
\begin{equation*}
    \int_{E}{\overline{\Projbeta(\vm{u}_h)}^T \symgrad{(\delta\vm{u}_h)} \, d\vm{x}} = \int_{E}{\vm{b}^T \delta \vm{u}_h \, d\vm{x}} + \int_{\Gamma_t}{\Bar{\vm{t}}^T \delta\vm{u}_h \, ds}.
\end{equation*}
On applying~\eqref{eq:mod_strain_disp_relation} and simplifying, we rewrite the first integral as
\begin{equation}
    \begin{split}
    \int_{E}{\overline{\Projbeta(\vm{u}_h)}^T \symgrad{(\delta\vm{u}_h)} \, d\vm{x}} = & \int_{E}{\overline{\Projbeta(\vm{u}_h)}^T\vm{C}^{-1}\overline{\Projbeta(\delta \vm{u}_h)} \, d\vx} \\
     &+ \alpha \int_{E}{(\vm{\partial}\overline{\Projbeta(\vm{u}_h)})^T(\vm{\partial}\overline{\Projbeta(\delta\vm{u}_h)}+\vm{b}) \, d\vx}. 
    \end{split}
\end{equation}
Now we have the equation
\begin{equation}
    \begin{split}
        \int_{E}{\overline{\Projbeta(\vm{u}_h)}^T\vm{C}^{-1}\overline{\Projbeta(\delta \vm{u}_h)} \, d\vx} \ + & \ \alpha \int_{E}{(\vm{\partial}\overline{\Projbeta(\vm{u}_h)})^T(\vm{\partial}\overline{\Projbeta(\delta\vm{u}_h)}+\vm{b}) \, d\vx} \\
        &= \int_{E}{\vm{b}^T \delta \vm{u}_h \, d\vm{x}}  + \int_{\Gamma_t \cap \partial E}{\Bar{\vm{t}}^T \delta\vm{u}_h \, ds}.
    \end{split}
\end{equation}
After expanding $\overline{\Projbeta} = \vm{P}\vm{\beta}$ and simplifying, we construct the element stiffness matrix 
\begin{subequations}
    \begin{align}
        \vm{K}_E &= \vm{L}^T\bigl(\vm{H}+\alpha\vm{H}_p\bigr)^{-1}\vm{L},
        \intertext{and the element force vector is given by}
        \vm{f}_E &= \int_{E}{\vm{\varphi}^T\vm{b} \, d\vx} + \int_{\Gamma_t \cap \partial E}{\vm{\varphi}^T\Bar{\vm{t}} \, ds} + \alpha \vm{L}^T\bigl(\vm{H}+\alpha\vm{H}_p\bigr)^{-1}\vm{L}_p,
    \end{align}
\end{subequations}
where $\vm{H}, \vm{H}_p, \vm{L}, \vm{L}_p$ are given in~\eqref{eq:penalty_H_L_matrices}.
\subsection{Choice of penalty parameter and stress basis functions}
\ac{Similar to constructing a suitable stabilization term in the standard virtual element method, one drawback of the penalty equilibrium method is the need for a properly designed penalty parameter $\alpha$.} In~\cite{Wu:1995:OOA}, $\alpha = \frac{\kappa}{E_Y}$, where $E_Y$ is the Young's modulus and $\kappa>10^3$ is a large dimensionless number, is suggested. However, we found that this can result in the loss of stability for problems with highly distorted meshes. It is suggested in~\cite{Cao:2003:HAB}, that the penalty term should be scaled by a term that depends on the geometry in order to attain consistent units. In particular, $\alpha$ is chosen to be of the form $\alpha = \frac{\kappa\ell_0^2}{E_Y}$, where $\ell_0$ is a characteristic length that is dependent on the element geometry. In our tests, we let $\ell_0$ be the minimum length from the element centroid to the nodes and $\kappa=10^4$. From our numerical experiments on benchmark problems, we found that using the penalty parameter $\alpha = \frac{\kappa\ell_0^2}{E_Y}$ resulted in higher accuracy and superconvergence on sufficiently refined uniform 
meshes. 
However, in the case when $E_Y$ is small ($\alpha$ too large), low energy modes appear in the element stiffness matrix. From~\cite{Long:2006:ANS}, it is suggested that a reasonable value of $\frac{\kappa}{E_Y}$ is between $10$ and $100$; therefore, we fix an arbitrary upper bound of $10$ in our tests. That is, we use 
\begin{equation}\label{eq:penalty_param}
    \alpha = \min\Bigl\{10,\frac{\kappa}{E_Y} \Bigr\}\ell_0^2.
\end{equation} 
For the penalty equilibrium stress-hybrid method, we do not require the \ac{equilibrated} stress basis functions given in~\eqref{eq:airy_basis}. Instead, we seek the smallest number of stress basis functions that still retains stability and is not overly stiff in bending. For a six-noded element, a minimum of nine terms are needed; but it was found in~\cite{Chen:2023:SFV} that the $9\beta$ complete linear basis is not sufficient for stability. We choose the $12\beta$ expansion with complete bilinear polynomials given by
\begin{equation}
    \vm{P} = \begin{bmatrix}\label{eq:penalty_basis}
        1 & 0 & 0 & \xi & 0 & 0 & \eta & 0 & 0 & \xi\eta & 0 & 0 \\
        0 & 1 & 0 & 0 & \xi & 0 & 0 & \eta & 0 & 0 & \xi\eta & 0 \\
        0 & 0 & 1 & 0 & 0 & \xi & 0 & 0 & \eta & 0 & 0 & \xi\eta
    \end{bmatrix}.
\end{equation}

\section{Numerical results}\label{sec:results}
We present a series of numerical examples in linear elasticity under plane strain conditions on a variety of meshes. The unstructured triangular meshes used in these examples are generated using DistMesh~\cite{Persson:2004:ASM}.
We compute the errors of the displacement in the $L^2$ norm and energy seminorm, and the $L^2$ error of the hydrostatic stress (denoted by $\tilde{p}$). The convergence rates of CT FEM, B-bar VEM, SH-VEM, and PSH-VEM are computed using the following discrete error measures: 
 \begin{subequations}\label{eq:error_norms}
\begin{align}
\|\vm{u}-\vm{u}_h\|_{\vm{L}^2(\Omega)} &=\sqrt{\sum_{E}{\int_{E}{|\vm{u}-\proj{\vm{u}_h}}|^2 \,d\vm{x}}},  \\ 
\|\tilde{p}-\tilde{p}_h\|_{\vm{L}^2(\Omega)} &=\sqrt{\sum_{E}{\int_{E}{|\tilde{p}-\tilde{p}_h|^2 \,d\vm{x}}}},\\
\|\vm{u} - \vm{u}_h\|_{a} &= \sqrt{\sum_{E}{\int_{E}{(\overline{\vsigma}-\Projbetavect)^T\,\vm{C}^{-1}(\overline{\vsigma}-\Projbetavect)\, d\vm{x}}}}.
\end{align}
\end{subequations}
We use the scaled boundary cubature (SBC) method~\cite{chin:2021:cmame} to compute the \ac{coefficients of the} matrices $\vm{H}$ in~\eqref{eq:HandL_matrix}, $\vm{H}_p$, $\vm{L}$ and $\vm{L}_p$ in~\eqref{eq:penalty_H_L_matrices}, and the integrals appearing in~\eqref{eq:error_norms}. The \ac{coefficients of the} matrices $\vm{H}$, $\vm{H}_p$ and $\vm{L}$ are integrals of polynomial functions and are exactly computed by the SBC method.
The integrals in~\eqref{eq:error_norms} are in general not integrals of polynomials but they can be computed to arbitrary accuracy with the SBC method.

\subsection{Eigenvalue analysis}\label{subsec:eigenanalysis}
\begin{figure}[!bht]
    \centering
    \begin{subfigure}{0.32\textwidth}
    \centering
    \includegraphics[width=\textwidth]{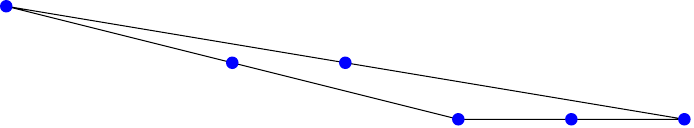}
    \caption{}
    \end{subfigure}
     \hfill
     \begin{subfigure}{0.32\textwidth}
         \centering
         \includegraphics[width=\textwidth]{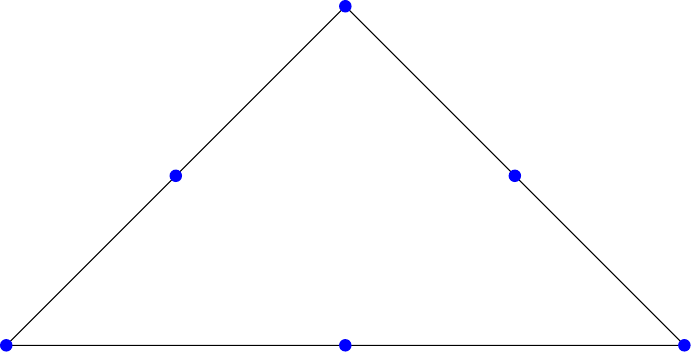}
         \caption{}
     \end{subfigure}
     \hfill
     \begin{subfigure}{0.32\textwidth}
         \centering
         \includegraphics[width=\textwidth]{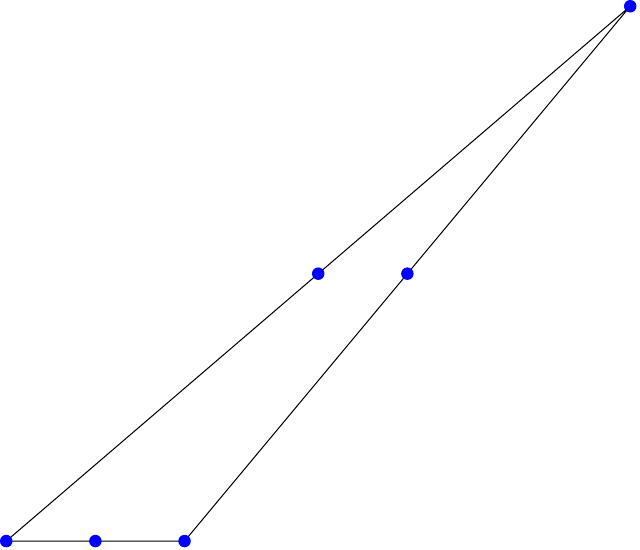}
         \caption{}
     \end{subfigure}
    \caption{(a)-(c) Sequence of six-noded triangular elements with vertices at $\{(-1,0), \, (1,0), \, (\alpha,\beta)\}$ and nodes placed at the midpoint of each edge.}
    \label{fig:eig_single_elements}
\end{figure}

\begin{figure}[!bht]
     \centering
     \begin{subfigure}{.48\textwidth}
         \centering
         \includegraphics[width=\textwidth]{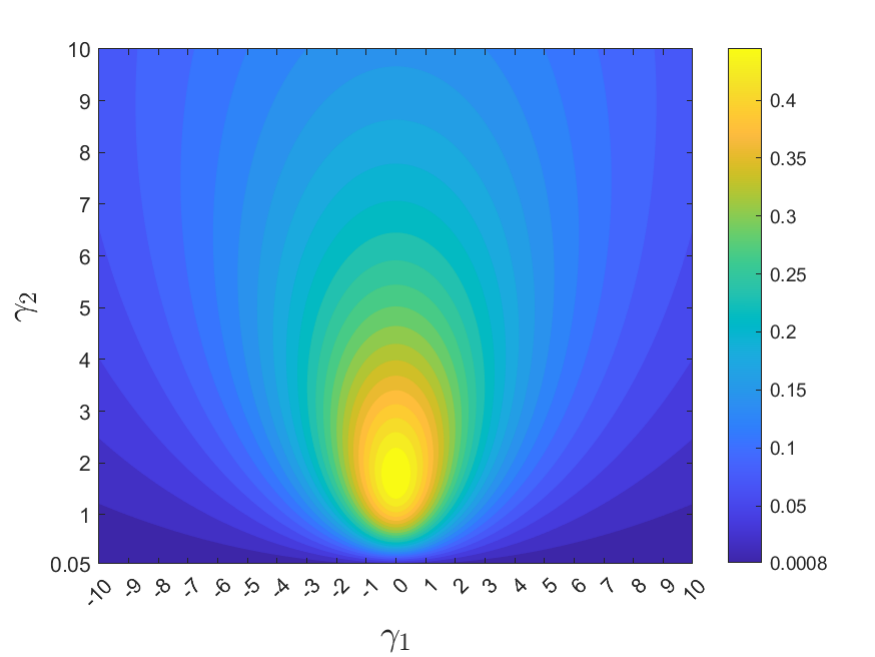}
         \caption{}\label{fig:eig_analysis_1_a}
     \end{subfigure}
     \hfill
     \begin{subfigure}{.48\textwidth}
         \centering
         \includegraphics[width=\textwidth]{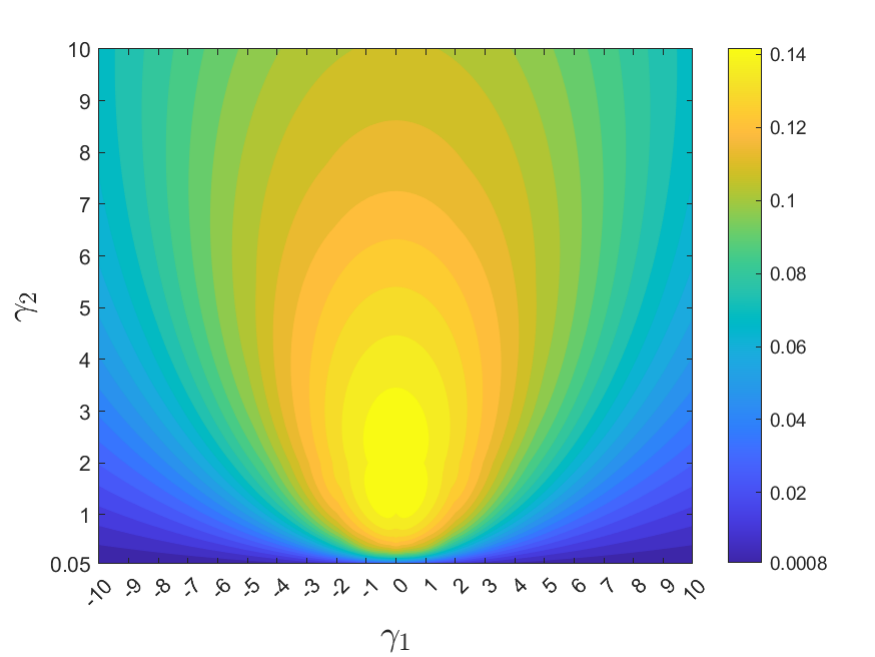}
         \caption{}\label{fig:eig_analysis_1_b}
     \end{subfigure}
     \vfill
     \begin{subfigure}{.48\textwidth}
         \centering
         \includegraphics[width=\textwidth]{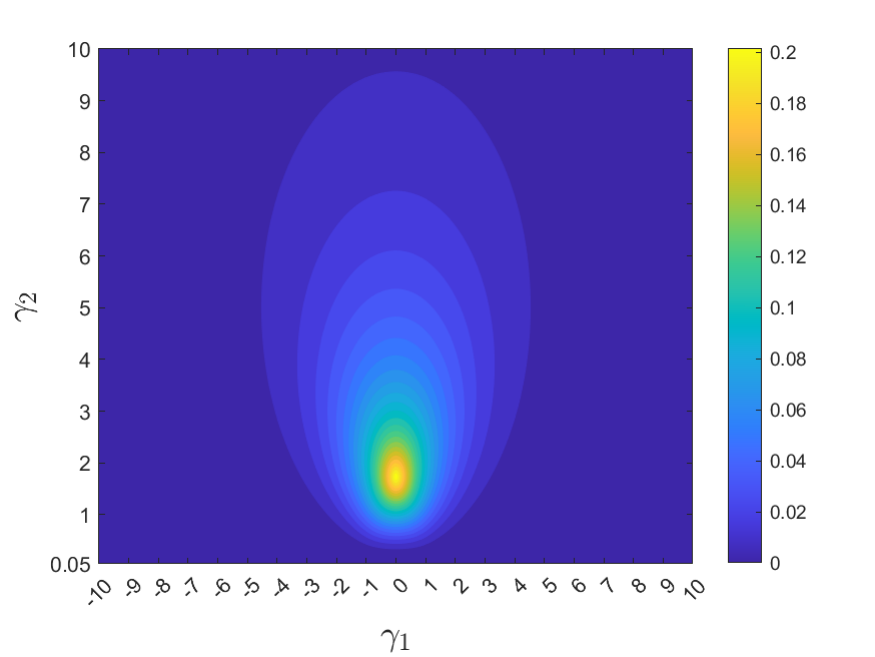}
         \caption{}\label{fig:eig_analysis_1_c}
     \end{subfigure}
      \hfill
     \begin{subfigure}{.48\textwidth}
         \centering
         \includegraphics[width=\textwidth]{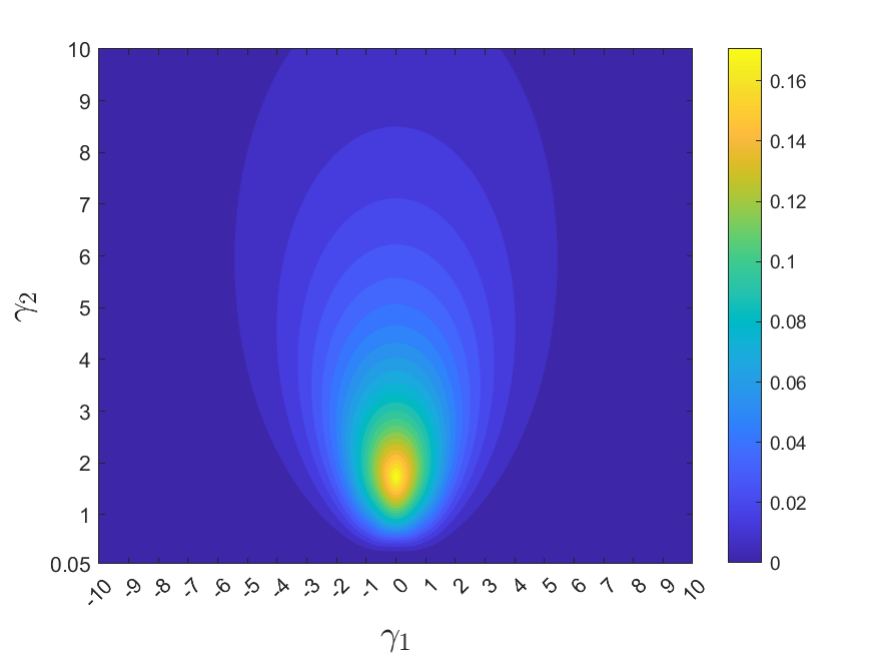}
         \caption{}\label{fig:eig_analysis_1_d}
     \end{subfigure}
     \vfill
     \begin{subfigure}{.48\textwidth}
         \centering
         \includegraphics[width=\textwidth]{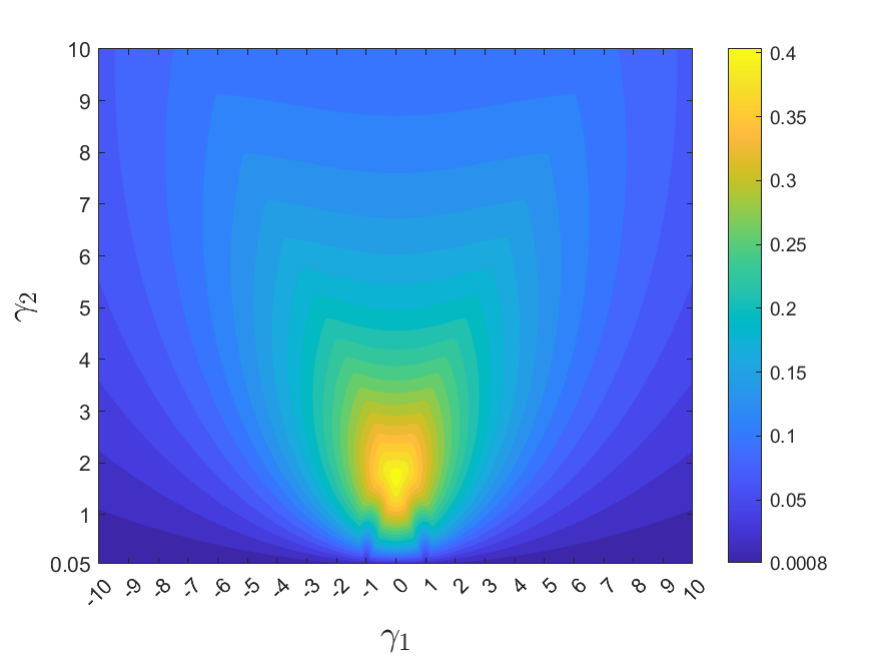}
         \caption{}\label{fig:eig_analysis_1_e}
     \end{subfigure}
     \hfill
     \begin{subfigure}{.48\textwidth}
         \centering
         \includegraphics[width=\textwidth]{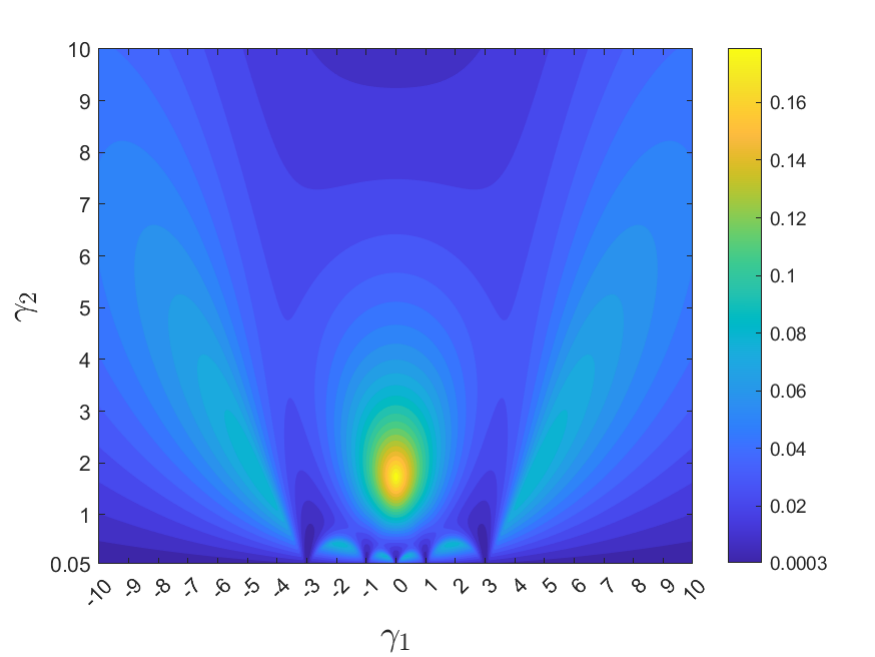}
         \caption{}\label{fig:eig_analysis_1_f}
     \end{subfigure}
        \caption{Contour plots of the fourth-lowest eigenvalue as a function of $(\gamma_1,\gamma_2)$ for (a) standard VEM, (b) B-bar VEM, (c) CT FEM, (d) $11\beta$ SH-VEM, (e) $13\beta$ SH-VEM, and (f) $15\beta$ SH-VEM.}
        \label{fig:eig_analysis_1}
\end{figure}

We examine the eigenvalues of the SH-VEM for general triangular elements to determine the stability of the method. The material has Young's modulus $E_Y = 1$ and Poisson's ratio $\nu=0.49995$. We assess the eigenvalues of the standard VEM~\cite{elasticdaveiga}, B-bar VEM~\cite{Park:2020:meccanica}, composite triangle FEM~\cite{Guo:2000:TCF} and the three formulations: a $11\beta$ and $15\beta$ that is based on the Airy stress function basis given in~\eqref{eq:airy_basis}, and a $13\beta$ hybrid formulation given in~\eqref{eq:hybrid_basis}. For a stable method, the element stiffness matrix should have three zero eigenvalues that correspond to the zero-energy modes and the next smallest eigenvalue should be positive and bounded away from zero. For this test, we construct six-noded triangular elements with vertices at $\{(-1,0), \, (1,0), \, (\gamma_1,\gamma_2)\}$, where $\gamma_1 \in [-10,10]$, $\gamma_2 \in [.05,10]$, and then nodes are placed at the midpoints of each edge (see~\fref{fig:eig_single_elements}). For each combination of $\gamma_1$ and $\gamma_2$, we compute the first non-rigid body eigenvalue, which corresponds to the fourth smallest eigenvalue of the element stiffness matrix. The contour plots of the first non-rigid body eigenvalue are given in Figure~\ref{fig:eig_analysis_1}. \ac{In~\fref{fig:eig_analysis_1_d}, the plot shows that the minimum value of the eigenvalue for the $11\beta$ formulation continues to decrease as the element becomes distorted from an equilateral triangle (when $(\gamma_1,\gamma_2)$ is further away from $(0,\sqrt{3})$). Eventually, the value of the eigenvalue approaches zero and develops into a spurious mode as the elements become highly distorted. In~\fref{fig:eig_analysis_1_e}, the eigenvalue in the $13\beta$ formulation has a similar range to the standard VEM (see~\fref{fig:eig_analysis_1_a}) and further numerical tests show that the $13\beta$ formulation locks in the incompressible limit. Therefore, both the $11\beta$ Airy stress and $13\beta$ hybrid formulations are not considered in the remaining examples.
The $15\beta$ SH-VEM and B-bar VEM produced similar ranges for their eigenvalues and both do not produce any spurious zero-energy modes for the given values of $(\gamma_1,\gamma_2)$ (see Figures~\ref{fig:eig_analysis_1_b} and~\ref{fig:eig_analysis_1_f}). In all the examples that follow, $15\beta$ SH-VEM is denoted as SH-VEM.} 

\subsection{Eigenvalue analysis of the penalty formulation}
\begin{figure}[!bht]
     \centering
     \begin{subfigure}{0.48\textwidth}
         \centering
         \includegraphics[width=\textwidth]{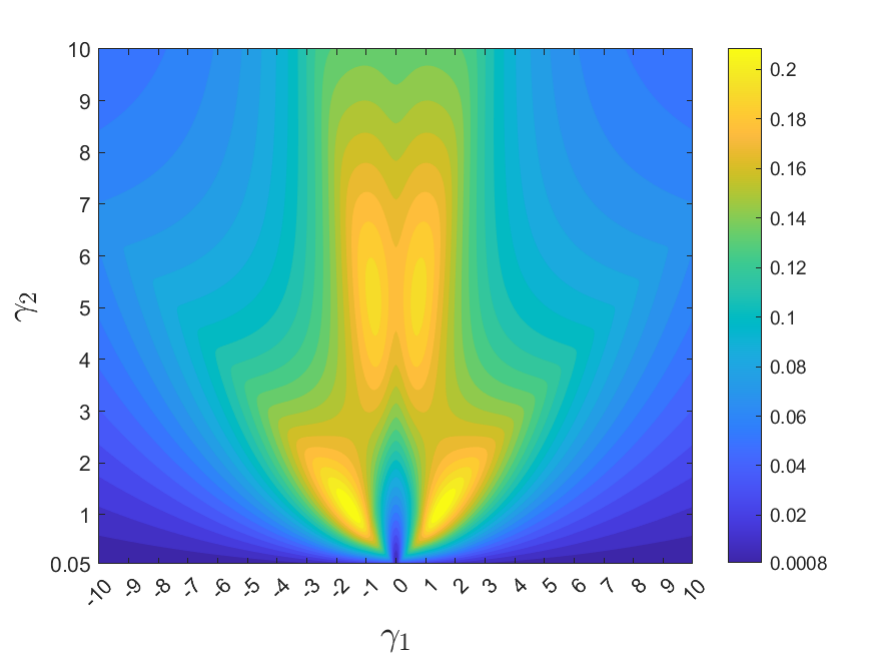}
         \caption{}\label{fig:eigen_analysis_penalty_a}
     \end{subfigure}
     \hfill
     \begin{subfigure}{0.48\textwidth}
         \centering
         \includegraphics[width=\textwidth]{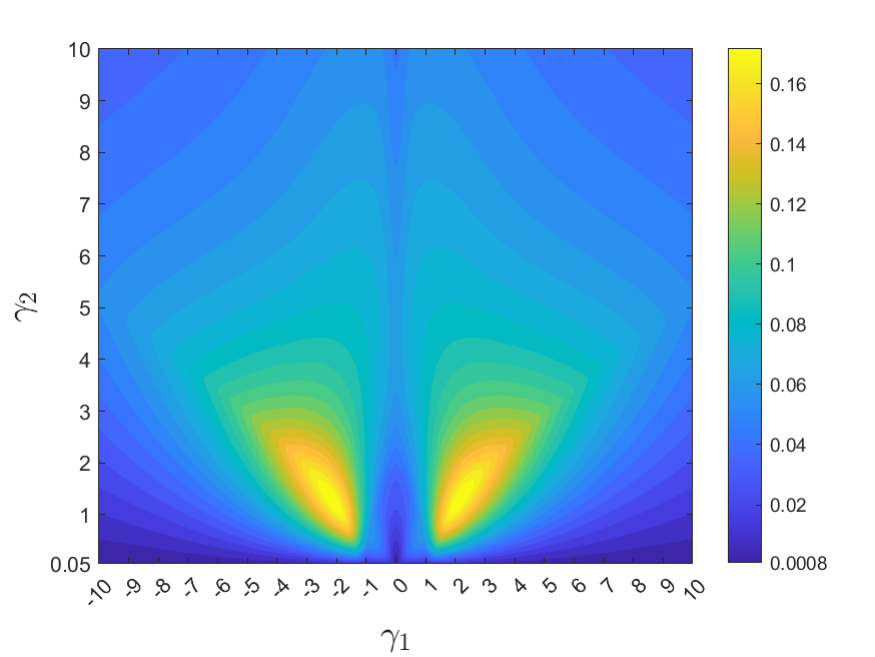}
         \caption{}\label{fig:eigen_analysis_penalty_b}
     \end{subfigure}
     \vfill
     \begin{subfigure}{0.48\textwidth}
         \centering
         \includegraphics[width=\textwidth]{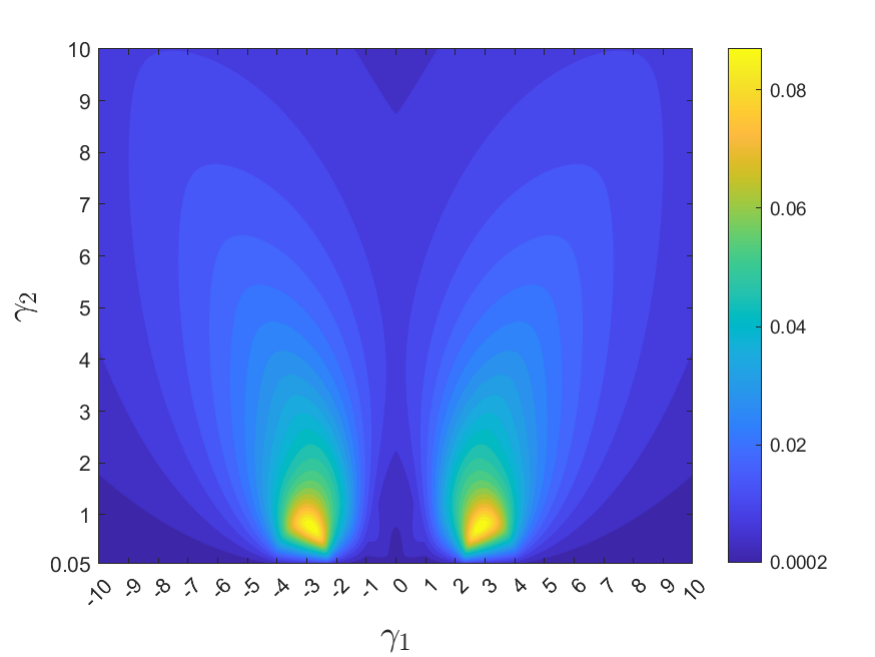}
         \caption{}\label{fig:eigen_analysis_penalty_c}
     \end{subfigure}
     \hfill
     \begin{subfigure}{0.48\textwidth}
         \centering
         \includegraphics[width=\textwidth]{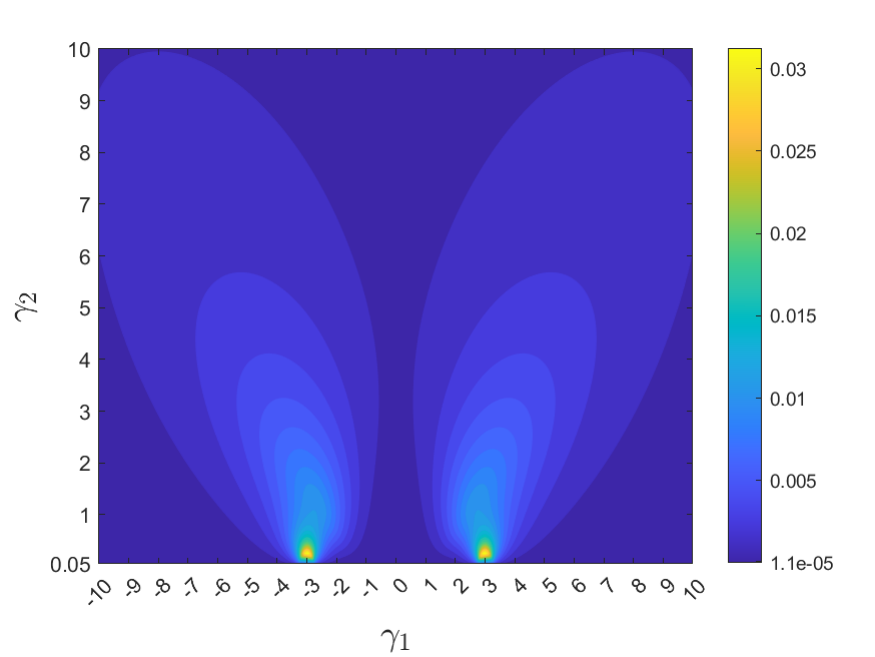}
         \caption{}\label{fig:eigen_analysis_penalty_d}
     \end{subfigure}
        \caption{Contour plots of the fourth smallest eigenvalue as a function of $(\gamma_1,\gamma_2)$ using the penalty parameters (a) $\alpha = \frac{\ell_0^2}{10} $, (b) $\alpha=\ell_0^2 $, (c) $\alpha =10\ell_0^2$, and (d) $\alpha=100\ell_0^2$.}
        \label{fig:eigen_analysis_penalty}
\end{figure}

We test the stability of the PSH-VEM with different values of the penalty parameter $\alpha$. The material has Young's modulus $E_Y = 1$ and Poisson's ratio $\nu=0.49995$. We choose four penalty parameters $\alpha \in \bigl\{\frac{\ell_0^2}{10}, \ell_0^2, 10\ell_0^2, 100\ell_0^2\bigr\}$ and examine the first non-rigid body eigenvalue of the element stiffness matrix of the penalty stress-hybrid formulation. For each $\alpha$, we repeat the eigenvalue analysis presented in the previous section. In~\fref{fig:eigen_analysis_penalty}, the countour plots of the the first non-rigid body eigenvalue as a function of $(\gamma_1,\gamma_2)$ are shown. The plots reveal that as $\alpha$ increases, the maximum value of the eigenvalue decreases \ac{(see Figures~\ref{fig:eigen_analysis_penalty_a} and~\ref{fig:eigen_analysis_penalty_d})}. \ac{This implies that as $\alpha \to \infty$, the first non-rigid body eigenvalue becomes zero for all values of $(\gamma_1,\gamma_2)$; therefore, the element stiffness matrix will be rank deficient and lose stability in the limit.} However, for the tested values of $\alpha$, no spurious eigenvalues appear even for highly distorted elements.  

\subsection{Eigenvalue analysis for near incompressibility}\label{subsec:eigenanalysis_incompressible}
It is known that in the incompressible limit $(\nu \to 0.5)$ that the element stiffness matrix should only have one eigenvalue that tends to infinity~\cite{Simo:1993:IVA}. Elements with more than one infinite eigenvalue will experience volumetric locking. We examine the eigenvalues of the element stiffness matrix for the standard VEM, CT FEM, B-bar VEM, $15\beta$ SH-VEM, and $12\beta$ PSH-VEM on a single element. The material has Youngs modulus $E_Y=1$ and Poisson's ratio $\nu=0.4999999$.  
In~\tref{tab:eig_largest_regular}, the five largest eigenvalues of each method is presented for a regular six-noded triangular element (see~\fref{fig:example_single_elements_c}) and~\tref{tab:eig_largest_nonconvex} shows the eigenvalues for a six-noded nonconvex element (see~\fref{fig:example_single_elements_a}). The tables show that the standard virtual element approach has all five largest eigenvalues tending to infinity, which leads to severe volumetric locking. The composite element has three diverging eigenvalues, which can result in the mild locking behavior. The B-bar VEM, SH-VEM, and PSH-VEM have only a single large eigenvalue for both the regular and the nonconvex element.     

\begin{table}[!bht]
    \centering
    \resizebox{.8\textwidth}{!}{
    \begin{tabular}{|c|c|c|c|c|c|}
    \hline
    Eigenvalue & VEM & CT FEM & B-bar VEM & $15\beta$ SH-VEM & $12\beta$ PSH-VEM \\
    \hline
       1&  $1.1 \times 10^6$ & $8.1 \times 10^{-1}$ & $2.1\times 10^{-1}$ & $6.3\times 10^{-1}$ & $4.3 \times 10^{-1}$\\
       2&  $1.2 \times 10^6$& $1.2\times 10^0$ & $2.2\times 10^{-1}$ & $8.8\times 10^{-1}$ & $8.2 \times 10^{-1}$\\ 
       3&  $1.7 \times 10^6$& $3.7\times 10^5$ & $7.5\times 10^{-1}$ & $1.9\times 10^0$ & $8.4 \times 10^{-1}$ \\
       4&  $1.8 \times 10^6$& $4.2\times 10^5$ & $9.6\times 10^{-1}$ & $4.5\times 10^{1}$ & $2.8 \times 10^{0}$ \\
       5&  $5.0 \times 10^6$& $4.6\times 10^{5}$ & $4.2\times 10^{6}$ & $4.2\times 10^{6}$ & $4.2 \times 10^{6}$ \\
         \hline
    \end{tabular}}
    \caption{Comparison of the five largest eigenvalues of the element stiffness matrix on a six-noded triangular element.}
    \label{tab:eig_largest_regular}
\end{table}
\begin{table}[!bht]
    \centering
    \resizebox{.8\textwidth}{!}{
    \begin{tabular}{|c|c|c|c|c|c|}
    \hline
     Eigenvalue & VEM & CT FEM & B-bar VEM & $15\beta$ SH-VEM & $12\beta$ PSH-VEM \\
    \hline
        1& $1.2 \times 10^6$ & $5.8 \times 10^{-1}$ & $2.9\times 10^{-1}$ & $5.3\times 10^{-1}$ & $3.6 \times 10^{-1}$ \\
        2& $1.3 \times 10^6$& $1.9\times 10^0$ & $3.9\times 10^{-1}$ & $1.5\times 10^0$ & $1.3 \times 10^0$\\ 
        3& $5.8 \times 10^6$& $4.8\times 10^5$ & $2.1\times 10^0$ & $6.6\times 10^0$ & $2.0 \times 10^0$\\
        4& $5.9 \times 10^6$& $5.6\times 10^5$ & $2.2\times 10^{0}$ & $1.8\times 10^1$ & $5.1 \times 10^0$ \\
        5& $1.2 \times 10^7$& $6.8\times 10^{6}$ & $6.7\times 10^{6}$ & $6.7\times 10^{6}$ & $6.7 \times 10^{6}$ \\
         \hline
    \end{tabular}}
    \caption{Comparison of the five largest eigenvalues of the element stiffness matrix on a six-noded nonconvex element.}
    \label{tab:eig_largest_nonconvex}
\end{table}

\subsection{Thin cantilever beam}
\begin{figure}[!bht]
     \centering
     \begin{subfigure}{\textwidth}
         \centering
         \includegraphics[width=\textwidth]{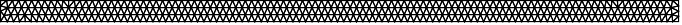}
         \caption{}
     \end{subfigure}
     \vfill
     \begin{subfigure}{\textwidth}
         \centering
         \includegraphics[width=\textwidth]{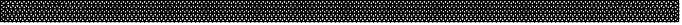}
         \caption{}
     \end{subfigure}
     \vfill
     \begin{subfigure}{\textwidth}
         \centering
         \includegraphics[width=\textwidth]{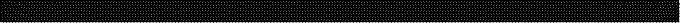}
         \caption{}
     \end{subfigure}
        \caption{Unstructured triangular meshes for the cantilever beam problem. (a) 600 elements, (b) 1600 elements, and (c) 3000 elements.  }
        \label{fig:beammesh_unstructured}
\end{figure}

\begin{figure}[!bht]
     \centering
     \begin{subfigure}{0.32\textwidth}
         \centering
         \includegraphics[width=\textwidth]{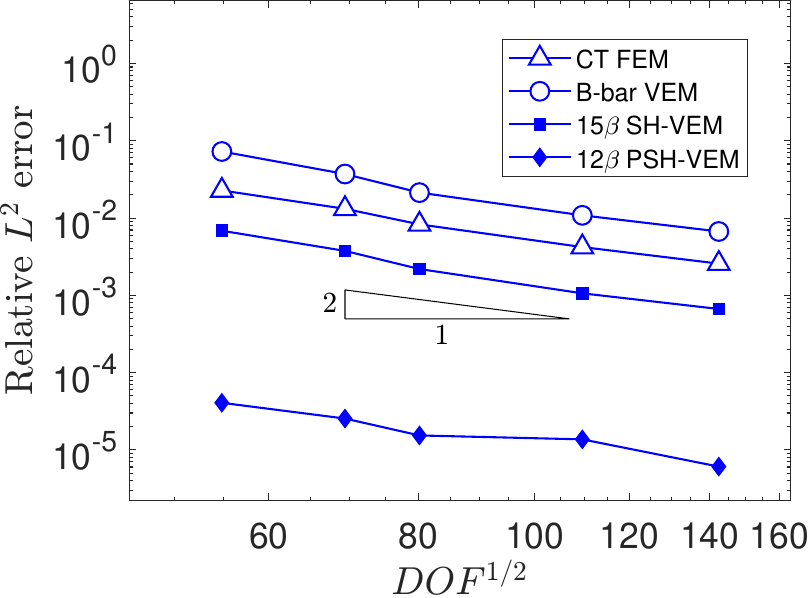}
         \caption{}
     \end{subfigure}
     \hfill
     \begin{subfigure}{0.32\textwidth}
         \centering
         \includegraphics[width=\textwidth]{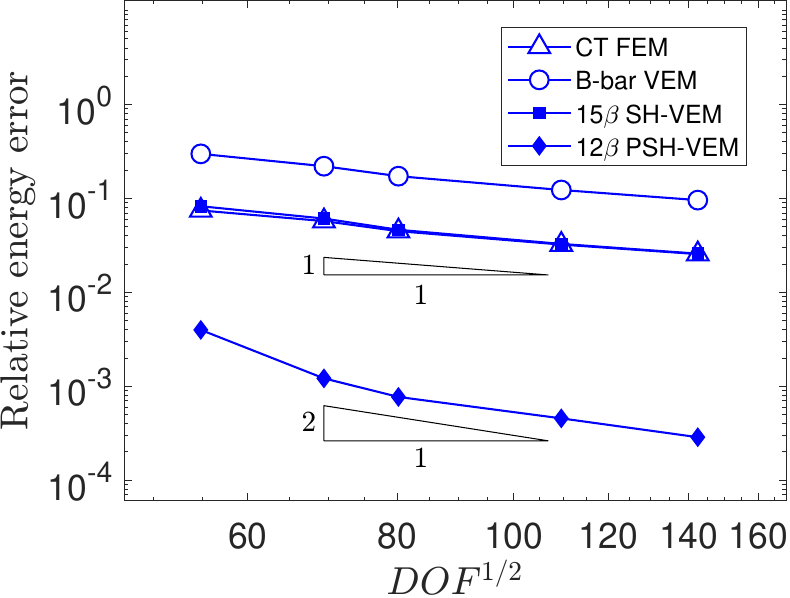}
         \caption{}
     \end{subfigure}
     \hfill
     \begin{subfigure}{0.32\textwidth}
         \centering
         \includegraphics[width=\textwidth]{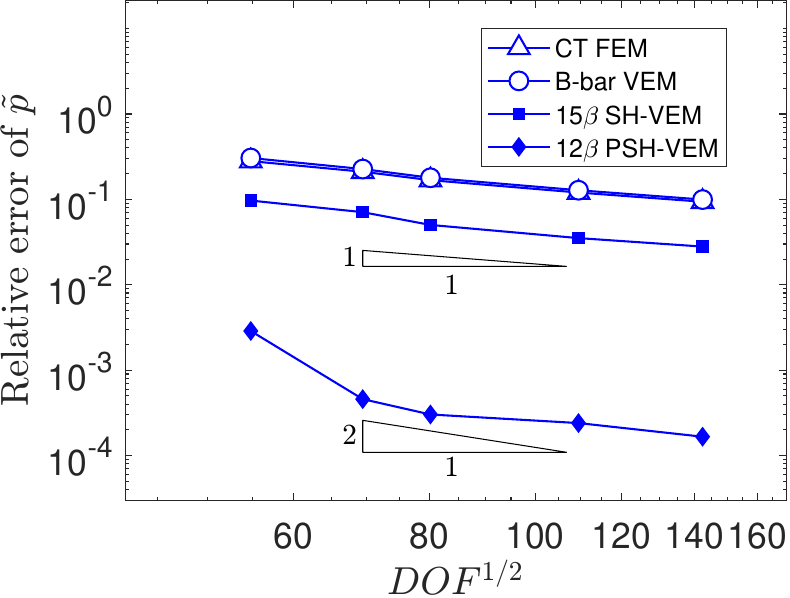}
         \caption{}
     \end{subfigure}
        \caption{Comparison of CT FEM, B-bar VEM, SH-VEM and PSH-VEM for the thin cantilever beam problem on unstructured meshes (\ac{see~\fref{fig:beammesh_unstructured}}). (a) $L^2$ error of displacement, (b) energy error, and (c) $L^2$ error of hydrostatic stress. }
        \label{fig:beam_convergence_unstructured}
\end{figure}

\begin{figure}[!bht]
    \centering
    \begin{subfigure}{0.48\textwidth}
    \centering
    \includegraphics[width=\textwidth]{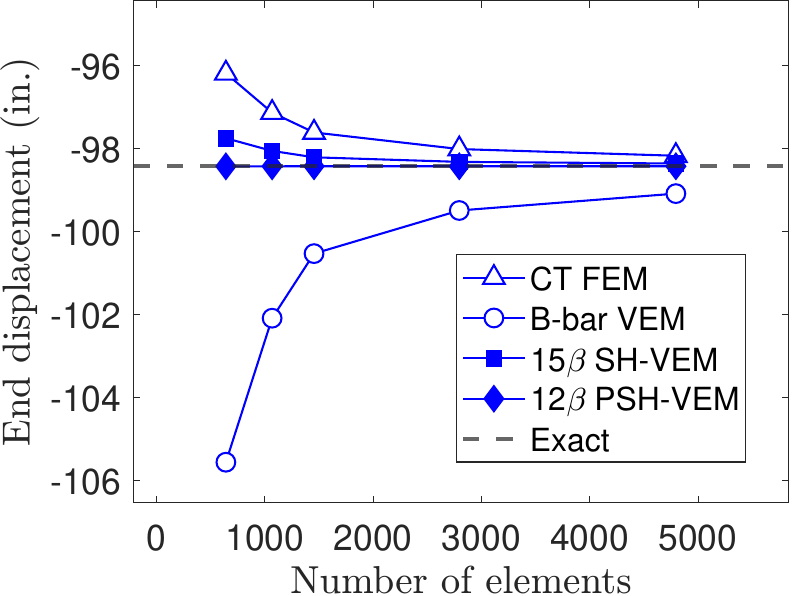}
    \caption{}\label{fig:tip_dispacement_unstructured_a}
    \end{subfigure}
    \hfill
     \begin{subfigure}{0.48\textwidth}
    \centering
    \includegraphics[width=\textwidth]{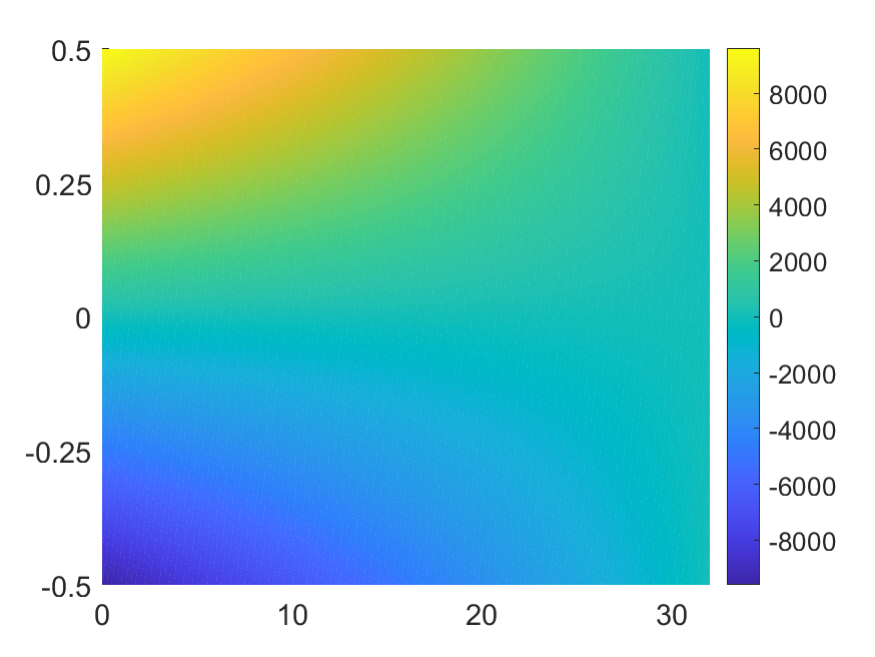}
    \caption{}\label{fig:tip_dispacement_unstructured-b}
    \end{subfigure}
    \caption{ (a) Convergence of the end displacement for the cantilever beam problem. The mesh consists of unstructured triangles (\ac{see~\fref{fig:beammesh_unstructured}}). (b) Contour plot of the hydrostatic stress for PSH-VEM.}
    \label{fig:tip_dispacement_unstructured} 
\end{figure}

We consider the problem of a thin cantilever beam subjected to a shear end load~\cite{timoshenko1951theory}. The material has Youngs modulus $E_Y = 1\times 10^5$ psi and Poisson's ratio $\nu=0.49995$. The beam has a length of $L=32$ inch, a height of $D=1$ inch and a unit thickness. A fixed displacement is applied on the left boundary, while a shear end load of $P=-100$ lbf is applied on the right boundary. We first use a set of unstructured triangular meshes. In~\fref{fig:beammesh_unstructured}, we show examples of the unstructred meshes and in~\fref{fig:beam_convergence_unstructured} we show the convergence of the four methods in the displacement $L^2$ norm, energy seminorm, and $L^2$ norm of hydrostatic stress. \fref{fig:tip_dispacement_unstructured} shows the convergence of the end displacement and the contour plot of hydrostatic stress for the SH-VEM. For methods that experience volumetric locking, non-physical oscillations will appear in the hydrostatic stress field. The plots reveal that the two stress-hybrid approaches yield better results, with the PSH-VEM having superior convergence and accuracy in the displacement, energy, and hydrostatic stress.


\begin{figure}[!htb]
\begin{minipage}{.48\textwidth}
     \centering
     \begin{subfigure}{\textwidth}
         \centering
         \includegraphics[width=\textwidth]{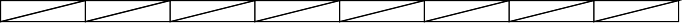}
         \caption{}
     \end{subfigure}
     \vfill
     \begin{subfigure}{\textwidth}
         \centering
         \includegraphics[width=\textwidth]{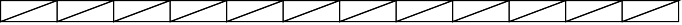}
         \caption{}
     \end{subfigure}
     \vfill
     \begin{subfigure}{\textwidth}
         \centering
         \includegraphics[width=\textwidth]{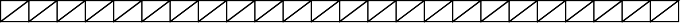}
         \caption{}
     \end{subfigure}
          \vfill
     \begin{subfigure}{\textwidth}
         \centering
         \includegraphics[width=\textwidth]{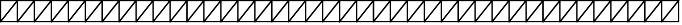}
         \caption{}
     \end{subfigure}
     \end{minipage}
     \hfill
     \begin{minipage}{.48\textwidth}
     \centering
     \begin{subfigure}{\textwidth}
         \centering
         \includegraphics[width=\textwidth]{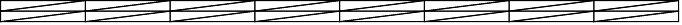}
         \caption{}
     \end{subfigure}
     \vfill
     \begin{subfigure}{\textwidth}
         \centering
         \includegraphics[width=\textwidth]{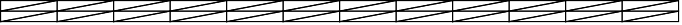}
         \caption{}
     \end{subfigure}
     \vfill
     \begin{subfigure}{\textwidth}
         \centering
         \includegraphics[width=\textwidth]{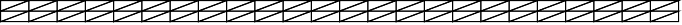}
         \caption{}
     \end{subfigure}
          \vfill
     \begin{subfigure}{\textwidth}
         \centering
         \includegraphics[width=\textwidth]{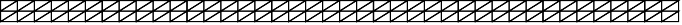}
         \caption{}
     \end{subfigure}
     \end{minipage}
        \caption{Structured meshes for the cantilever beam problem.
        Mesh is refined along the length with (a)-(d) 1 element along the
        height and (e)-(h) 2 elements along the height.}
        \label{fig:beammesh_structured_v1}
\end{figure}

\begin{figure}[!bht]
    \centering
    \begin{subfigure}{0.48\textwidth}
    \centering
    \includegraphics[width=\textwidth]{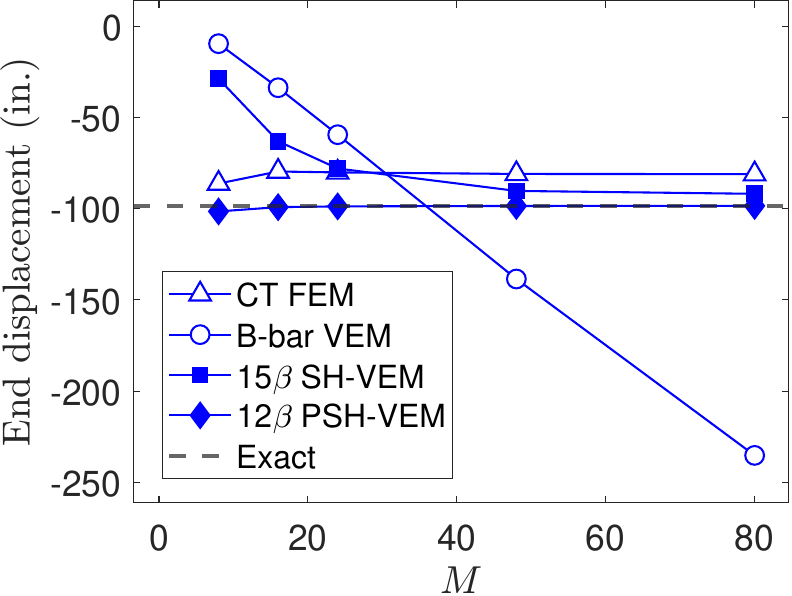}
    \caption{}\label{fig:tip_dispacement_structured_v1-a}
    \end{subfigure}
    \hfill
     \begin{subfigure}{0.48\textwidth}
    \centering
    \includegraphics[width=\textwidth]{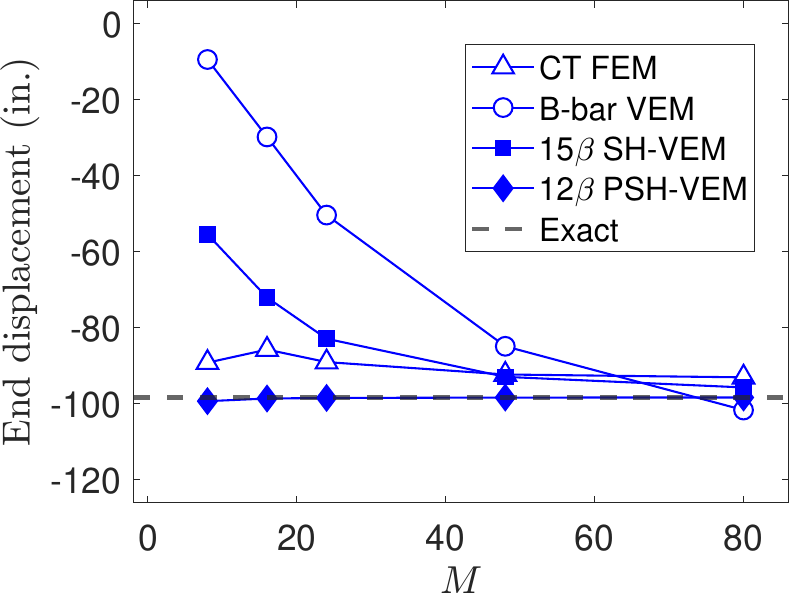}
    \caption{}\label{fig:tip_dispacement_structured_v1-b}
    \end{subfigure}
    \caption{ Convergence of the end displacement for the
    cantilever beam problem. The mesh consists of
    $M\times N$ right triangles, where $M$ is the number of elements along the length of the beam (\ac{see~\fref{fig:beammesh_structured_v1}}). (a) $N = 1$ and (b) $N = 2$.}
    \label{fig:tip_dispacement_structured_v1} 
\end{figure}

In the following tests, we examine the four methods on meshes with $N=1$ and $N=2$ elements along the height of the beam. For many formulations, these meshes will lead to overly stiff displacements and shear locking. The first set of meshes is constructed by taking a uniform quadrilateral mesh and splitting each element into two right triangles (see~\fref{fig:beammesh_structured_v1}). The convergence of the end displacement of the four methods is presented in~\fref{fig:tip_dispacement_structured_v1}. For $N=1$, the two stress-hybrid methods are converging to the exact solution, while CT FEM experiences shear locking and B-bar VEM diverges from the exact solution (see~\fref{fig:tip_dispacement_structured_v1-a}) and for $N=2$, B-bar VEM is still overly flexible (see~\fref{fig:tip_dispacement_structured_v1-b}). In both cases, the PSH-VEM shows far superior accuracy in displacement even on coarse meshes.  

\begin{figure}[!bht]
\begin{minipage}{.48\textwidth}
     \centering
     \begin{subfigure}{\textwidth}
         \centering
         \includegraphics[width=\textwidth]{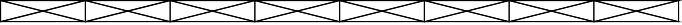}
         \caption{}
     \end{subfigure}
     \vfill
     \begin{subfigure}{\textwidth}
         \centering
         \includegraphics[width=\textwidth]{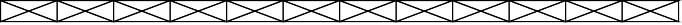}
         \caption{}
     \end{subfigure}
     \vfill
     \begin{subfigure}{\textwidth}
         \centering
         \includegraphics[width=\textwidth]{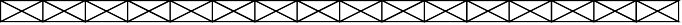}
         \caption{}
     \end{subfigure}
          \vfill
     \begin{subfigure}{\textwidth}
         \centering
         \includegraphics[width=\textwidth]{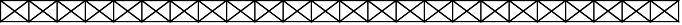}
         \caption{}
     \end{subfigure}
     \end{minipage}
     \hfill
     \begin{minipage}{.48\textwidth}
     \centering
     \begin{subfigure}{\textwidth}
         \centering
         \includegraphics[width=\textwidth]{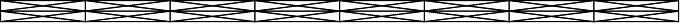}
         \caption{}
     \end{subfigure}
     \vfill
     \begin{subfigure}{\textwidth}
         \centering
         \includegraphics[width=\textwidth]{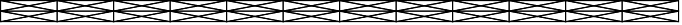}
         \caption{}
     \end{subfigure}
     \vfill
     \begin{subfigure}{\textwidth}
         \centering
         \includegraphics[width=\textwidth]{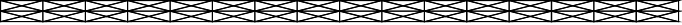}
         \caption{}
     \end{subfigure}
          \vfill
     \begin{subfigure}{\textwidth}
         \centering
         \includegraphics[width=\textwidth]{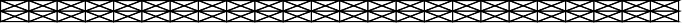}
         \caption{}
     \end{subfigure}
     \end{minipage}
        \caption{Structured meshes for the cantilever beam problem.
        Mesh is refined along the length with (a)-(d) 1 element along the
        height and (e)-(h) 2 elements along the height.}
        \label{fig:beammesh_structured_v2}
\end{figure}

\begin{figure}[!bht]
    \centering
    \begin{subfigure}{0.48\textwidth}
    \centering
    \includegraphics[width=\textwidth]{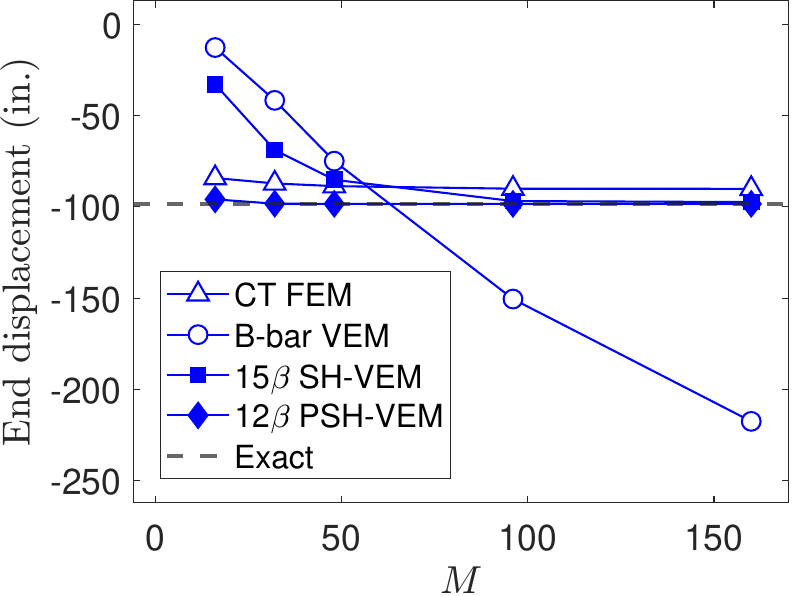}
    \caption{}\label{fig:tip_dispacement_structured_v2-a}
    \end{subfigure}
    \hfill
     \begin{subfigure}{0.48\textwidth}
    \centering
    \includegraphics[width=\textwidth]{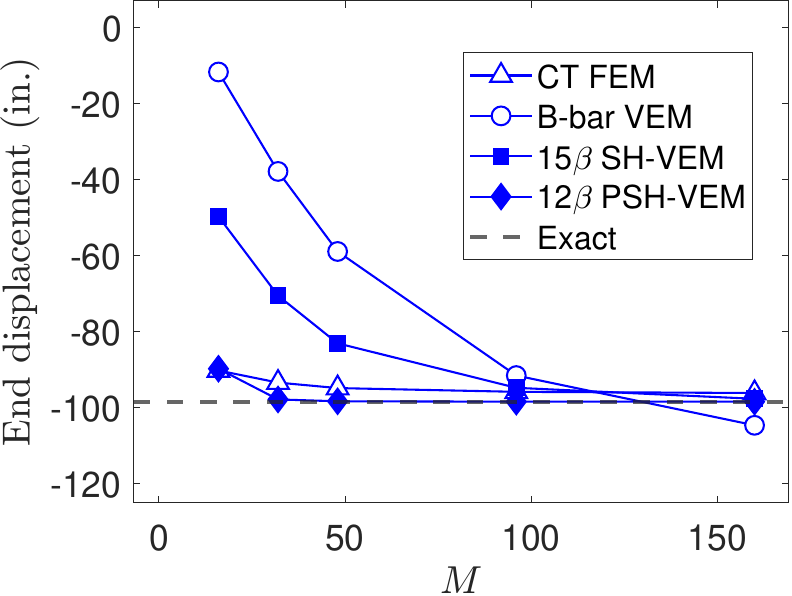}
    \caption{}\label{fig:tip_dispacement_structured_v2-b}
    \end{subfigure}
    \caption{ Convergence of the end displacement for the
    cantilever beam problem. The mesh consists of
    $M\times N$ triangles, where $M$ is the number of elements along the length of the beam (\ac{see~\fref{fig:beammesh_structured_v2}}). (a) $N = 1$ and (b) $N = 2$.}
    \label{fig:tip_dispacement_structured_v2} 
\end{figure}

The second set of meshes is constructed by taking a uniform quadrilateral mesh and splitting each element along the two diagonals to form four triangles. The meshes are shown in~\fref{fig:beammesh_structured_v2} and the convergence of the end displacement is depicted in~\fref{fig:tip_dispacement_structured_v2}. The plots show that SH-VEM and PSH-VEM converge to the exact solution for both $N=1$ and $N=2$; however, the SH-VEM is much stiffer and less accurate than the penalty approach. CT FEM suffers from locking, while B-bar fails to converge for the case of a single element along the height. 

\begin{figure}[!bht]
\begin{minipage}{.35\textwidth}
     \centering
     \begin{subfigure}{\textwidth}
         \centering
         \includegraphics[width=\textwidth]{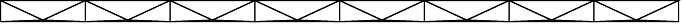}
         \caption{}
     \end{subfigure}
     \vfill
     \begin{subfigure}{\textwidth}
         \centering
         \includegraphics[width=\textwidth]{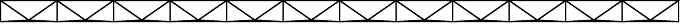}
         \caption{}
     \end{subfigure}
     \vfill
     \begin{subfigure}{\textwidth}
         \centering
         \includegraphics[width=\textwidth]{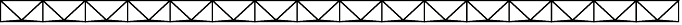}
         \caption{}
     \end{subfigure}
          \vfill
     \begin{subfigure}{\textwidth}
         \centering
         \includegraphics[width=\textwidth]{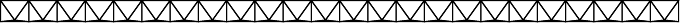}
         \caption{}
     \end{subfigure}
     \end{minipage}
     \hfill
     \begin{minipage}{.35\textwidth}
     \centering
     \begin{subfigure}{\textwidth}
         \centering
         \includegraphics[width=\textwidth]{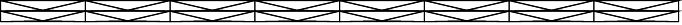}
         \caption{}
     \end{subfigure}
     \vfill
     \begin{subfigure}{\textwidth}
         \centering
         \includegraphics[width=\textwidth]{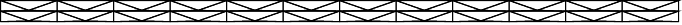}
         \caption{}
     \end{subfigure}
     \vfill
     \begin{subfigure}{\textwidth}
         \centering
         \includegraphics[width=\textwidth]{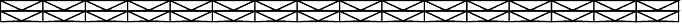}
         \caption{}
     \end{subfigure}
          \vfill
     \begin{subfigure}{\textwidth}
         \centering
         \includegraphics[width=\textwidth]{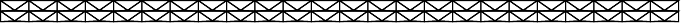}
         \caption{}
     \end{subfigure}
     \end{minipage}
     \hfill
     \begin{minipage}{.28\textwidth}
     \centering
     \begin{subfigure}{\textwidth}
         \centering
         \includegraphics[width=\textwidth]{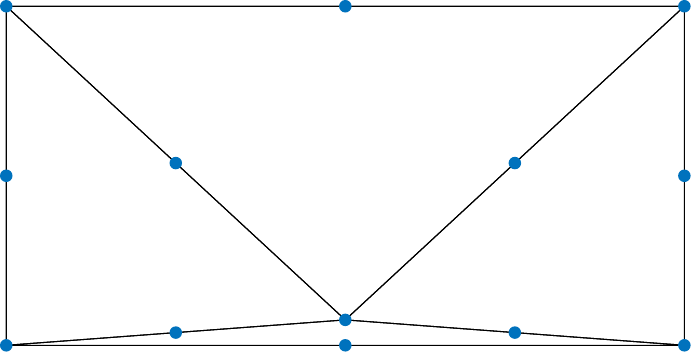}
         \caption{}
     \end{subfigure}
     \end{minipage}
        \caption{Structured nearly degenerate meshes for the cantilever beam problem.
        Mesh is refined along the length with (a)-(d) 1 element along the
        height and (e)-(h) 2 elements along the height. (i) Magnification of a single element split into four six-noded triangles.}
        \label{fig:beammesh_degen_v1}
\end{figure}

\begin{figure}[!bht]
    \centering
    \begin{subfigure}{0.48\textwidth}
    \centering
    \includegraphics[width=\textwidth]{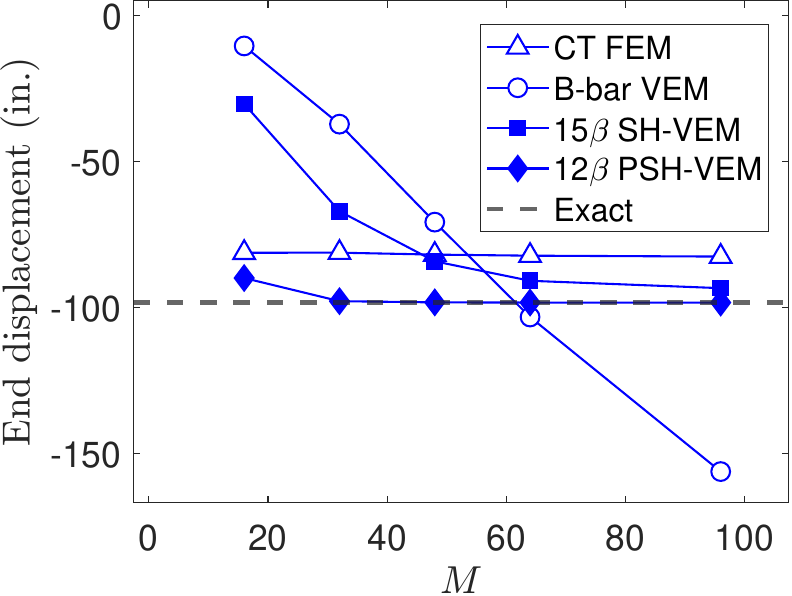}
    \caption{}\label{fig:tip_dispacement_degen_v1-a}
    \end{subfigure}
    \hfill
     \begin{subfigure}{0.48\textwidth}
    \centering
    \includegraphics[width=\textwidth]{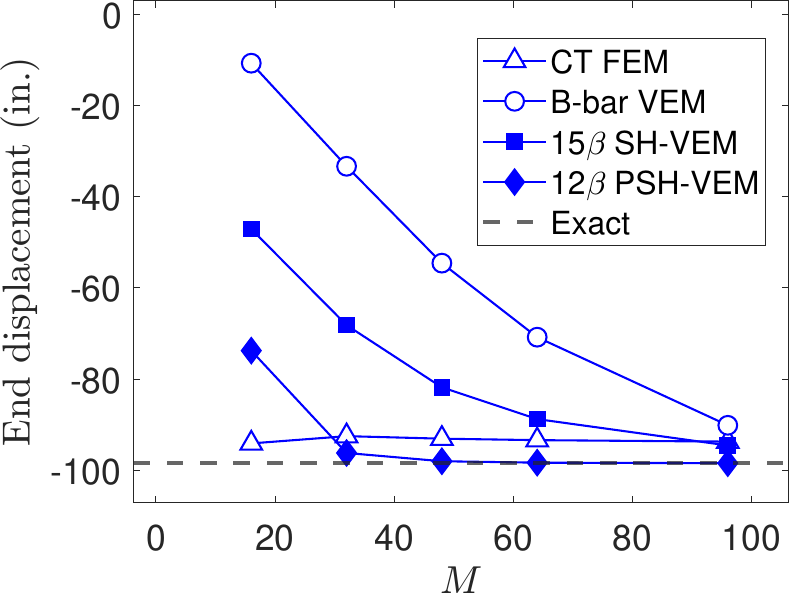}
    \caption{}\label{fig:tip_dispacement_degen_v1-b}
    \end{subfigure}
    \caption{ Convergence of the end displacement for the
    cantilever beam problem. The mesh consists of
    $M\times N$ triangles with some nearly degenerate, where $M$ is the number of elements along the length of the beam (\ac{see~\fref{fig:beammesh_degen_v1}}). (a) $N = 1$ and (b) $N = 2$.}
    \label{fig:tip_dispacement_degen_v1} 
\end{figure}

The third set of meshes is constructed by taking the previous mesh and collapsing one of the triangles to a nearly degenerate triangle (see~\fref{fig:beammesh_degen_v1}). The plots showing the convergence of the end displacement is given in~\fref{fig:tip_dispacement_degen_v1}. The plots reveal that for these meshes, the stress-hybrid methods are converging to the exact solution, while CT FEM and B-bar VEM do not converge in the case of a single element along the height. When using two elements along the height, the B-bar formulation is tending to the exact solution but is not accurate, while CT FEM still suffers from locking and does not converge. The PSH-VEM again attains a much more accurate solution than the other methods; however, it appears to be stiff for the coarsest mesh when $N=2$.   


\begin{figure}[!bht]
\begin{minipage}{.48\textwidth}
     \centering
     \begin{subfigure}{\textwidth}
         \centering
         \includegraphics[width=\textwidth]{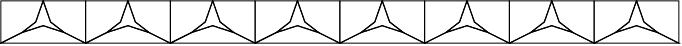}
         \caption{}
     \end{subfigure}
     \vfill
     \begin{subfigure}{\textwidth}
         \centering
         \includegraphics[width=\textwidth]{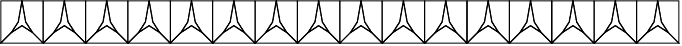}
         \caption{}
     \end{subfigure}
     \vfill
     \begin{subfigure}{\textwidth}
         \centering
         \includegraphics[width=\textwidth]{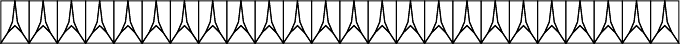}
         \caption{}
     \end{subfigure}
          \vfill
     \begin{subfigure}{\textwidth}
         \centering
         \includegraphics[width=\textwidth]{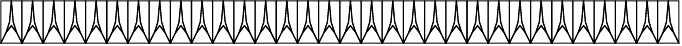}
         \caption{}
     \end{subfigure}
     \end{minipage}
     \hfill
     \begin{minipage}{.48\textwidth}
     \centering
     \begin{subfigure}{\textwidth}
         \centering
         \includegraphics[width=\textwidth]{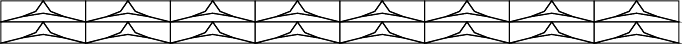}
         \caption{}
     \end{subfigure}
     \vfill
     \begin{subfigure}{\textwidth}
         \centering
         \includegraphics[width=\textwidth]{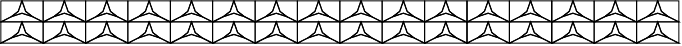}
         \caption{}
     \end{subfigure}
     \vfill
     \begin{subfigure}{\textwidth}
         \centering
         \includegraphics[width=\textwidth]{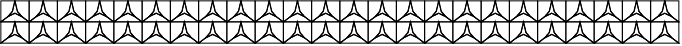}
         \caption{}
     \end{subfigure}
          \vfill
     \begin{subfigure}{\textwidth}
         \centering
         \includegraphics[width=\textwidth]{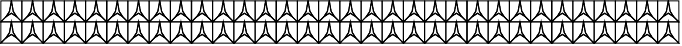}
         \caption{}
     \end{subfigure}
     \end{minipage}
        \caption{Meshes with nonconvex elements for the cantilever beam problem.
        Mesh is refined along the length with (a)-(d) 1 element along the
        height and (e)-(h) 2 elements along the height.}
        \label{fig:beammesh_nonconvex}
\end{figure}

\begin{figure}[!bht]
    \centering
    \begin{subfigure}{0.48\textwidth}
    \centering
    \includegraphics[width=\textwidth]{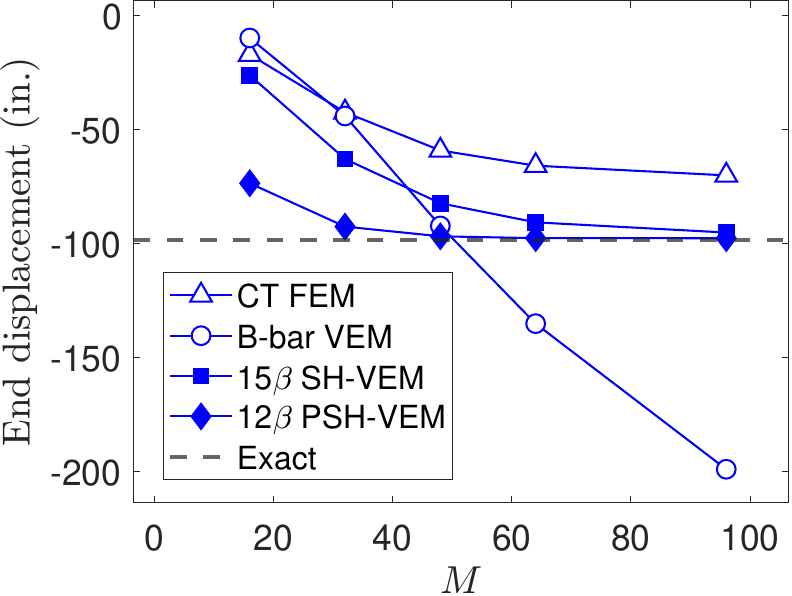}
    \caption{}\label{fig:tip_dispacement_nonconvex-a}
    \end{subfigure}
    \hfill
     \begin{subfigure}{0.48\textwidth}
    \centering
    \includegraphics[width=\textwidth]{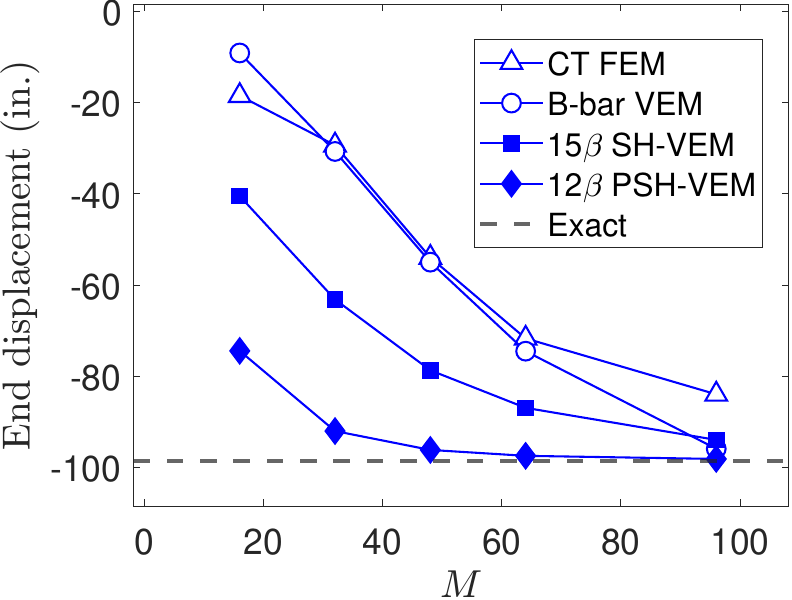}
    \caption{}\label{fig:tip_dispacement_nonconvex-b}
    \end{subfigure}
    \caption{ Convergence of the end displacement for the
    cantilever beam problem. The mesh consists of
    $M\times N$ convex and nonconvex elements, where $M$ is the number of elements along the length of the beam (\ac{see~\fref{fig:beammesh_nonconvex}}). (a) $N = 1$ and (b) $N = 2$.}
    \label{fig:tip_dispacement_nonconvex} 
\end{figure}

One benefit of the virtual element formulation is that it allows for very general element shapes. In particular, for the final set of meshes we use a mixture of distorted nonconvex hexagons and convex hexagons. A few representative meshes are shown in~\fref{fig:beammesh_nonconvex}. In~\fref{fig:tip_dispacement_nonconvex}, we present the convergence of the end displacement for the four methods. In both cases $N=1$ and $N=2$, SH-VEM and PSH-VEM are convergent, but CT FEM is overly stiff.    
\begin{remark}
    The solutions produced by $15\beta$ SH-VEM are much stiffer than those found using the $5\beta$ formulation on quadrilaterals in~\cite{chen:2023:shv}. In particular, for uniform rectangular meshes with $N=1$, the $5\beta$ SH-VEM converges to nearly
    the exact solution with just $M=8$ elements along the length of the beam. However, the PSH-VEM offers coarse mesh accuracy that is similar to the $5\beta$ SH-VEM.   
\end{remark}

\begin{figure}[!htb]
    \centering
    \includegraphics[width=0.55\textwidth]{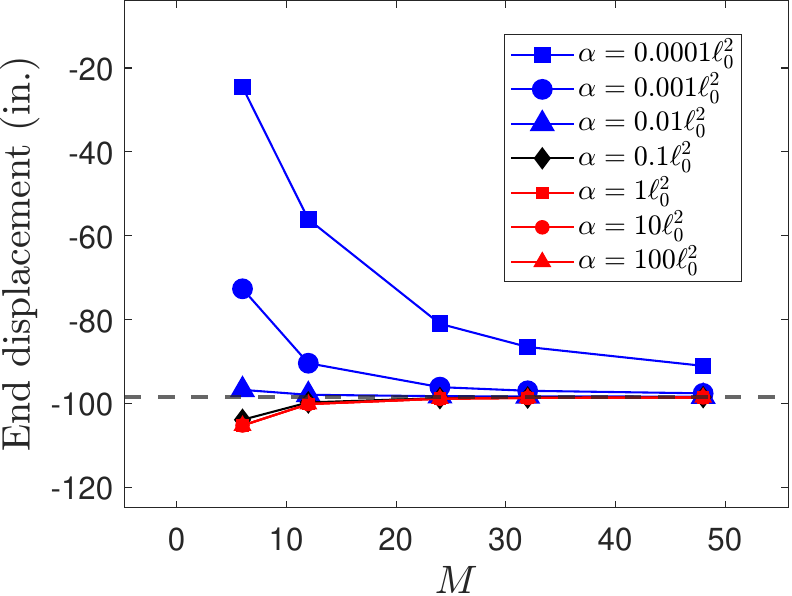}
    \caption{Comparison of the convergence of end displacement for PSH-VEM when using different choices of penalty parameter $\alpha$. The mesh consists of $M\times 1$ right triangles, where $M$ is the number of elements along the length of the beam (\ac{see~\fref{fig:beammesh_structured_v1}}).}
    \label{fig:beam_displacement_penalty_vary}
\end{figure}

The penalty stress-hybrid approach provides the best accuracy for the cantilever beam problem; however, it relies on the choice of a suitable penalty parameter $\alpha$. Therefore, we examine the sensitivity of the bending solution of PSH-VEM to the penalty parameter. For this test, we use the same material properties $E_Y=1\times 10^5$ psi and $\nu=0.49995$. The problem is solved on a mesh consisting of $M\times 1$ right triangles (see~\fref{fig:beammesh_structured_v1}). The default penalty parameter is given by $\alpha = \frac{10^4\ell_0^2}{E_Y} = 10^{-1}\ell_0^2$. To test the sensitivity of the parameter, we vary the penalty parameter three orders of magnitude on each side. That is, $\alpha$ is varied from $10^{-4}\ell_0^2$ to $10^2\ell_0^2$. In~\fref{fig:beam_displacement_penalty_vary}, the convergence of the end displacements for different values of $\alpha$ is presented. The plot shows that increasing $\alpha$ above the value $10^{-1}\ell_0^2$ does not greatly affect the solution; however, decreasing $\alpha$ makes the solution stiffer. If we decrease the value of $\alpha$ by three orders of magnitude to $\alpha = 10^{-4}\ell_0^2$, the solution becomes overly stiff and produces much larger errors than the other values of $\alpha$. 

\subsection{Cook's membrane}
\begin{figure}[!bht]
     \centering
     \begin{subfigure}{.32\textwidth}
         \centering
         \includegraphics[width=\textwidth]{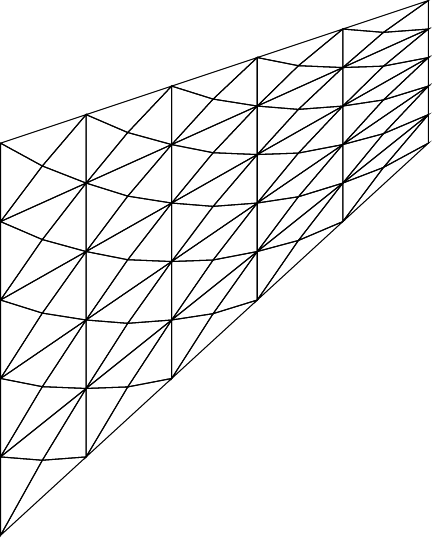}
         \caption{}
     \end{subfigure}
     \hfill
     \begin{subfigure}{.32\textwidth}
         \centering
         \includegraphics[width=\textwidth]{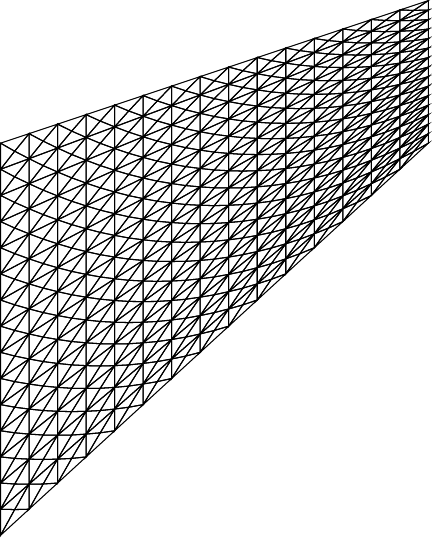}
         \caption{}
     \end{subfigure}
     \hfill
     \begin{subfigure}{.32\textwidth}
         \centering
         \includegraphics[width=\textwidth]{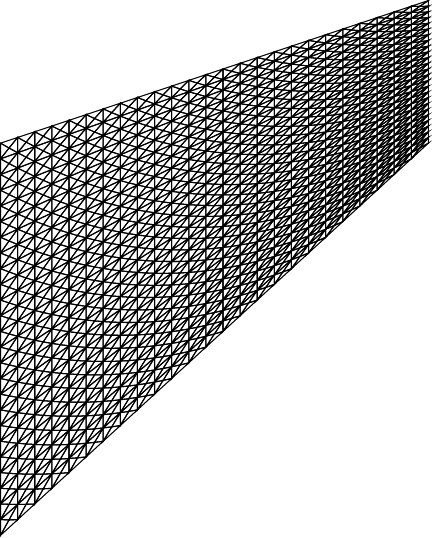}
         \caption{}
     \end{subfigure}
        \caption{Structured triangular meshes for the Cook's membrane problem. (a) 100 elements, (b) 1000 elements, and (c) 2500 elements.  }
        \label{fig:cooksmesh_structured}
\end{figure}

\begin{figure}[!bht]
    \centering
    \begin{subfigure}{0.48\textwidth}
    \centering
    \includegraphics[width=\textwidth]{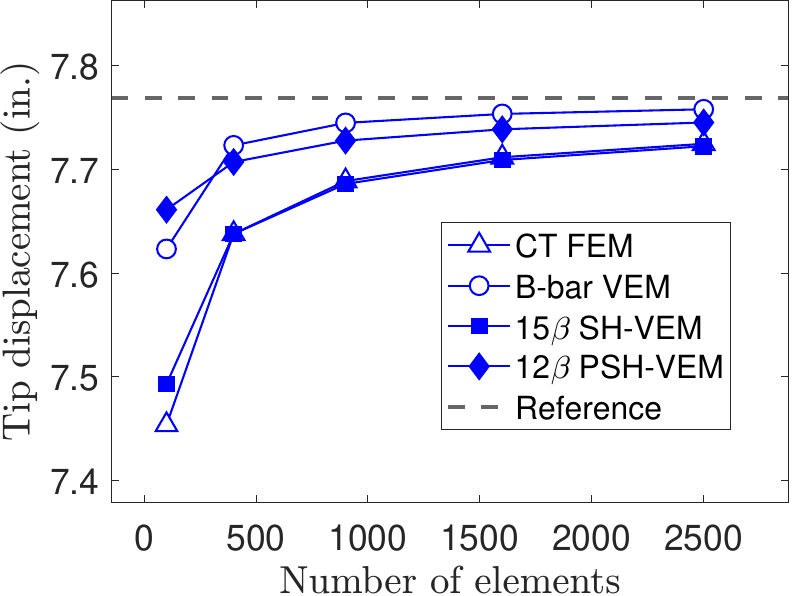}
    \caption{}\label{fig:cook_tip_dispacement_structured-a}
    \end{subfigure}
     \hfill
     \begin{subfigure}{0.48\textwidth}
         \centering
         \includegraphics[width=\textwidth]{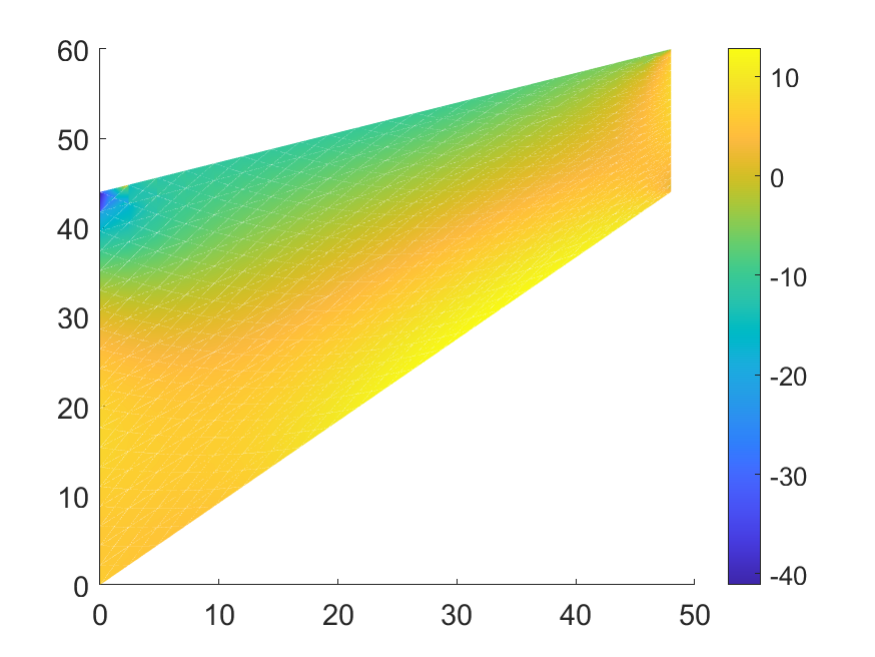}
         \caption{}\label{fig:cook_tip_dispacement_structured-b}
     \end{subfigure}
    \caption{(a) Convergence of the tip displacement for Cook's membrane problem. The mesh consists of structured triangles (\ac{see~\fref{fig:cooksmesh_structured}}). (b) Contour plot of the hydrostatic stress for PSH-VEM.}
    \label{fig:cook_tip_dispacement_structured}
\end{figure}

Next, we consider the problem of the Cook's membrane under shear load~\cite{Cook:1974:jsd}. This problem is commonly used to test the combined bending and shearing capabilities of an element. The material has Youngs modulus $E_Y=250$ psi and Poisson's ratio $\nu=0.49995$. The left boundary of the membrane is clamped and a shear load of $F=6.25$ lbf per unit length is applied along the right edge. The first set of meshes comprise of structured triangular meshes; a few sample meshes are shown in~\fref{fig:cooksmesh_structured}.~\fref{fig:cook_tip_dispacement_structured} presents the convergence of the tip displacement of the four methods and the hydrostatic stress of the PSH-VEM. The plot in~\fref{fig:cook_tip_dispacement_structured-a} shows that, unlike in the cantilever beam problem, the CT FEM, SH-VEM, and PSH-VEM all perform worse than the B-bar formulation. \fref{fig:cook_tip_dispacement_structured-b} shows that the PSH-VEM is able to produce a relatively smooth hydrostatic stress field with no visible signs of volumetric locking.   

\begin{figure}[!bht]
     \centering
     \begin{subfigure}{.32\textwidth}
         \centering
         \includegraphics[width=\textwidth]{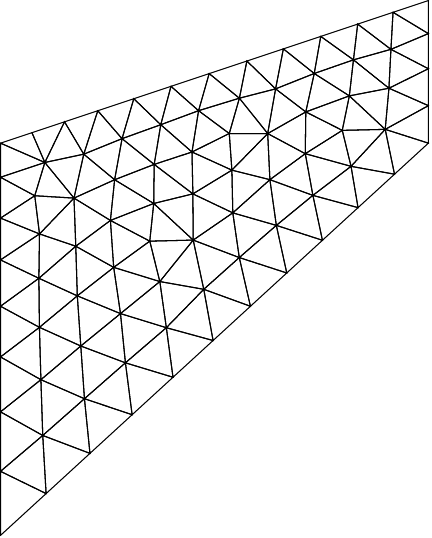}
         \caption{}
     \end{subfigure}
     \hfill
     \begin{subfigure}{.32\textwidth}
         \centering
         \includegraphics[width=\textwidth]{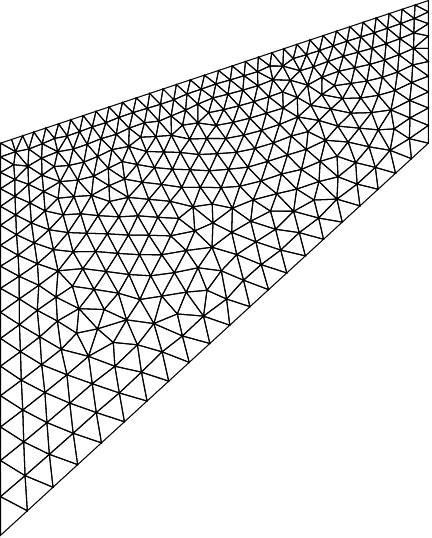}
         \caption{}
     \end{subfigure}
     \hfill
     \begin{subfigure}{.32\textwidth}
         \centering
         \includegraphics[width=\textwidth]{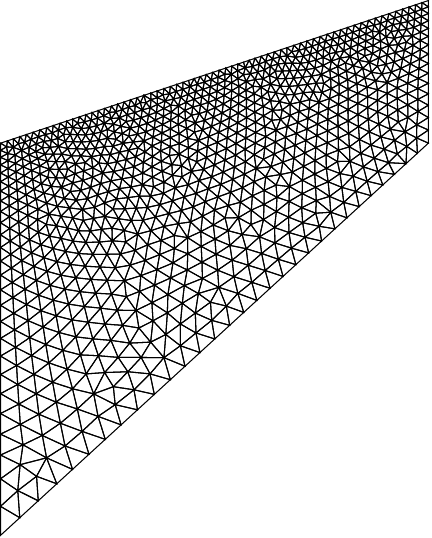}
         \caption{}
     \end{subfigure}
        \caption{Unstructured triangular meshes for the Cook's membrane problem. (a) 100 elements, (b) 700 elements, and (c) 2000 elements.  }
        \label{fig:cooksmesh_unstructured}
\end{figure}

\begin{figure}[!bht]
    \centering
    \begin{subfigure}{0.48\textwidth}
    \centering
    \includegraphics[width=\textwidth]{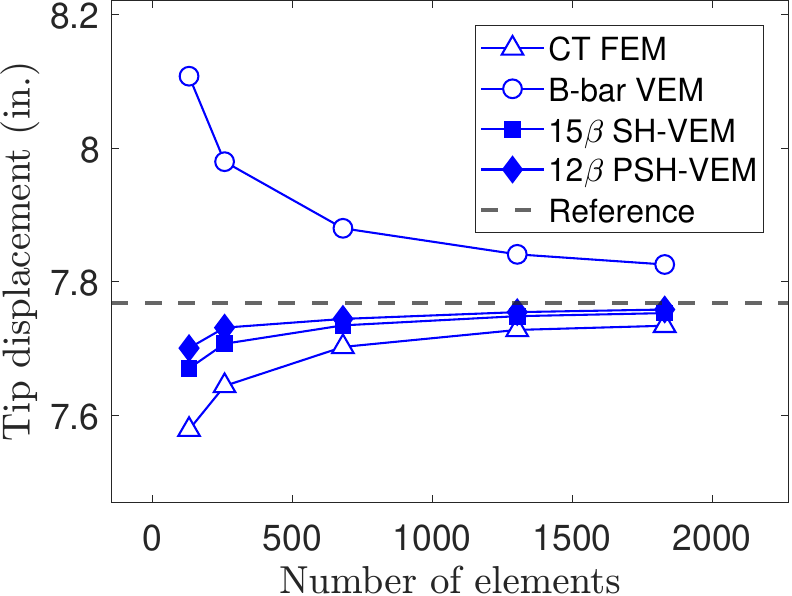}
    \caption{}
    \end{subfigure}
     \hfill
     \begin{subfigure}{0.48\textwidth}
         \centering
         \includegraphics[width=\textwidth]{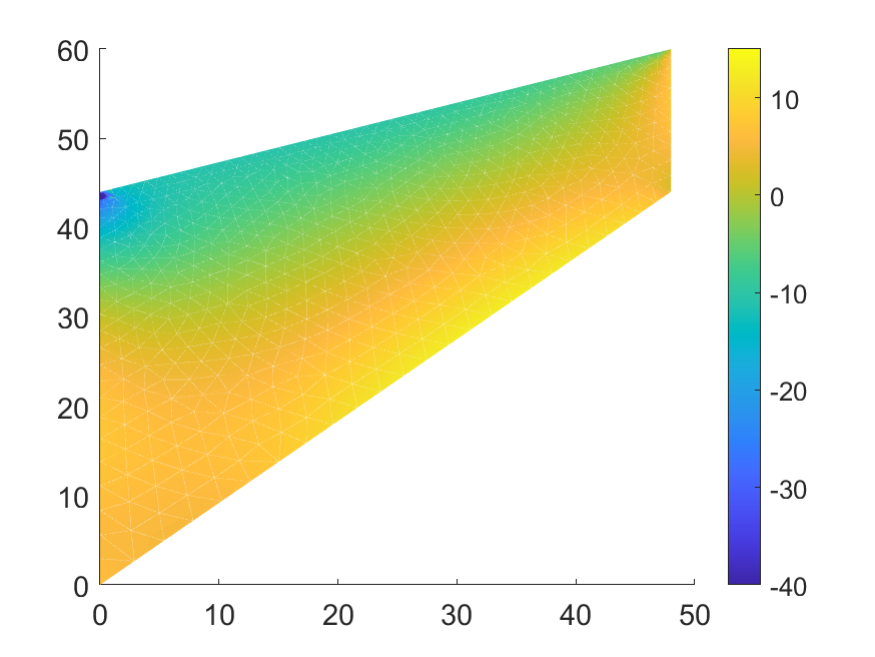}
         \caption{}
     \end{subfigure}
    \caption{(a) Convergence of the tip displacement for Cook's membrane problem. The mesh consists of unstructured triangles (\ac{see~\fref{fig:cooksmesh_unstructured}}). (b) Contour plot of the hydrostatic stress for PSH-VEM.}
    \label{fig:cook_tip_dispacement_unstructured}
\end{figure}

The second set of meshes that we test consists of unstructured triangles. In~\fref{fig:cooksmesh_unstructured}, sample meshes are shown. The convergence of tip displacement and the contour plot of PSH-VEM is given in~\fref{fig:cook_tip_dispacement_unstructured}. The plots show that for the unstructured mesh, SH-VEM and PSH-VEM have the fastest convergence, while B-Bar VEM becomes very flexible and converges from above. The hydrostatic stress contours of PSH-VEM remain smooth and agree with the contours found on structured meshes. 
\subsection{Plate with a circular hole}
\begin{figure}[!bht]
     \centering
     \begin{subfigure}{.32\textwidth}
         \centering
         \includegraphics[width=\textwidth]{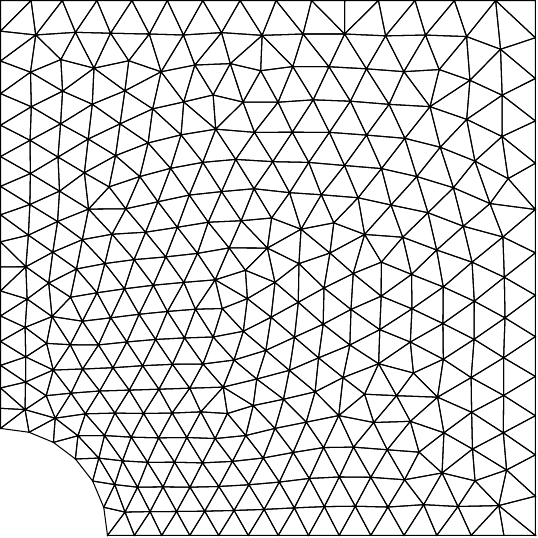}
         \caption{}
     \end{subfigure}
     \hfill
     \begin{subfigure}{.32\textwidth}
         \centering
         \includegraphics[width=\textwidth]{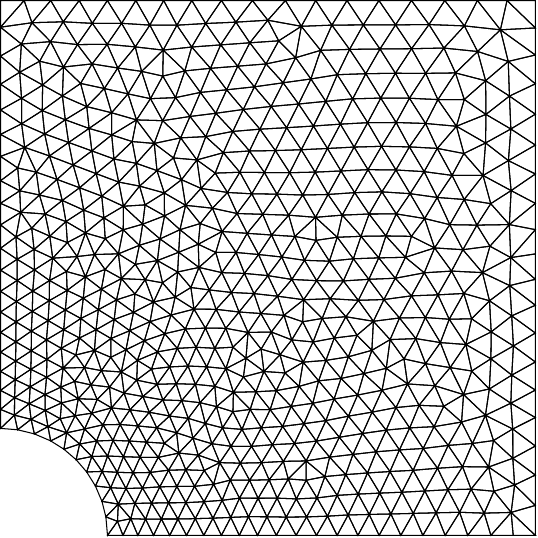}
         \caption{}
     \end{subfigure}
     \hfill
     \begin{subfigure}{.32\textwidth}
         \centering
         \includegraphics[width=\textwidth]{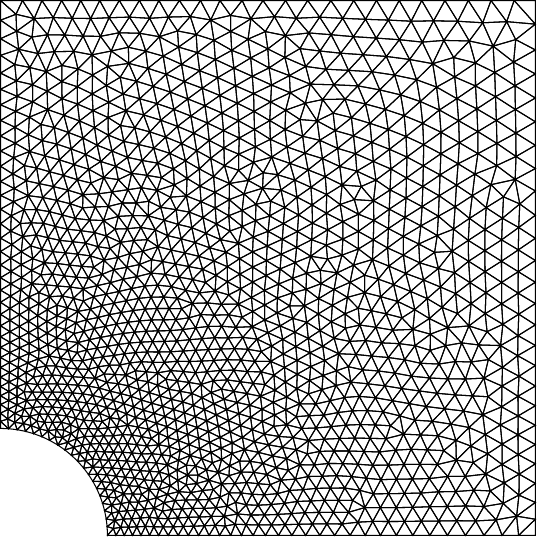}
         \caption{}
     \end{subfigure}
        \caption{Unstructured triangular meshes for the plate with a hole problem. (a) 500 elements, (b) 1000 elements, and (c) 3000 elements.  }
        \label{fig:platemesh_unstructured}
\end{figure}

\begin{figure}[!bht]
     \centering
     \begin{subfigure}{0.32\textwidth}
         \centering
         \includegraphics[width=\textwidth]{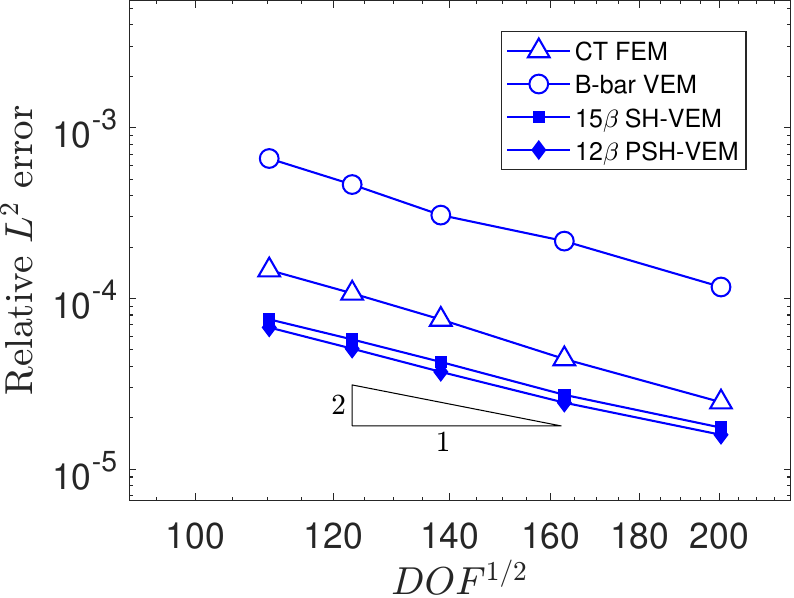}
         \caption{}
     \end{subfigure}
     \hfill
     \begin{subfigure}{0.32\textwidth}
         \centering
         \includegraphics[width=\textwidth]{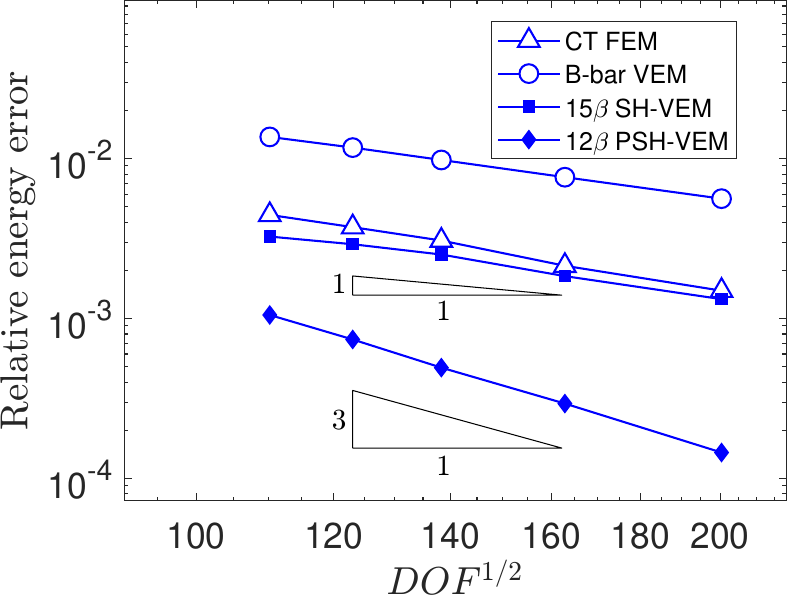}
         \caption{}
     \end{subfigure}
     \hfill
     \begin{subfigure}{0.32\textwidth}
         \centering
         \includegraphics[width=\textwidth]{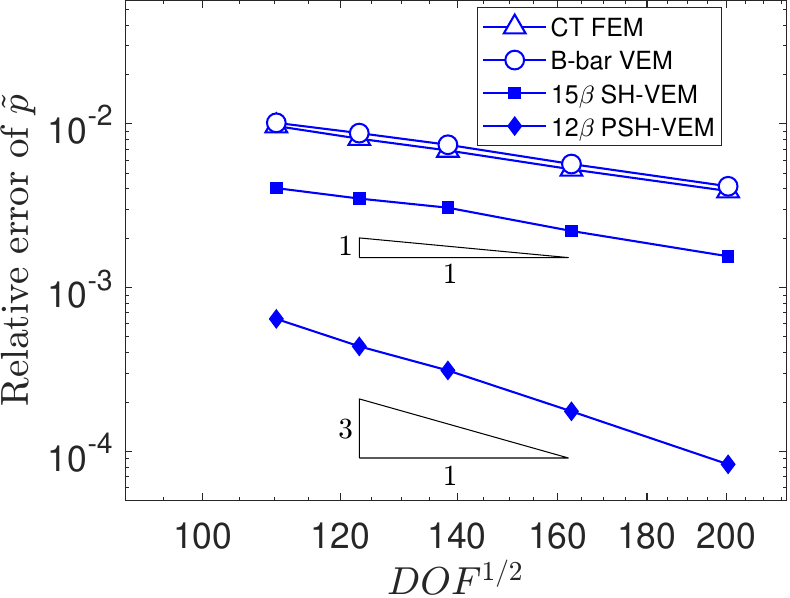}
         \caption{}
     \end{subfigure}
        \caption{Comparison of CT FEM, B-bar VEM, SH-VEM, and PSH-VEM for the plate with a hole problem on unstructured meshes (\ac{see~\fref{fig:platemesh_unstructured}}). (a) $L^2$ error of displacement, (b) energy error, and (c) $L^2$ error of hydrostatic stress. }
        \label{fig:plate_convergence_unstructured}
\end{figure}

\begin{figure}[!bht]
    \centering
    \begin{subfigure}{0.48\textwidth}
    \centering
    \includegraphics[width=\textwidth]{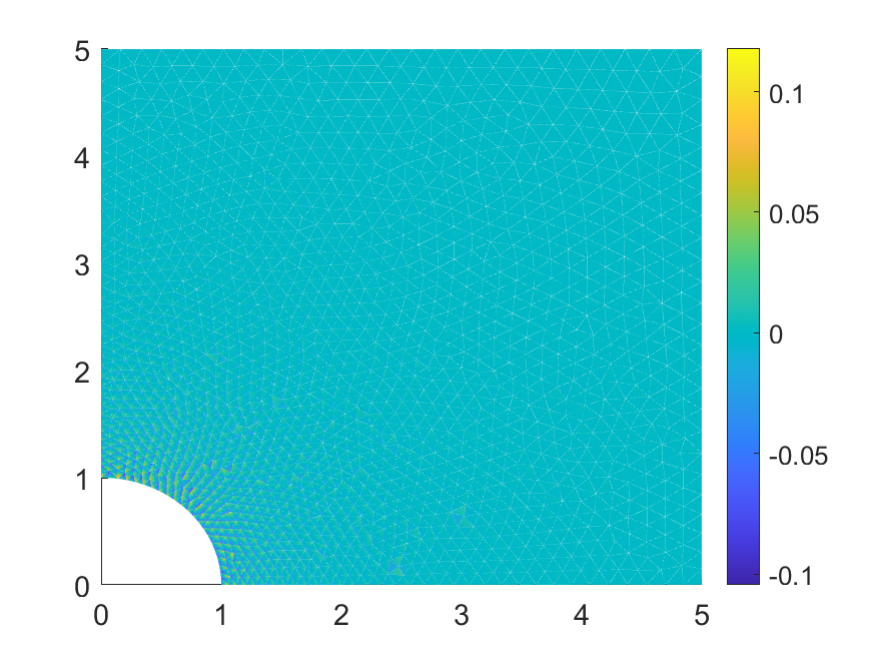}
    \caption{}
    \end{subfigure}
     \hfill
     \begin{subfigure}{0.48\textwidth}
         \centering
         \includegraphics[width=\textwidth]{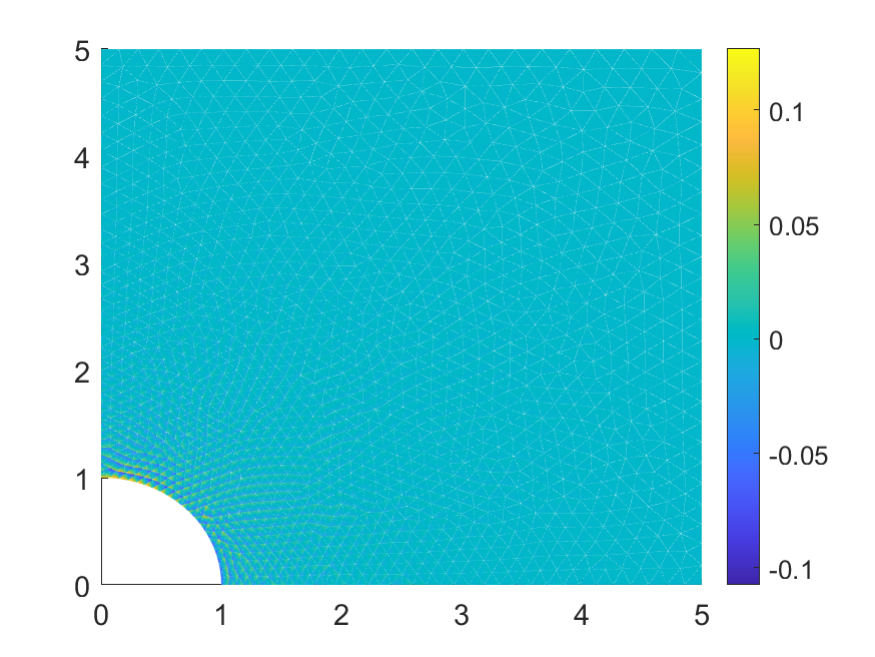}
         \caption{}
     \end{subfigure}
     \vfill
     \begin{subfigure}{0.48\textwidth}
    \centering
    \includegraphics[width=\textwidth]{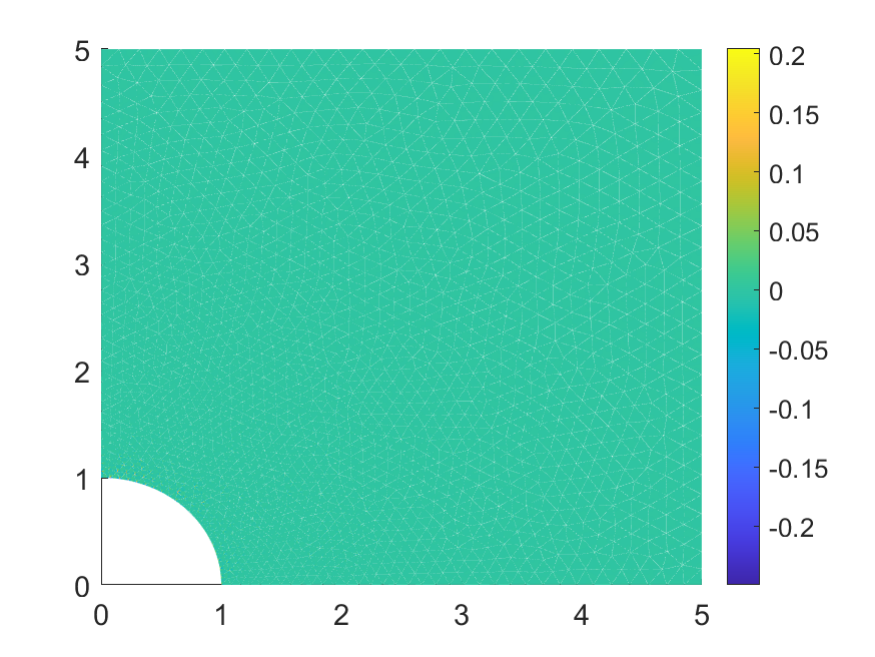}
    \caption{}
    \end{subfigure}
     \hfill
     \begin{subfigure}{0.48\textwidth}
         \centering
         \includegraphics[width=\textwidth]{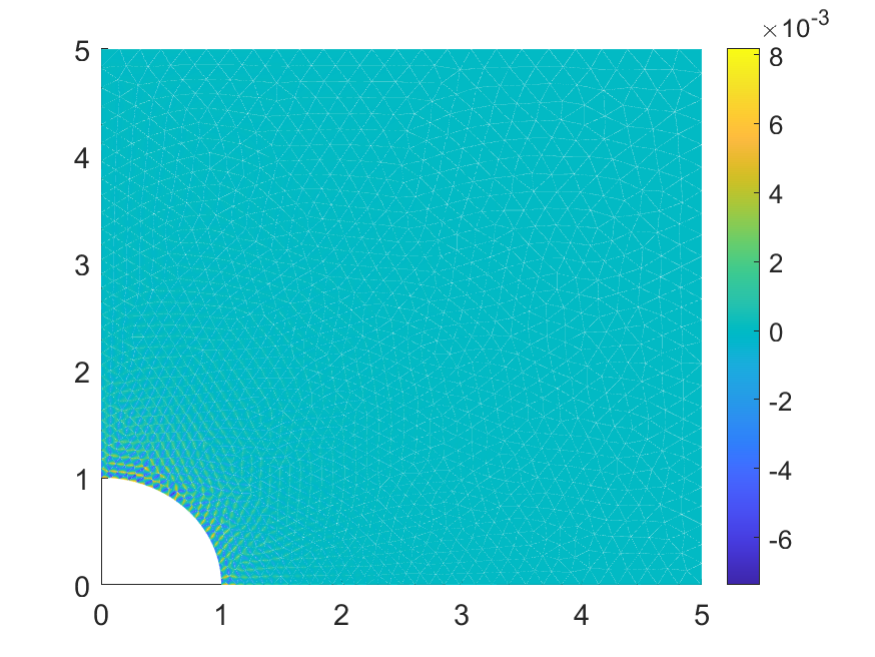}
         \caption{}
     \end{subfigure}
    \caption{Contour plots of the error $\tilde{p}-\tilde{p}_h$ in hydrostatic stress on unstructured meshes (\ac{see~\fref{fig:platemesh_unstructured}}) for the plate with a circular hole problem. (a) CT FEM, (b) B-bar VEM, (c) SH-VEM, and (d) PSH-VEM.}
    \label{fig:plate_hole_contour}
\end{figure}

We consider the problem of an infinite plate with a circular hole of radius $a=1$ inch under uniform tension~\cite{timoshenko1951theory}. The hole is assumed to be traction-free and a far-field load $\sigma_0=1$ psi is applied in the $x$-direction. Since the plate is symmetric, we model the problem as a quarter plate with length $L=5$ inch and apply exact tractions along the boundary. The material has Young's modulus $E_Y= 2\times 10^7$ psi and Poisson's ratio $\nu=0.49995$. The first test of this problem is on unstructured triangular meshes with representative meshes shown in~\fref{fig:platemesh_unstructured}. In~\fref{fig:plate_convergence_unstructured}, we show the convergence results of the four methods that are tested and find that the three methods CT FEM, B-bar VEM, and SH-VEM all have optimal convergence rates. However, the penalty stress-hybrid method has third order superconvergence in the energy seminorm and $L^2$ norm of hydrostatic stress.  In~\fref{fig:plate_hole_contour}, we plot the contours of the pointwise error $\tilde{p}-\tilde{p}_h$ of the hydrostatic stress. The plots reveal that the penalty stress-hybrid method has the smallest pointwise error and the remaining three methods produce similar error distributions. 

\begin{figure}[!bht]
     \centering
     \begin{subfigure}{.32\textwidth}
         \centering
         \includegraphics[width=\textwidth]{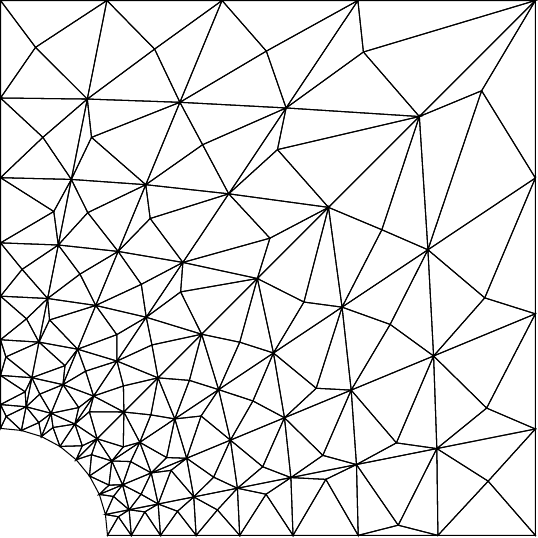}
         \caption{}
     \end{subfigure}
     \hfill
     \begin{subfigure}{.32\textwidth}
         \centering
         \includegraphics[width=\textwidth]{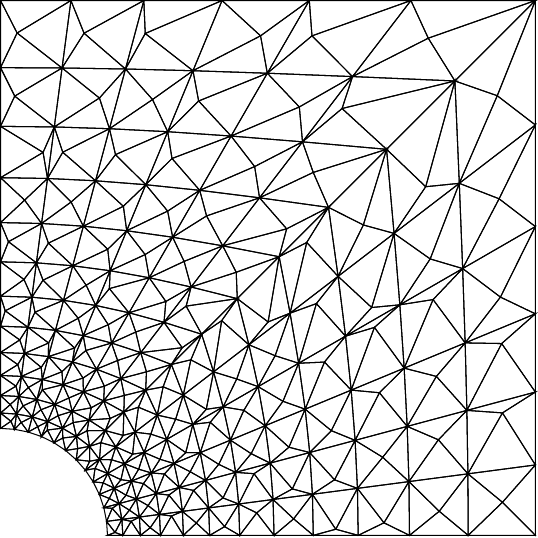}
         \caption{}
     \end{subfigure}
     \hfill
     \begin{subfigure}{.32\textwidth}
         \centering
         \includegraphics[width=\textwidth]{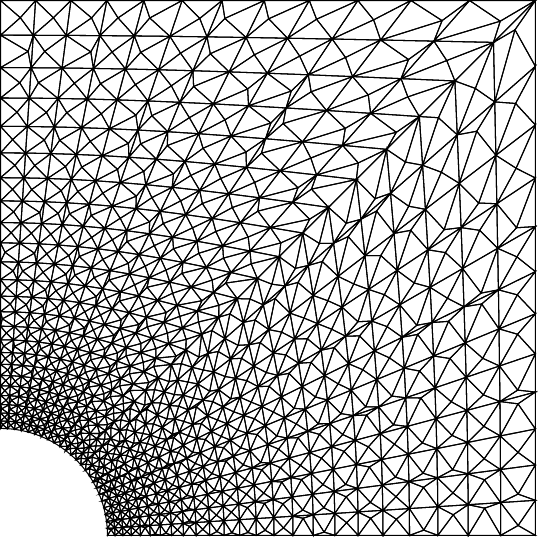}
         \caption{}
     \end{subfigure}
        \caption{Perturbed triangular meshes for the plate with a hole problem. (a) 250 elements, (b) 500 elements, and (c) 2500 elements.  }
        \label{fig:platemesh_perturbed}
\end{figure}

\begin{figure}[!bht]
     \centering
     \begin{subfigure}{0.32\textwidth}
         \centering
         \includegraphics[width=\textwidth]{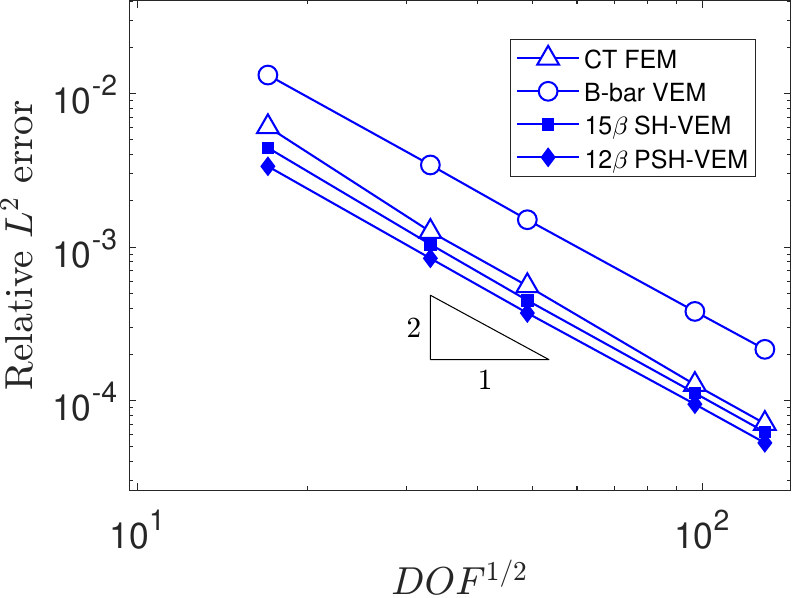}
         \caption{}
     \end{subfigure}
     \hfill
     \begin{subfigure}{0.32\textwidth}
         \centering
         \includegraphics[width=\textwidth]{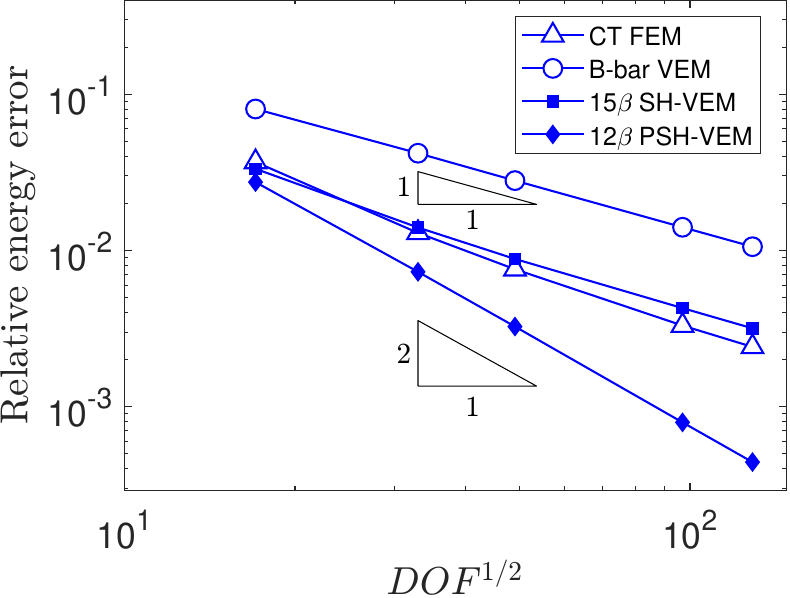}
         \caption{}
     \end{subfigure}
     \hfill
     \begin{subfigure}{0.32\textwidth}
         \centering
         \includegraphics[width=\textwidth]{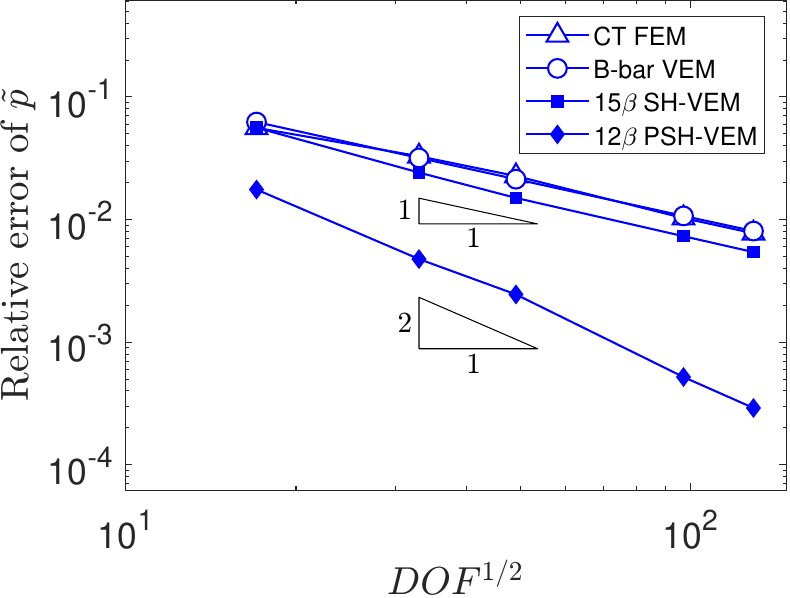}
         \caption{}
     \end{subfigure}
        \caption{Comparison of CT FEM, B-bar VEM, SH-VEM, and PSH-VEM for the plate with a hole problem on perturbed meshes (\ac{see~\fref{fig:platemesh_perturbed}}). (a) $L^2$ error of displacement, (b) energy error, and (c) $L^2$ error of hydrostatic stress. }
        \label{fig:plate_convergence_perturbed}
\end{figure}

\begin{figure}[!bht]
    \centering
    \begin{subfigure}{0.48\textwidth}
    \centering
    \includegraphics[width=\textwidth]{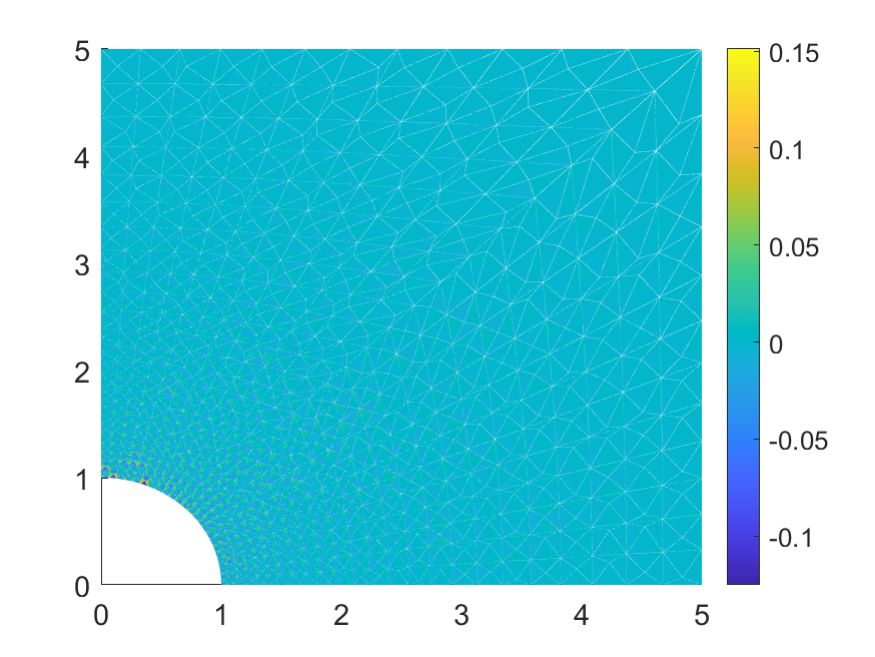}
    \caption{}
    \end{subfigure}
     \hfill
     \begin{subfigure}{0.48\textwidth}
         \centering
         \includegraphics[width=\textwidth]{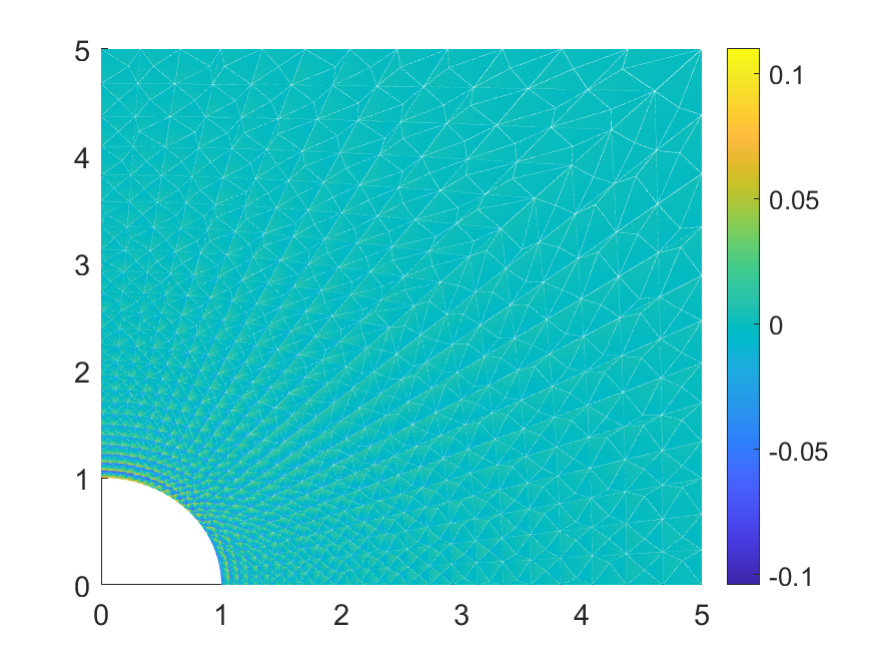}
         \caption{}
     \end{subfigure}
     \vfill
      \begin{subfigure}{0.48\textwidth}
    \centering
    \includegraphics[width=\textwidth]{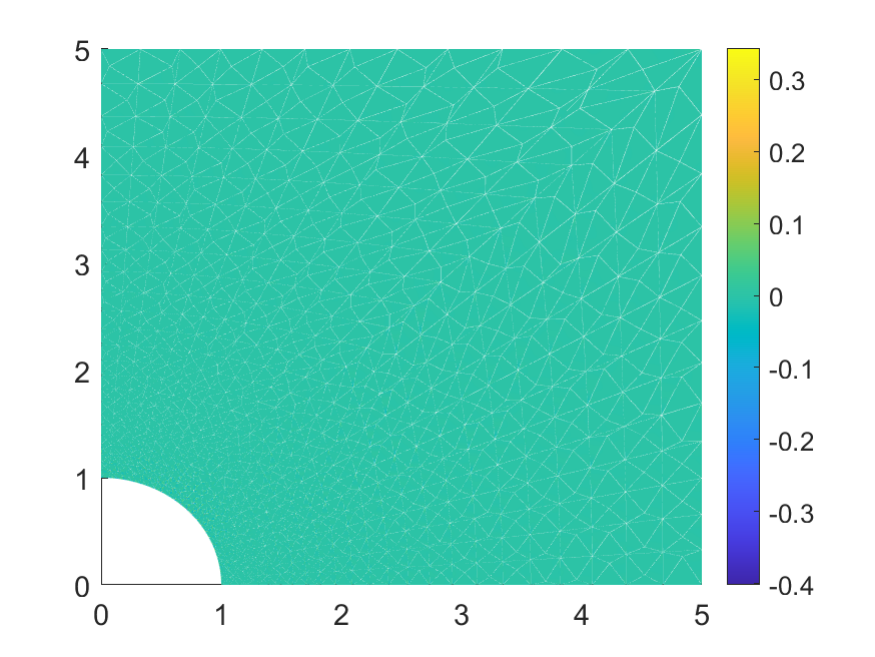}
    \caption{}
    \end{subfigure}
     \hfill
     \begin{subfigure}{0.48\textwidth}
         \centering
         \includegraphics[width=\textwidth]{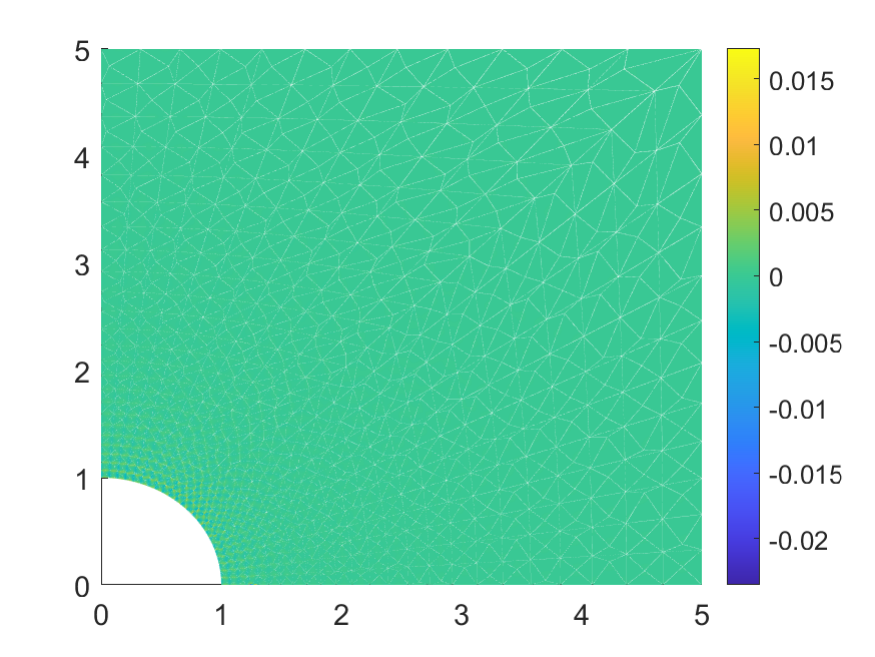}
         \caption{}
     \end{subfigure}
    \caption{Contour plots of the error $\tilde{p}-\tilde{p}_h$ in hydrostatic stress on perturbed meshes (\ac{see~\fref{fig:platemesh_perturbed}}) for the plate with a circular hole problem. (a) CT FEM, (b) B-bar VEM, (c) SH-VEM, and (d) PSH-VEM.}
    \label{fig:plate_hole_perturbed_contour}
\end{figure}

We now solve the problem on perturbed meshes. In~\cite{Ainsworth:2022:cmame}, it is shown that small perturbations of a regular triangular mesh can lead to locking and a reduction in the rate of convergence in the hydrostatic stress error. For this test, we start with a regular quadrilateral mesh and then cut along both diagonals to create four triangles. The point of intersection of the two diagonals are perturbed for each quadrilateral. Representative meshes are shown in~\fref{fig:platemesh_perturbed}. In~\fref{fig:plate_convergence_perturbed}, we plot the convergence of the four methods and find that CT FEM, B-bar VEM, and SH-VEM all still retain optimal rates of convergence. The penalty formulation exhibits second order superconvergence in the energy seminorm and $L^2$ norm of the hydrostatic stress. The contour plots in~\fref{fig:plate_hole_perturbed_contour} show that the four methods all have relatively smooth error distributions of the hydrostatic stress and display no signs of volumetric locking.


\subsection{Hollow cylinder under internal pressure}
\begin{figure}[!bht]
     \centering
     \begin{subfigure}{.32\textwidth}
         \centering
         \includegraphics[width=\textwidth]{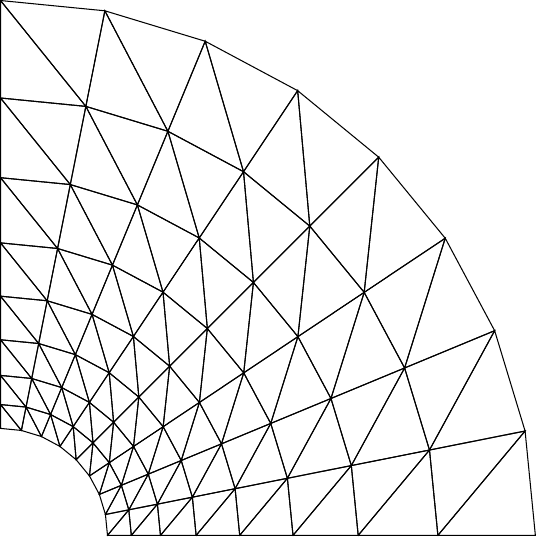}
         \caption{}
     \end{subfigure}
     \hfill
     \begin{subfigure}{.32\textwidth}
         \centering
         \includegraphics[width=\textwidth]{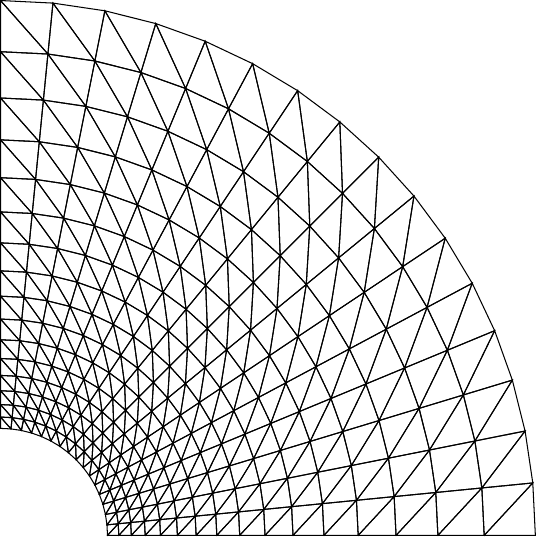}
         \caption{}
     \end{subfigure}
     \hfill
     \begin{subfigure}{.32\textwidth}
         \centering
         \includegraphics[width=\textwidth]{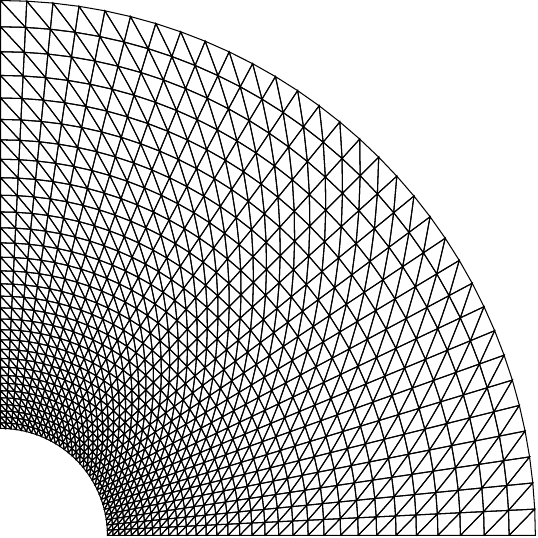}
         \caption{}
     \end{subfigure}
        \caption{Structured triangular meshes for the hollow cylinder problem. (a) 128 elements, (b) 512 elements, and (c) 2048 elements.  }
        \label{fig:cylinder_structured}
\end{figure}

\begin{figure}[!bht]
     \centering
     \begin{subfigure}{0.32\textwidth}
         \centering
         \includegraphics[width=\textwidth]{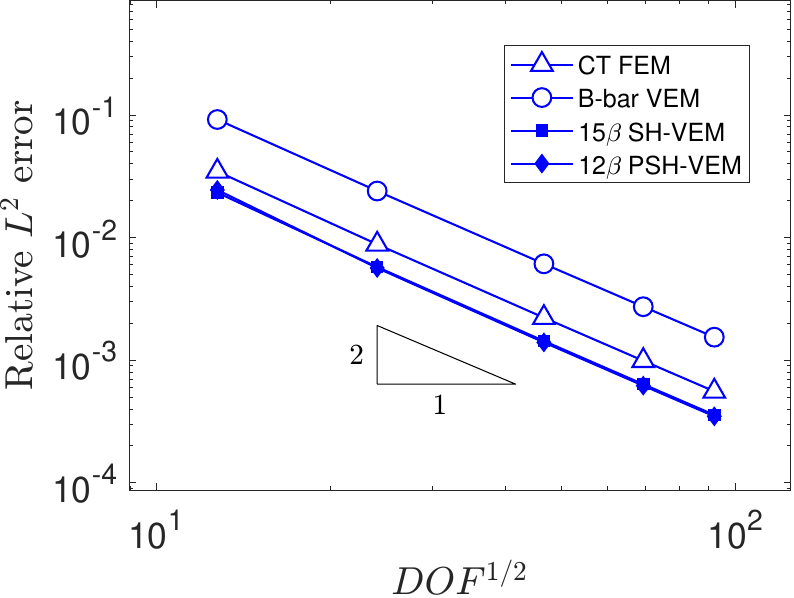}
         \caption{}
     \end{subfigure}
     \hfill
     \begin{subfigure}{0.32\textwidth}
         \centering
         \includegraphics[width=\textwidth]{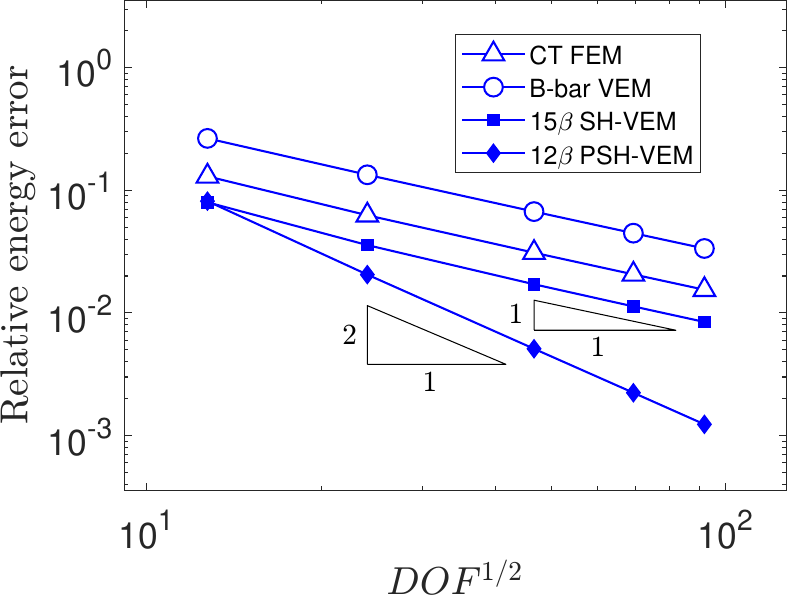}
         \caption{}
     \end{subfigure}
     \hfill
     \begin{subfigure}{0.32\textwidth}
         \centering
         \includegraphics[width=\textwidth]{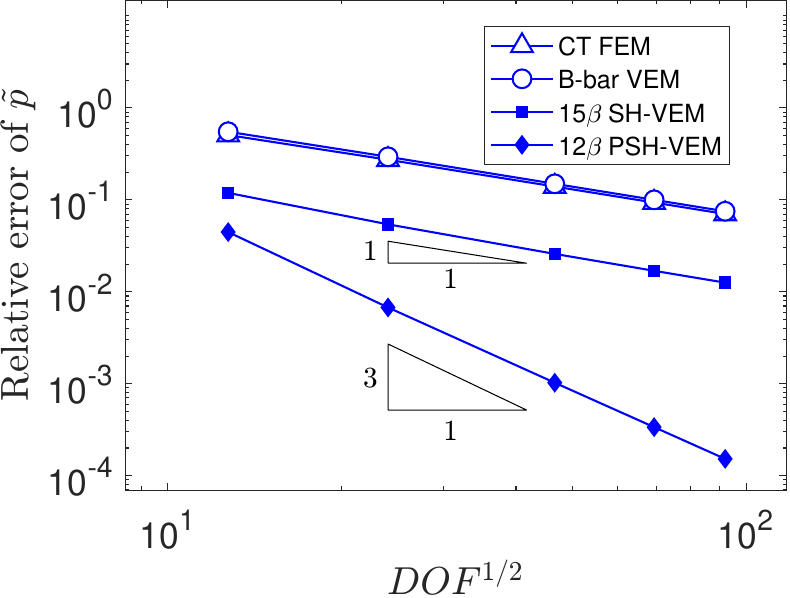}
         \caption{}
     \end{subfigure}
        \caption{Comparison of CT FEM, B-bar VEM, SH-VEM, and PSH-VEM for the pressurized cylinder problem on structured meshes (\ac{see~\fref{fig:cylinder_structured}}). (a) $L^2$ error of displacement, (b) energy error, and (c) $L^2$ error of hydrostatic stress. }
        \label{fig:cylinder_convergence_structured}
\end{figure}

\begin{figure}[!bht]
    \centering
    \begin{subfigure}{0.48\textwidth}
    \centering
    \includegraphics[width=\textwidth]{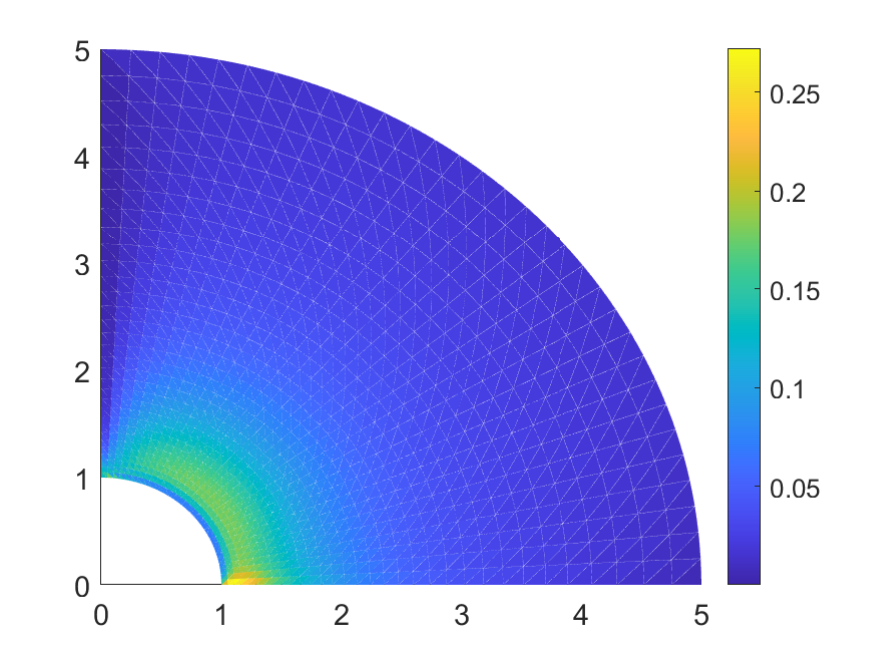}
    \caption{}
    \end{subfigure}
     \hfill
     \begin{subfigure}{0.48\textwidth}
         \centering
         \includegraphics[width=\textwidth]{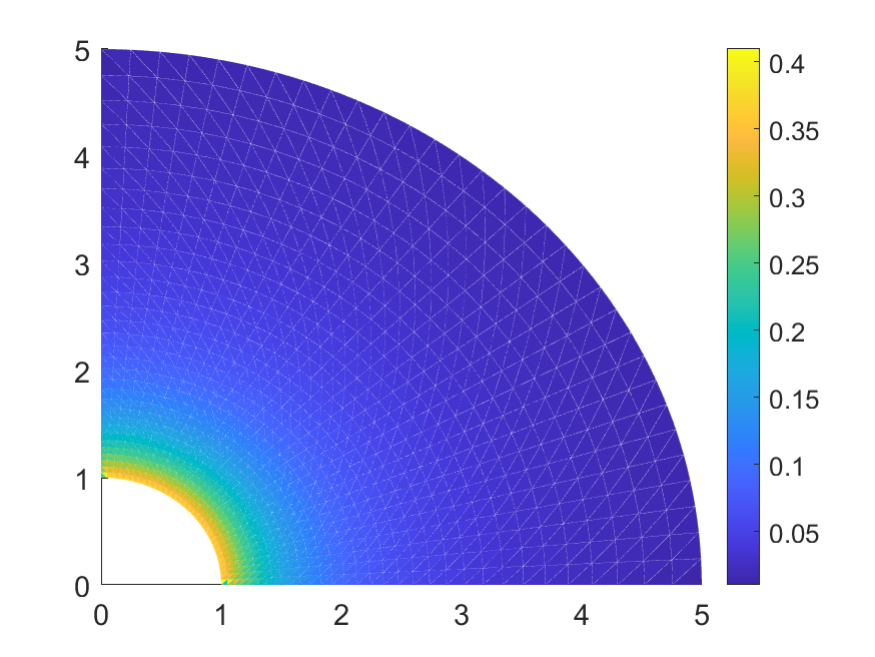}
         \caption{}
     \end{subfigure}
     \vfill
     \begin{subfigure}{0.48\textwidth}
         \centering
         \includegraphics[width=\textwidth]{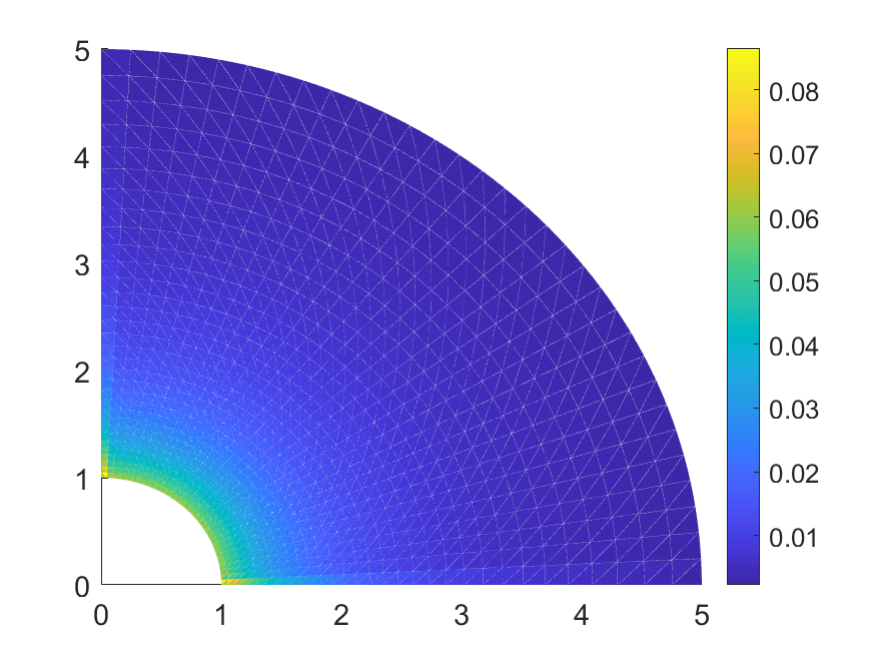}
         \caption{}
     \end{subfigure}
     \hfill
     \begin{subfigure}{0.48\textwidth}
         \centering
         \includegraphics[width=\textwidth]{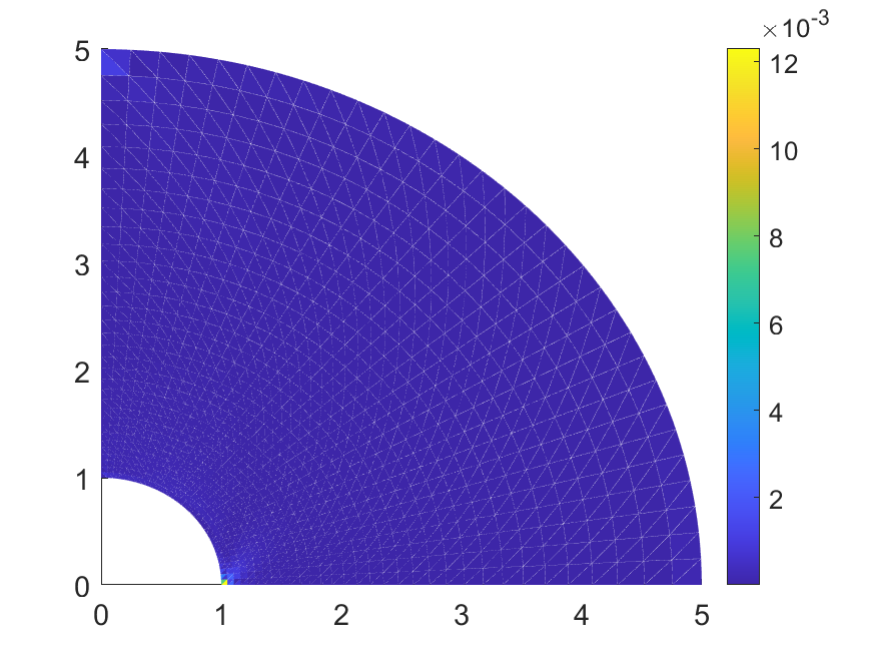}
         \caption{}
     \end{subfigure}
    \caption{Contour plots of the relative error in the hydrostatic stress on structured meshes (\ac{see~\fref{fig:cylinder_structured}}) for the pressurized cylinder problem. The exact hydrostatic stress is 4166.528 psi. (a) CT FEM, (b) B-bar VEM, (c) SH-VEM, and (d) PSH-VEM. }
    \label{fig:cylinder_structured_pressure_contour}
\end{figure}

Next, we consider a hollow cylinder subject to internal pressure~\cite{timoshenko1951theory}. By symmetry, the problem is modeled by a quarter cylinder with an outer radius $b=5$ inch and an inner radius $a=1$ inch. The outer radius of the cylinder is assumed to be traction-free, while a uniform pressure of $p=10^5$ psi is applied on the inner radius. The material properties are $E_Y = 2\times 10^5$ psi and $\nu = 0.49995$. For this problem, the hydrostatic stress field is constant and we found that using nonconstant stress functions in the SH-VEM resulted in larger errors around the element corners. Therefore, we use an element averaged hydrostatic stress approximation for SH-VEM. We first solve the problem on structured triangular meshes; a few representative meshes are shown in~\fref{fig:cylinder_structured}. In~\fref{fig:cylinder_convergence_structured}, we show the rates of convergence in three error norms. The three methods CT FEM, B-bar VEM, and SH-VEM produce optimal convergence in the displacement $L^2$ norm, $L^2$ norm of hydrostatic stress as well as in the energy seminorm; while the PSH-VEM has second order superconvergence in the energy seminorm and third order in the $L^2$ norm of the hydrostatic stress. In~\fref{fig:cylinder_structured_pressure_contour}, we plot the contours of the relative error in the hydrostatic stress field for the four methods. The plots show that the maximum errors concentrate along the inner radius and improves when away from the boundary. The SH-VEM and PSH-VEM produce the smallest relative errors, with a maximum of around $8$ and $1.2$ percent, respectively.      

\begin{figure}[!bht]
     \centering
     \begin{subfigure}{.32\textwidth}
         \centering
         \includegraphics[width=\textwidth]{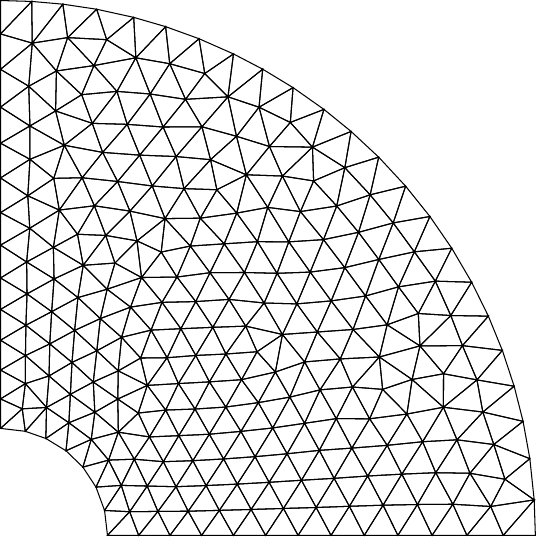}
         \caption{}
     \end{subfigure}
     \hfill
     \begin{subfigure}{.32\textwidth}
         \centering
         \includegraphics[width=\textwidth]{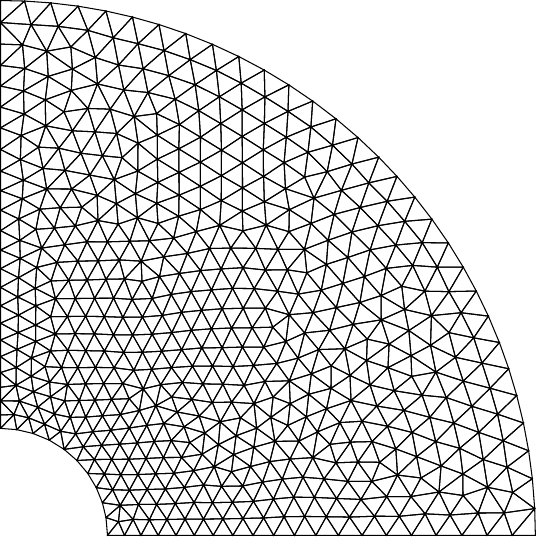}
         \caption{}
     \end{subfigure}
     \hfill
     \begin{subfigure}{.32\textwidth}
         \centering
         \includegraphics[width=\textwidth]{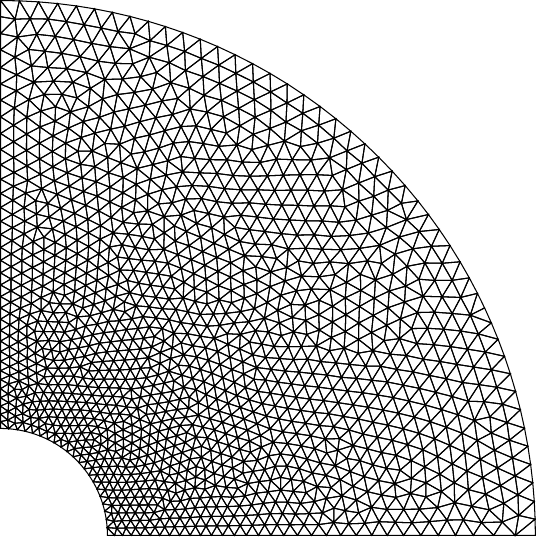}
         \caption{}
     \end{subfigure}
        \caption{Unstructured triangular meshes for the hollow cylinder problem. (a) 500 elements, (b) 1000 elements, and (c) 2500 elements.  }
        \label{fig:cylinder_unstructured}
\end{figure}

\begin{figure}[!bht]
     \centering
     \begin{subfigure}{0.32\textwidth}
         \centering
         \includegraphics[width=\textwidth]{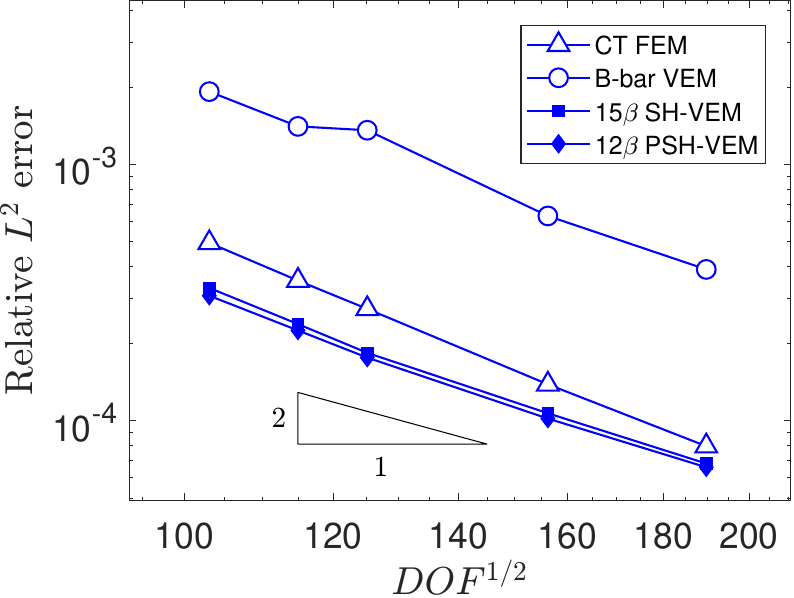}
         \caption{}
     \end{subfigure}
     \hfill
     \begin{subfigure}{0.32\textwidth}
         \centering
         \includegraphics[width=\textwidth]{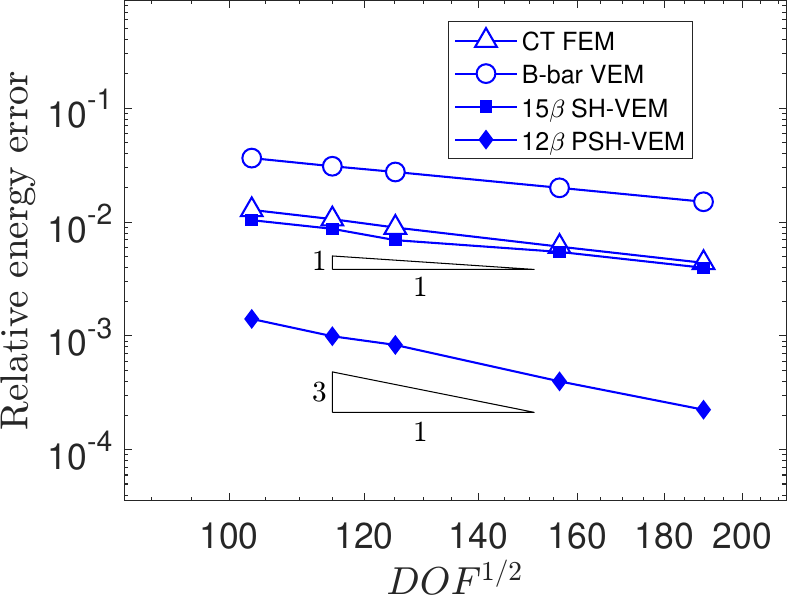}
         \caption{}
     \end{subfigure}
     \hfill
     \begin{subfigure}{0.32\textwidth}
         \centering
         \includegraphics[width=\textwidth]{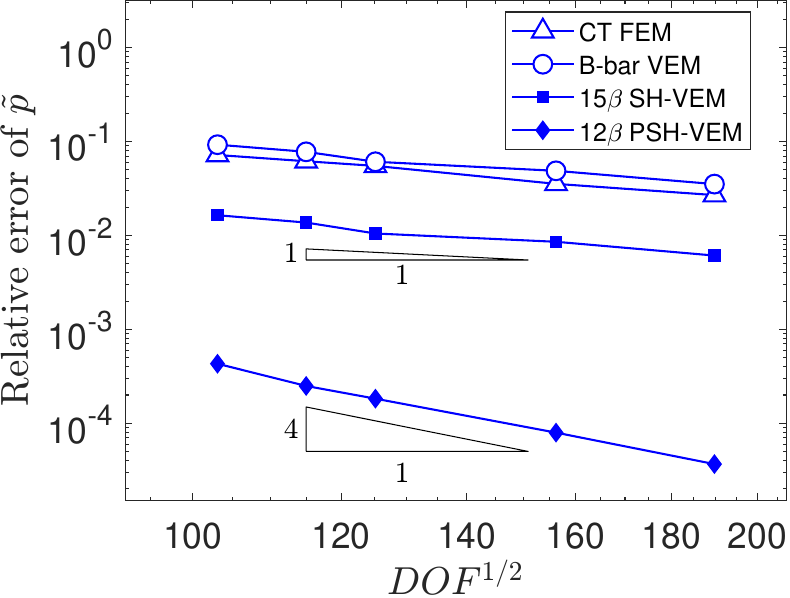}
         \caption{}
     \end{subfigure}
        \caption{Comparison of CT FEM, B-bar VEM, SH-VEM, and PSH-VEM for the pressurized cylinder problem on unstructured meshes (\ac{see~\fref{fig:cylinder_unstructured}}). (a) $L^2$ error of displacement, (b) energy error, and (c) $L^2$ error of hydrostatic stress. }
        \label{fig:cylinder_convergence_unstructured}
\end{figure}

\begin{figure}[!bht]
    \centering
    \begin{subfigure}{0.48\textwidth}
    \centering
    \includegraphics[width=\textwidth]{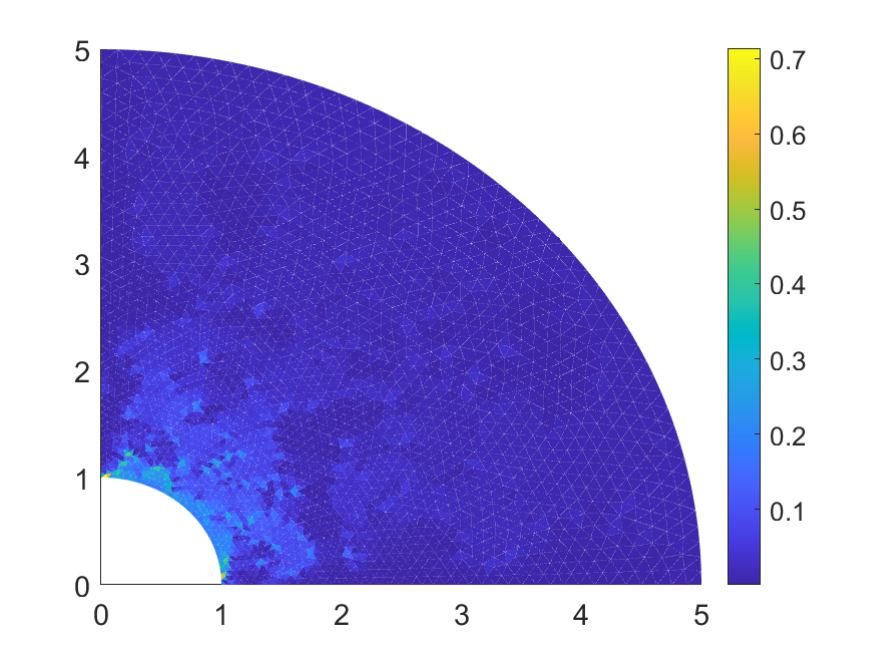}
    \caption{}
    \end{subfigure}
     \hfill
     \begin{subfigure}{0.48\textwidth}
         \centering
         \includegraphics[width=\textwidth]{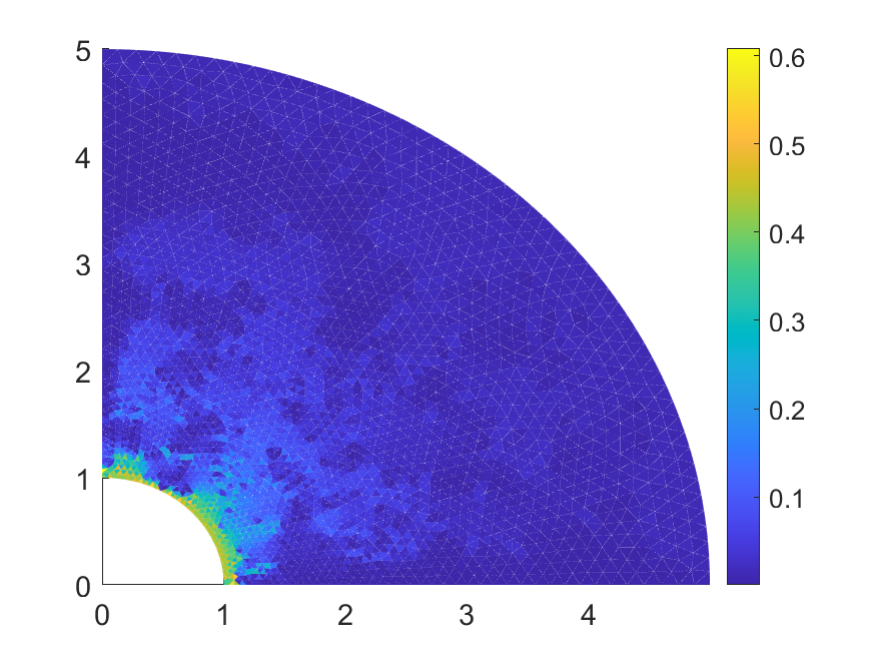}
         \caption{}
     \end{subfigure}
      \vfill
     \begin{subfigure}{0.48\textwidth}
         \centering
         \includegraphics[width=\textwidth]{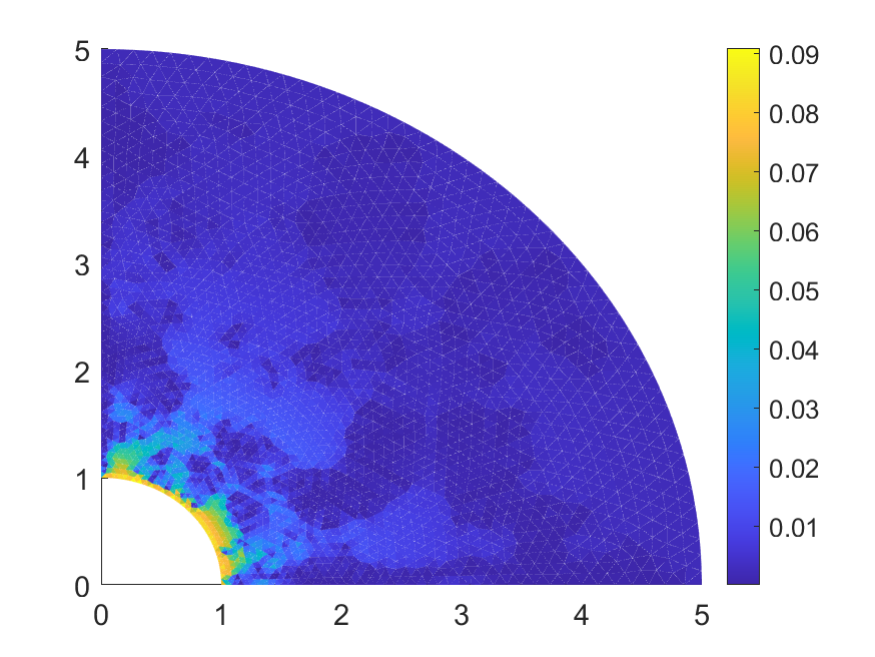}
         \caption{}
     \end{subfigure}
     \hfill
     \begin{subfigure}{0.48\textwidth}
         \centering
         \includegraphics[width=\textwidth]{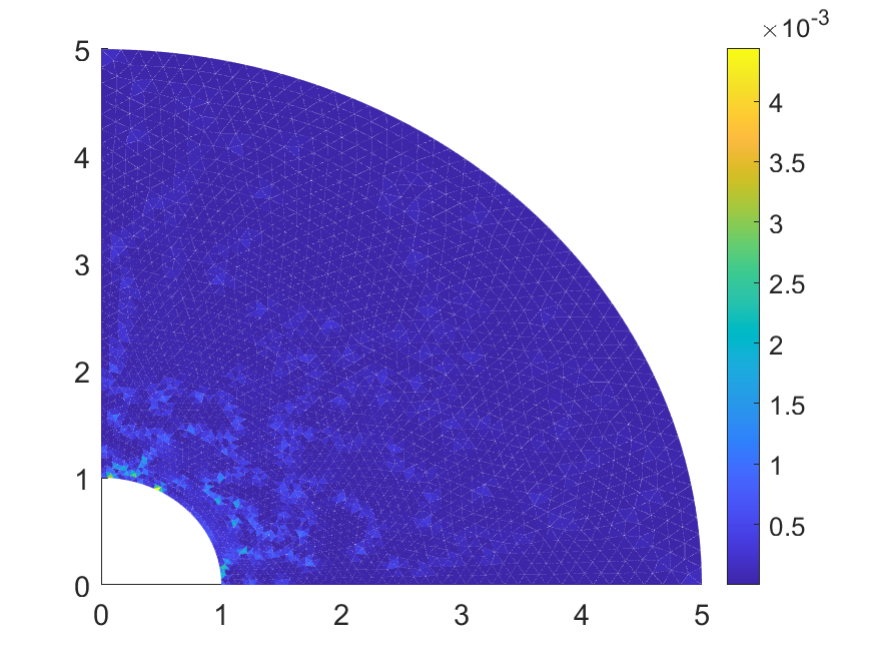}
         \caption{}
     \end{subfigure}
    \caption{Contour plots of the relative error in the hydrostatic stress on unstructured meshes (\ac{see~\fref{fig:cylinder_unstructured}}) for the pressurized cylinder problem. The exact hydrostatic stress is 4166.528 psi. (a) CT FEM, (b) B-bar VEM, (c) SH-VEM, and (d) PSH-VEM.}
    \label{fig:cylinder_unstructured_pressure_contour}
\end{figure}

We now test the hollow cylinder problem on an unstuctured triangular mesh; a few representative meshes are presented in~\fref{fig:cylinder_unstructured}. In~\fref{fig:cylinder_convergence_unstructured}, the convergence results are given and again show that CT FEM, B-bar VEM, and SH-VEM deliver optimal rates. The penalty approach attains third and fourth order superconvergence in the energy seminorm and hydrostatic stress $L^2$ norm, respectively. The contour plots in~\fref{fig:cylinder_unstructured_pressure_contour} show that for CT FEM, B-bar VEM, and SH-VEM, the errors are concentrated near the inner radius. The largest error from SH-VEM is $9$ percent, while both CT FEM and B-bar VEM produce much larger errors of $70$ and $60$ percent, respectively. For PSH-VEM, the errors are much smaller, with a maximum relative error of $0.4$ percent.      

\begin{figure}[!htb]
     \centering
     \begin{subfigure}{0.32\textwidth}
         \centering
         \includegraphics[width=\textwidth]{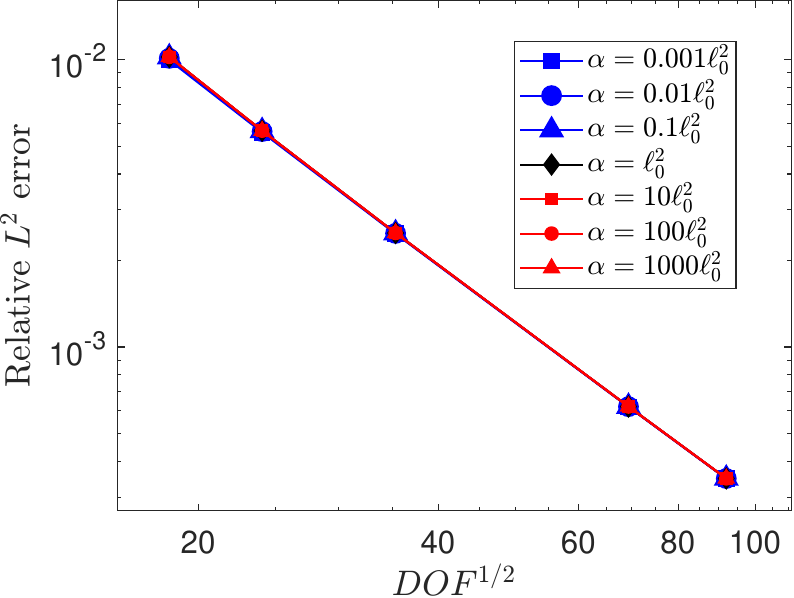}
         \caption{}
     \end{subfigure}
     \hfill
     \begin{subfigure}{0.32\textwidth}
         \centering
         \includegraphics[width=\textwidth]{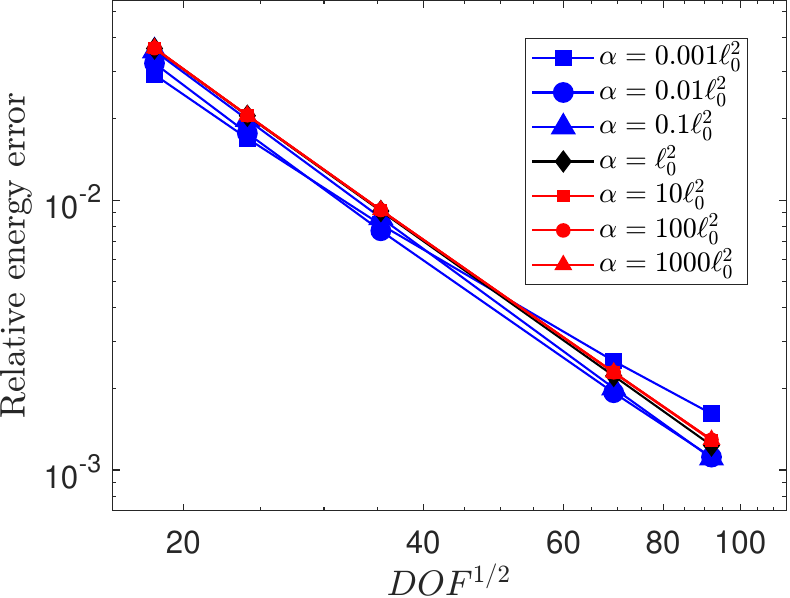}
         \caption{}
     \end{subfigure}
     \hfill
     \begin{subfigure}{0.32\textwidth}
         \centering
         \includegraphics[width=\textwidth]{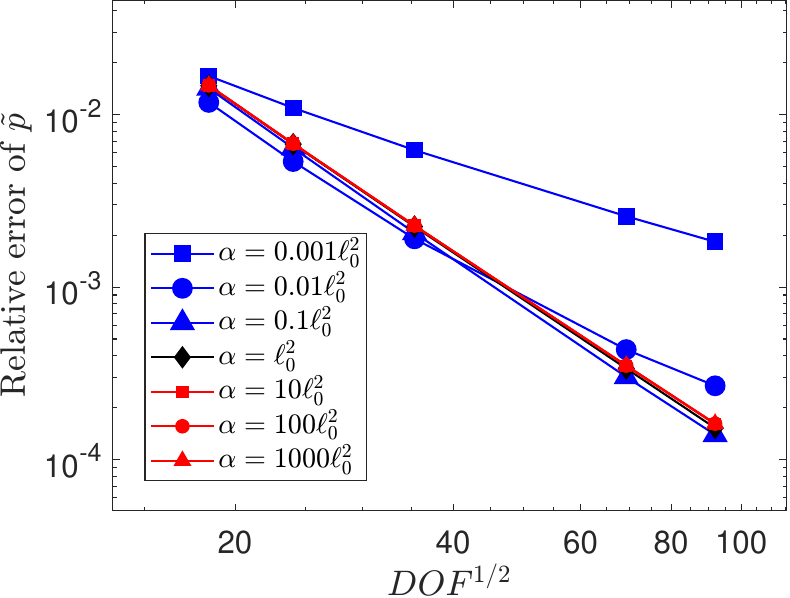}
         \caption{}
     \end{subfigure}
        \caption{Comparison of the convergence of PSH-VEM for different choices of the penalty parameter $\alpha$ for the pressurized cylinder problem on structured meshes (\ac{see~\fref{fig:cylinder_structured}}). (a) $L^2$ error of displacement, (b) energy error, and (c) $L^2$ error of hydrostatic stress. }
        \label{fig:cylinder_penalty_vary}
\end{figure}

Next, we examine the effects of the penalty parameter on the convergence rates of the pressurized cylinder problem. For simplicity, the material properties are set to be $E_Y = 1\times 10^4$ psi and $\nu=0.49995$. The problem is solved on structured triangular meshes (see~\fref{fig:cylinder_structured}). The initial penalty parameter is $\alpha = \frac{10^4\ell_0^2}{E_Y} = \ell_0^2$, so we test $\alpha$ in the range of $10^{-3}\ell_0^2$ to $10^3\ell_0^2$. In~\fref{fig:cylinder_penalty_vary}, we show the rates of convergence in the three error norms for different values of the penalty parameter $\alpha$. The plots show that varying $\alpha$ did not affect the convergence of the displacement errors; however, the convergence of the energy seminorm is slightly affected and the rate of convergence in the $L^2$ norm of the hydrostatic stress is significantly reduced as $\alpha$ is decreased.    

\subsection{Manufactured problem}
\begin{figure}[!bht]
     \centering
     \begin{subfigure}{.32\textwidth}
         \centering
         \includegraphics[width=\textwidth]{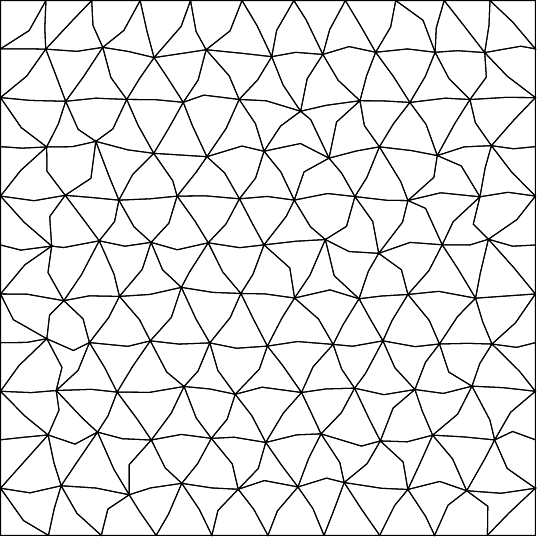}
         \caption{}
     \end{subfigure}
     \hfill
     \begin{subfigure}{.32\textwidth}
         \centering
         \includegraphics[width=\textwidth]{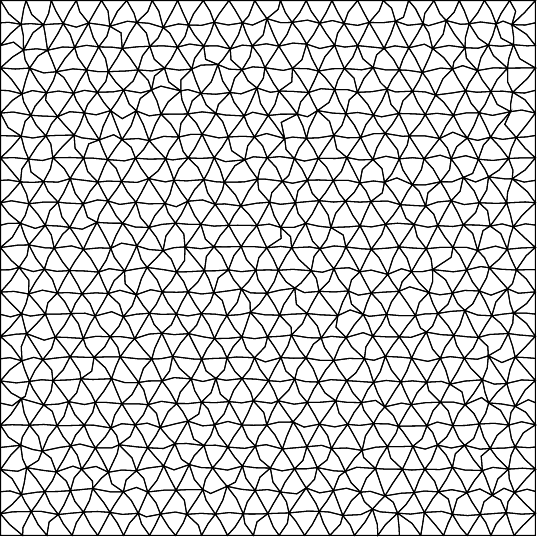}
         \caption{}
     \end{subfigure}
     \hfill
     \begin{subfigure}{.32\textwidth}
         \centering
         \includegraphics[width=\textwidth]{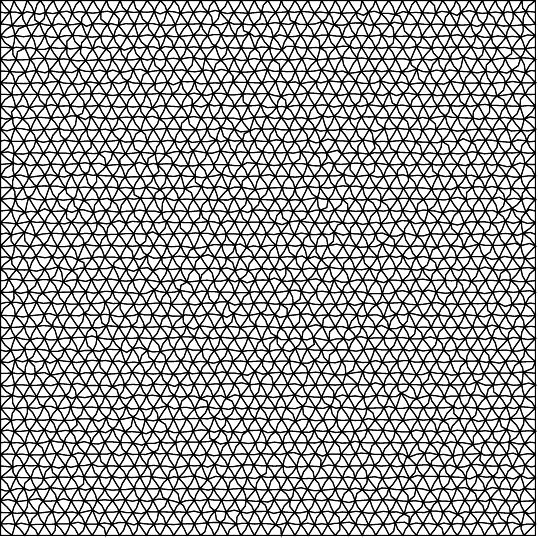}
         \caption{}
     \end{subfigure}
        \caption{Perturbed triangular meshes for the manufactured problem. (a) 200 elements, (b) 1000 elements, and (c) 3600 elements.  }
        \label{fig:manufactured_perturbed}
\end{figure}

\begin{figure}[!bht]
     \centering
     \begin{subfigure}{0.32\textwidth}
         \centering
         \includegraphics[width=\textwidth]{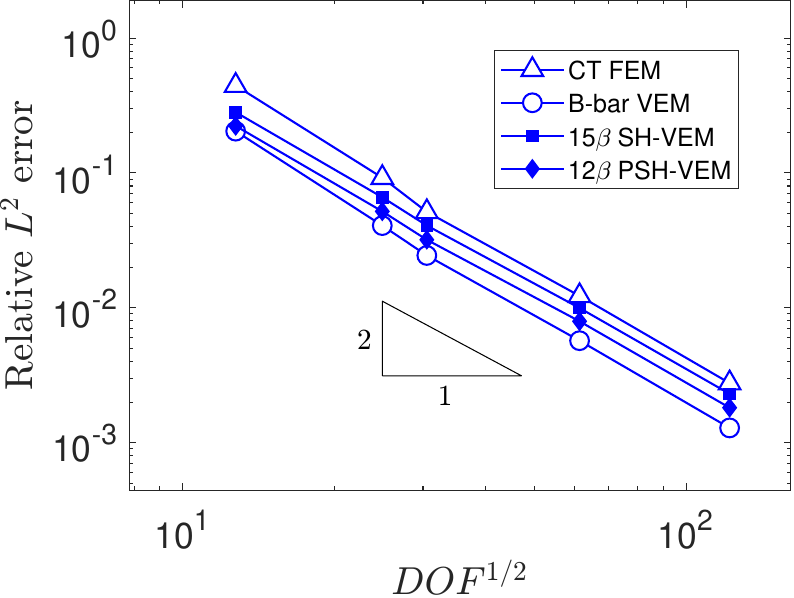}
         \caption{}
     \end{subfigure}
     \hfill
     \begin{subfigure}{0.32\textwidth}
         \centering
         \includegraphics[width=\textwidth]{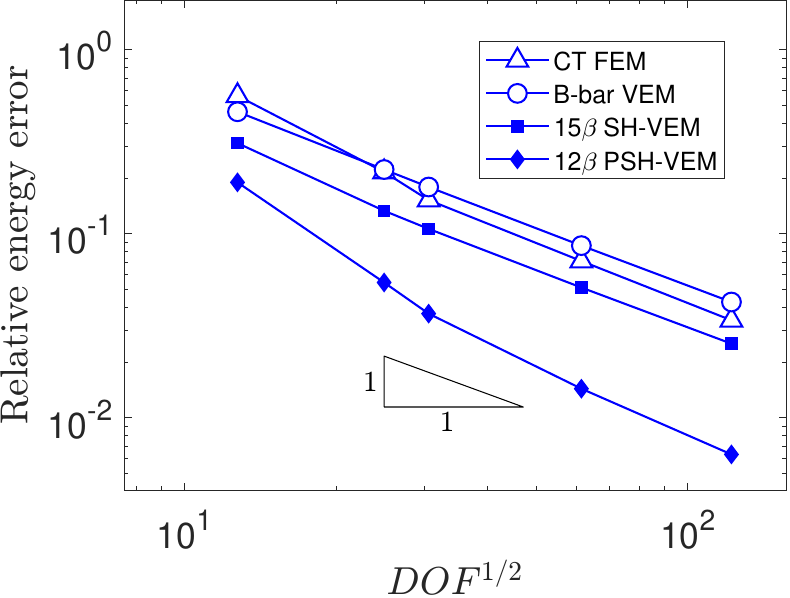}
         \caption{}
     \end{subfigure}
     \hfill
     \begin{subfigure}{0.32\textwidth}
         \centering
         \includegraphics[width=\textwidth]{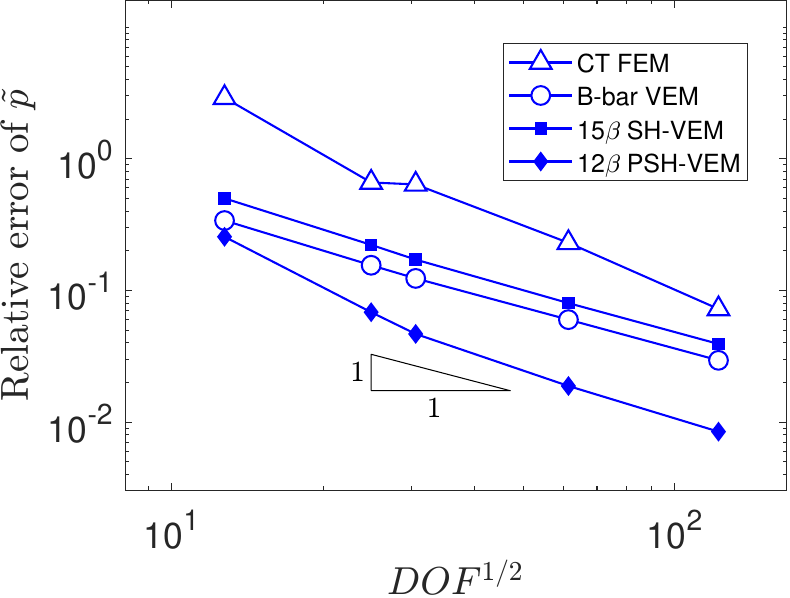}
         \caption{}
     \end{subfigure}
        \caption{Comparison of CT FEM, B-bar VEM, SH-VEM, and PSH-VEM for the manufactured problem on perturbed meshes (\ac{see~\fref{fig:manufactured_perturbed}}). (a) $L^2$ error of displacement, (b) energy error, and (c) $L^2$ error of hydrostatic stress. }
        \label{fig:manufactured_convergence_perturbed}
\end{figure}

\begin{figure}[!bht]
    \centering
    \begin{subfigure}{0.48\textwidth}
    \centering
    \includegraphics[width=\textwidth]{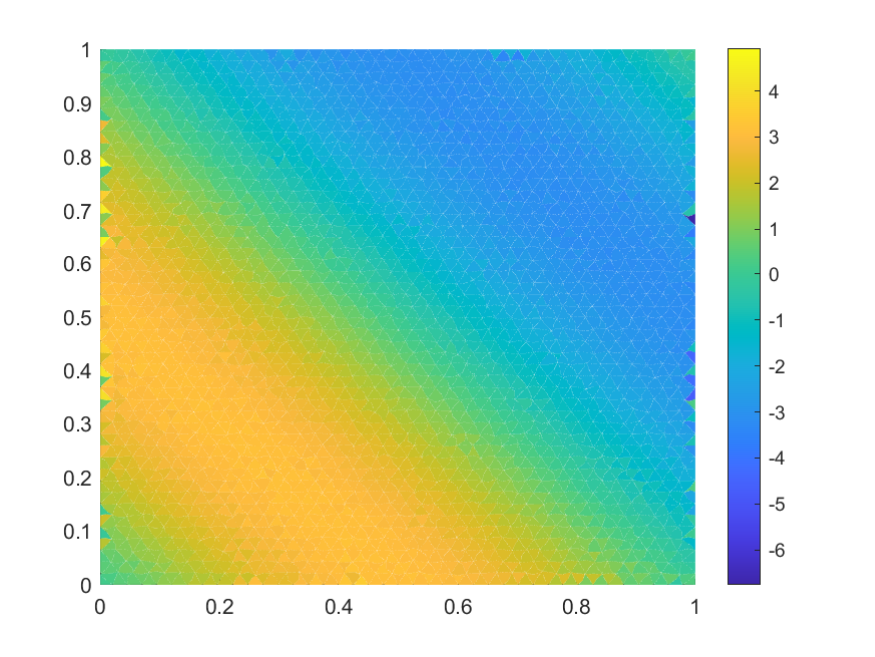}
    \caption{}
    \end{subfigure}
     \hfill
     \begin{subfigure}{0.48\textwidth}
         \centering
         \includegraphics[width=\textwidth]{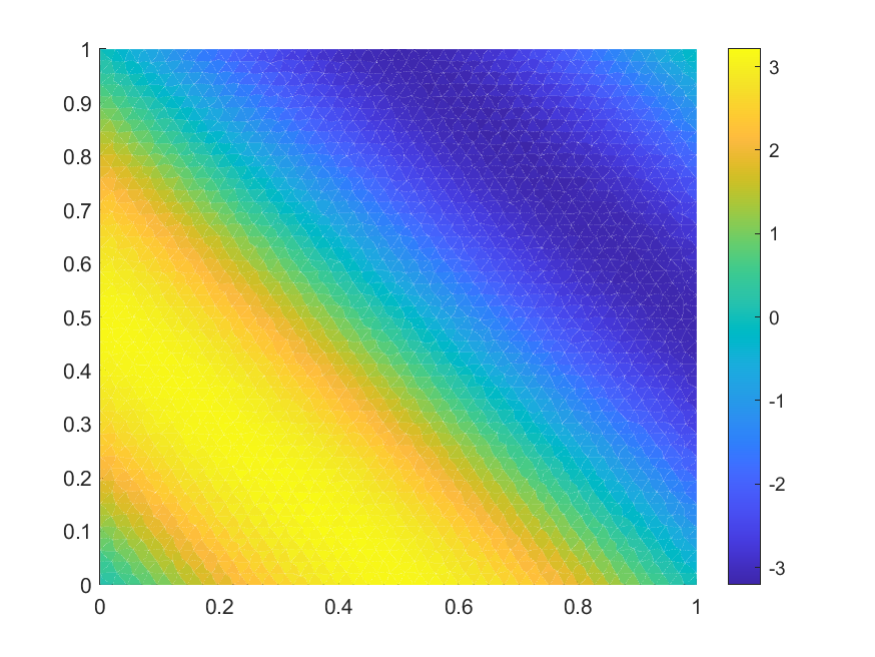}
         \caption{}
     \end{subfigure}
     \vfill
         \begin{subfigure}{0.48\textwidth}
    \centering
    \includegraphics[width=\textwidth]{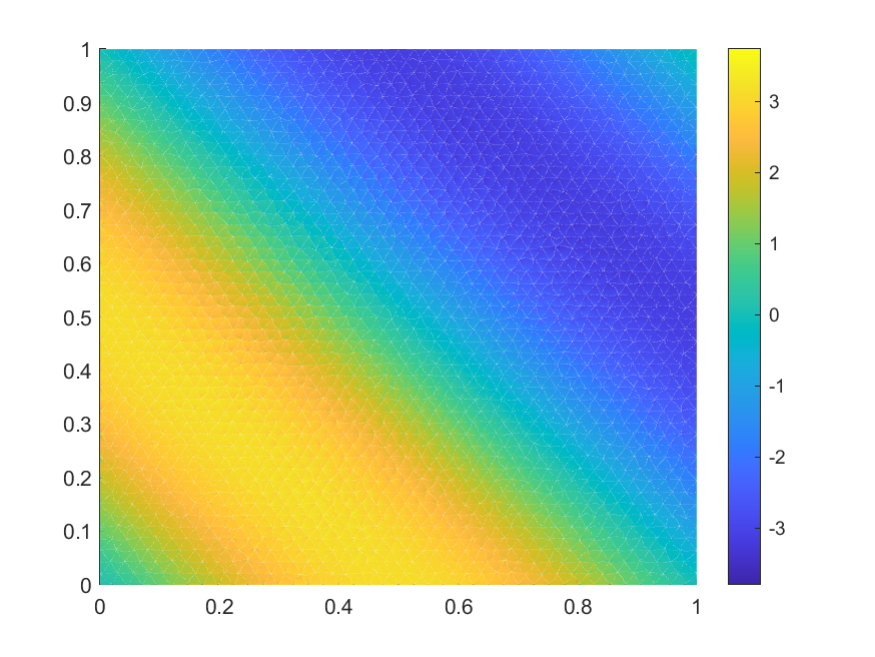}
    \caption{}
    \end{subfigure}
     \hfill
     \begin{subfigure}{0.48\textwidth}
         \centering
         \includegraphics[width=\textwidth]{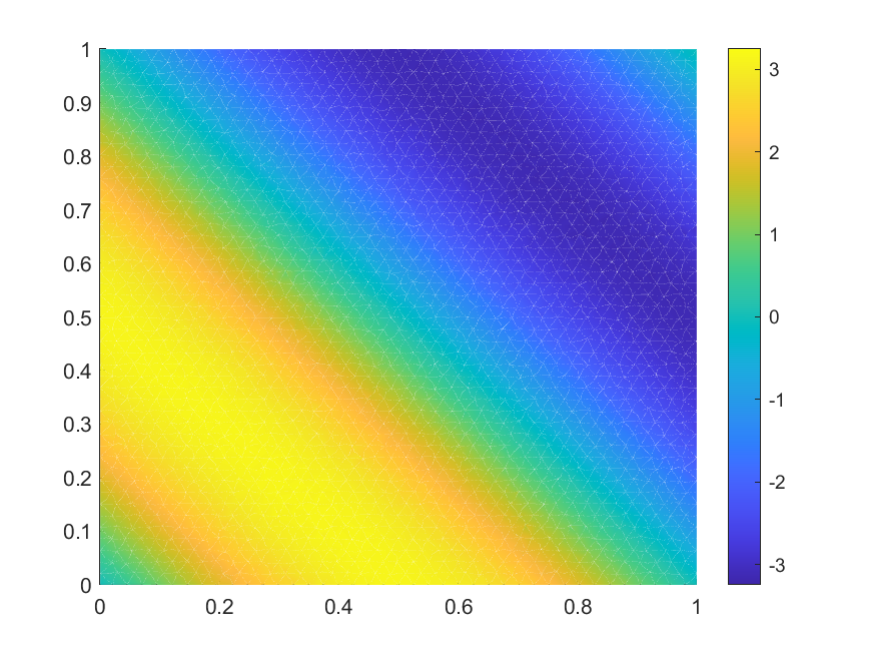}
         \caption{}
     \end{subfigure}
    \caption{Contour plots of the hydrostatic stress on perturbed meshes (\ac{see~\fref{fig:manufactured_perturbed}}) for the manufactured problem. (a) CT FEM, (b) B-bar VEM, (c) SH-VEM, and (d) PSH-VEM.}
    \label{fig:manufactured_perturbed_contour}
\end{figure}

The basis functions derived from the Airy stress functions only satisfy the equilibrium conditions without a body force. Therefore, we test the convergence on a problem with a nonzero body force. We consider a manufactured problem given in~\cite{Shilt:2020:SON} with the exact solution and loading given by:
{\allowdisplaybreaks
\begin{align*}
    u(\vx) &= \sin2\pi y(\cos2\pi x-1) + \frac{1}{1+\lambda}\sin\pi x\sin\pi y, \\
    v(\vx) &= \sin 2\pi x(1-\cos 2\pi y) + \frac{1}{1+\lambda}\sin\pi x\sin\pi y, \\
    \vm{b}(\vx) &= -\pi^2\begin{Bmatrix}
        \left(\frac{\lambda+\mu}{\lambda+1}\right)\cos \pi (x+y) -\mu (8\cos 2\pi x \sin 2\pi y -4\sin 2\pi y +\frac{2}{\lambda+1}\sin \pi x \sin \pi y)  \\ 
        \left(\frac{\lambda+\mu}{\lambda+1}\right)\cos \pi (x+y) -\mu (-8\cos 2\pi y \sin 2\pi x +4\sin 2\pi x +\frac{2}{\lambda+1}\sin \pi x \sin \pi y)
    \end{Bmatrix},
\end{align*}}
where $\lambda$ and $\mu$ are the first and second Lam\'e parameters given by
\begin{equation*}
    \lambda = \frac{E_Y \nu }{(1+\nu)(1-2\nu)}, \quad \mu = \frac{E_Y}{2(1+\nu)}. 
\end{equation*}
For this problem, the material properties are set to $E_Y = 1$ psi and $\nu=0.49995$.
The first set of meshes we use is a perturbed unstructured triangular mesh. A few sample meshes are shown in~\fref{fig:manufactured_perturbed}. The plots in~\fref{fig:manufactured_convergence_perturbed} show that the four methods converge optimally in displacement $L^2$ norm, energy seminorm, and hydrostatic stress $L^2$ norm. The PSH-VEM has smaller errors in both the energy and hydrostatic stress, but does not reach second order superconvergence. The contour plots in~\fref{fig:manufactured_perturbed_contour} reveal that large errors in CT FEM appear along the boundary, while the other three methods have smooth hydrostatic stress fields. 

\begin{figure}[!bht]
     \centering
     \begin{subfigure}{.32\textwidth}
         \centering
         \includegraphics[width=\textwidth]{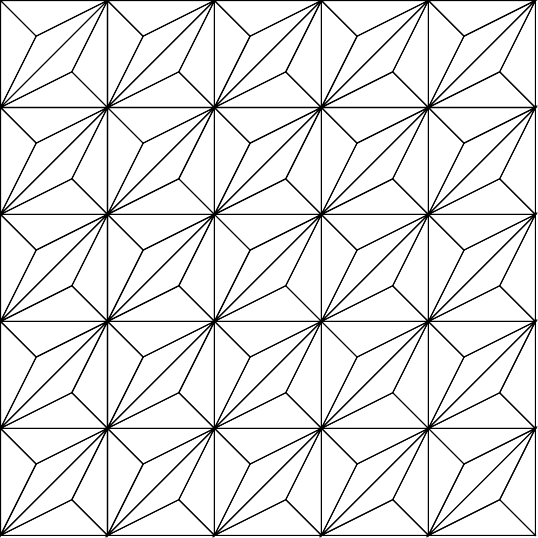}
         \caption{}
     \end{subfigure}
     \hfill
     \begin{subfigure}{.32\textwidth}
         \centering
         \includegraphics[width=\textwidth]{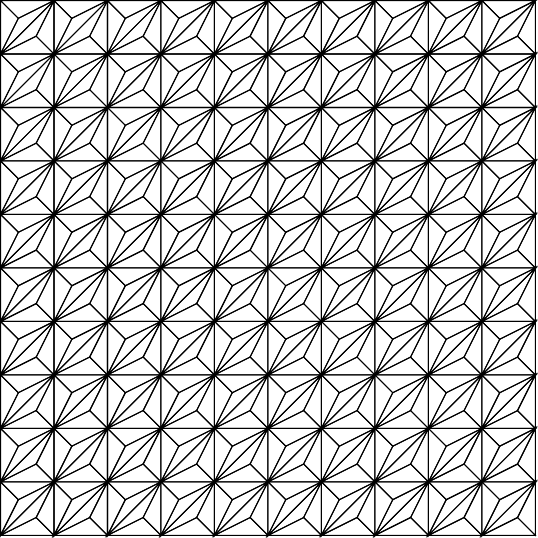}
         \caption{}
     \end{subfigure}
     \hfill
     \begin{subfigure}{.32\textwidth}
         \centering
         \includegraphics[width=\textwidth]{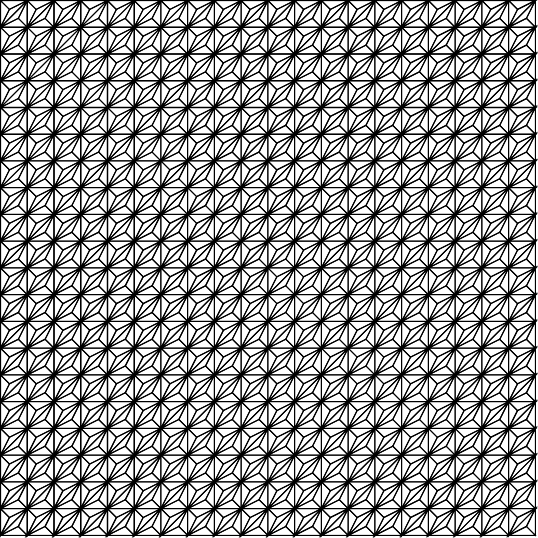}
         \caption{}
     \end{subfigure}
        \caption{Structured triangular meshes for the manufactured problem. (a) 150 elements, (b) 600 elements, and (c) 2400 elements.  }
        \label{fig:manufactured_structured}
\end{figure}

\begin{figure}[!bht]
     \centering
     \begin{subfigure}{0.32\textwidth}
         \centering
         \includegraphics[width=\textwidth]{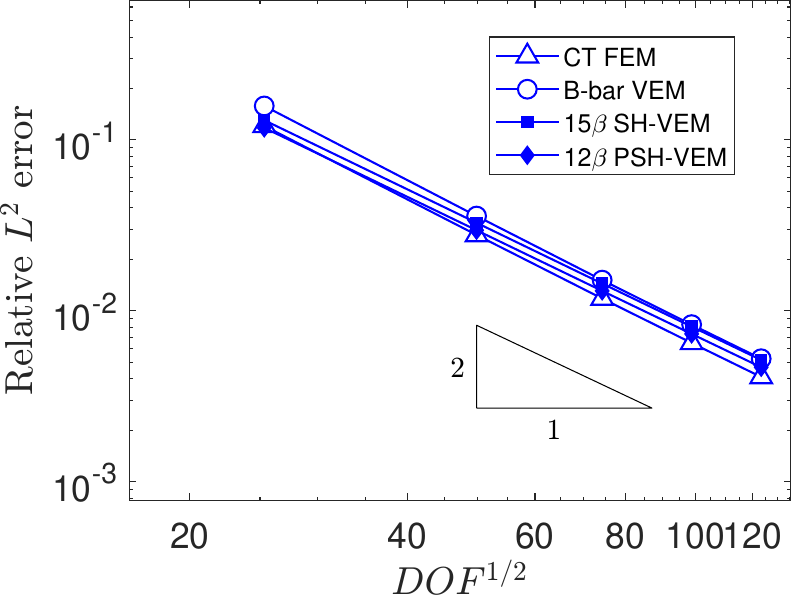}
         \caption{}
     \end{subfigure}
     \hfill
     \begin{subfigure}{0.32\textwidth}
         \centering
         \includegraphics[width=\textwidth]{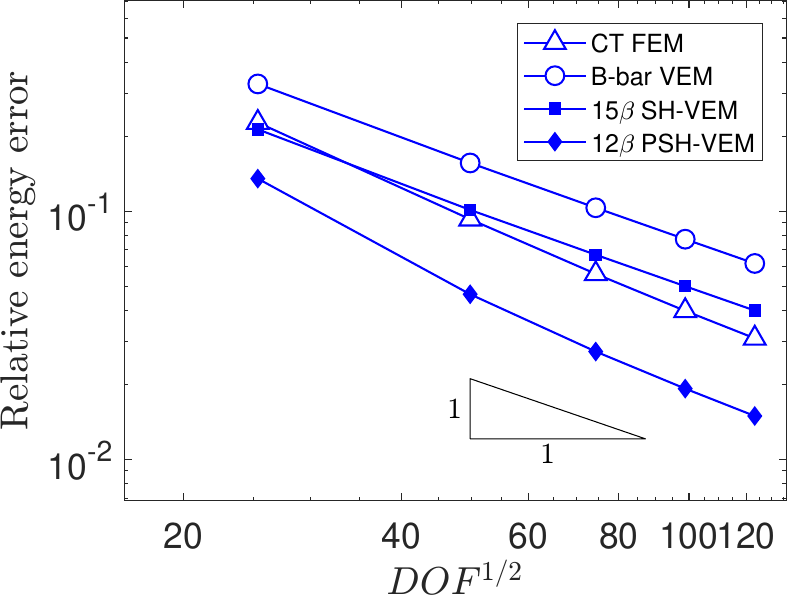}
         \caption{}
     \end{subfigure}
     \hfill
     \begin{subfigure}{0.32\textwidth}
         \centering
         \includegraphics[width=\textwidth]{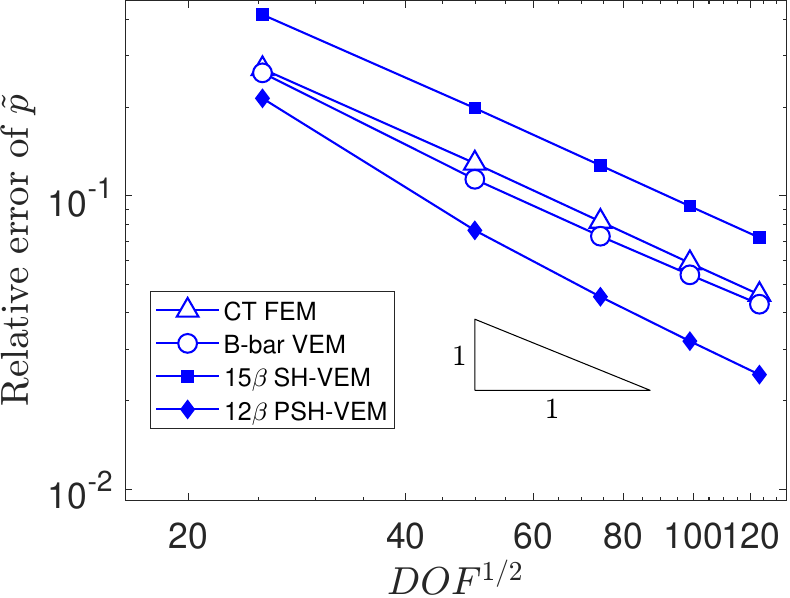}
         \caption{}
     \end{subfigure}
        \caption{Comparison of CT FEM, B-bar VEM, SH-VEM, and PSH-VEM for the manufactured problem on structured meshes (\ac{see~\fref{fig:manufactured_structured}}). (a) $L^2$ error of displacement, (b) energy error, and (c) $L^2$ error of hydrostatic stress. }
        \label{fig:manufactured_convergence_structured}
\end{figure}

\begin{figure}[!bht]
    \centering
    \begin{subfigure}{0.48\textwidth}
    \centering
    \includegraphics[width=\textwidth]{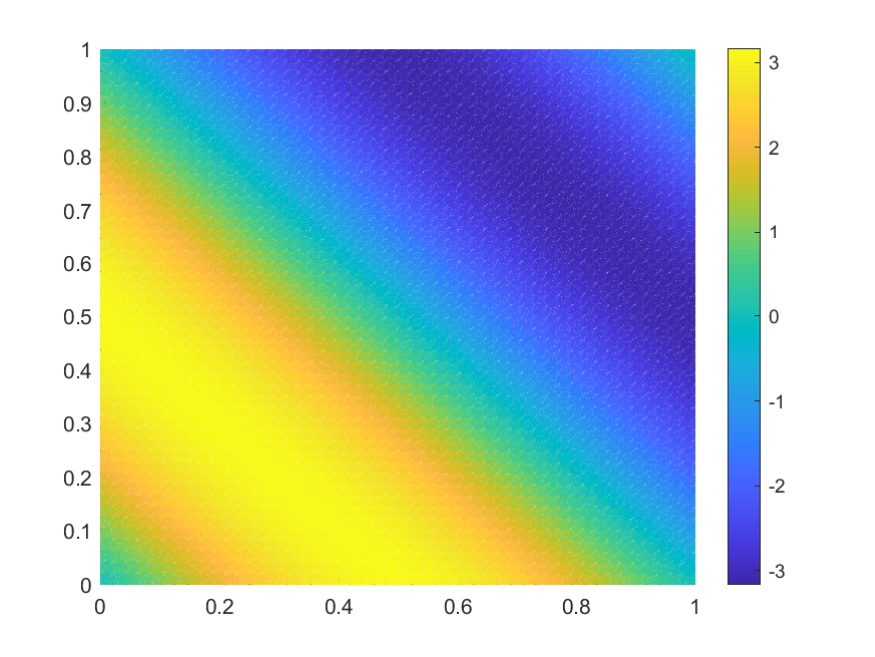}
    \caption{}
    \end{subfigure}
     \hfill
     \begin{subfigure}{0.48\textwidth}
         \centering
         \includegraphics[width=\textwidth]{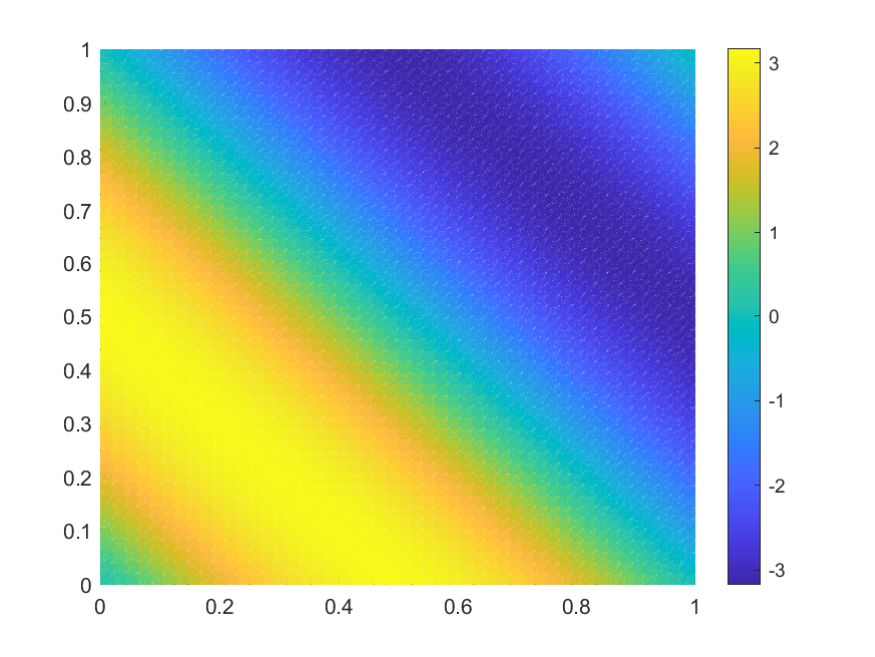}
         \caption{}
     \end{subfigure}
     \vfill
         \begin{subfigure}{0.48\textwidth}
    \centering
    \includegraphics[width=\textwidth]{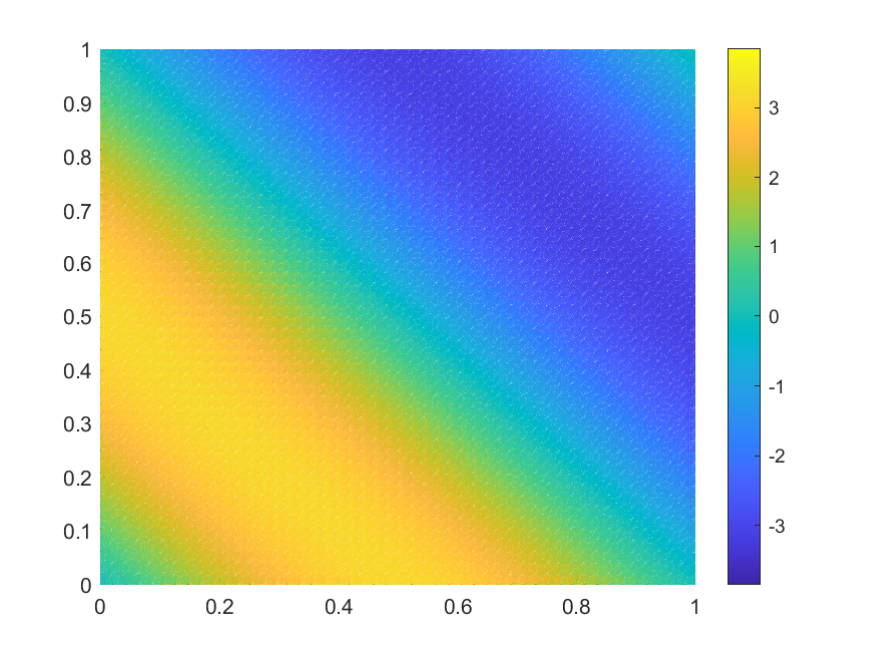}
    \caption{}
    \end{subfigure}
     \hfill
     \begin{subfigure}{0.48\textwidth}
         \centering
         \includegraphics[width=\textwidth]{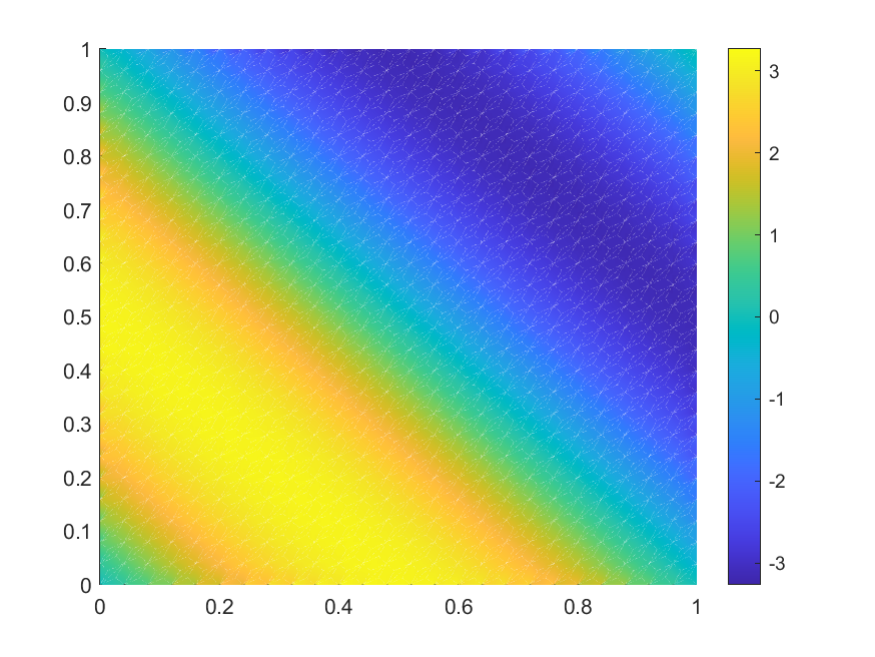}
         \caption{}
     \end{subfigure}
    \caption{Contour plots of the hydrostatic stress on structured meshes (\ac{see~\fref{fig:manufactured_structured}}) for the manufactured problem. (a) CT FEM, (b) B-bar VEM, (c) SH-VEM, and (d) PSH-VEM.}
    \label{fig:manufactured_structured_contour}
\end{figure}

Now, we solve the manufactured problem on a series of structured meshes. The meshes are generated by taking a uniform quadrilateral mesh, then splitting each element along the diagonal into two triangles. Then each triangular element is split into three pieces; a few representative meshes are shown in~\fref{fig:manufactured_structured}. As was the case with the perturbed meshes,~\fref{fig:manufactured_convergence_structured} shows that all four methods deliver optimal convergence rates. The PSH-VEM again has the smallest errors in energy and hydrostatic stress, while the SH-VEM has the largest errors in hydrostatic stress. The contour plots in~\fref{fig:manufactured_structured_contour} are relatively smooth and do not show large errors along the boundary. 

\begin{remark}
    From the sensitivity analysis of the beam and hollow cylinder, we know that if $\alpha$ is three orders of magnitude smaller than $\frac{10^4\ell_0^2}{\ac{E_Y}}$, then the convergence rates in energy seminorm and $L^2$ norm of the hydrostatic stress are affected. For this manufactured problem, when using a Young's modulus $E_Y=1$ psi and $\alpha$ given in~\eqref{eq:penalty_param}, the penalty parameter is set to the minimum value $\alpha = 10\ell_0^2$. This value is three orders of magnitude less than $\frac{10^4\ell_0^2}{\ac{E_Y}}$; therefore, the errors of energy and hydrostatic stress 
    are not expected to have higher order convergence rates even for uniform meshes. However, from our numerical tests, having a larger $E_Y$ or increasing the upper bound on the penalty parameter $\alpha$ results in superconvergent solutions on sufficiently regular meshes. 
\end{remark}

\subsection{Punch problem}
\begin{figure}[!htb]
     \centering
     \begin{subfigure}{0.45\textwidth}
         \centering
         \includegraphics[width=\textwidth]{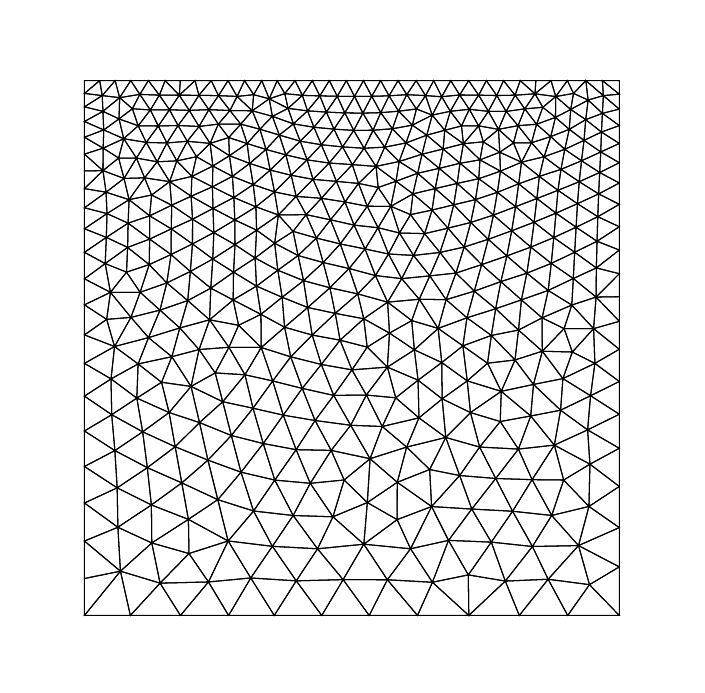}
         \caption{}\label{fig:punch_contour_unstructured_a}
     \end{subfigure}
     \hfill
     \begin{subfigure}{0.48\textwidth}
         \centering
         \includegraphics[width=\textwidth]{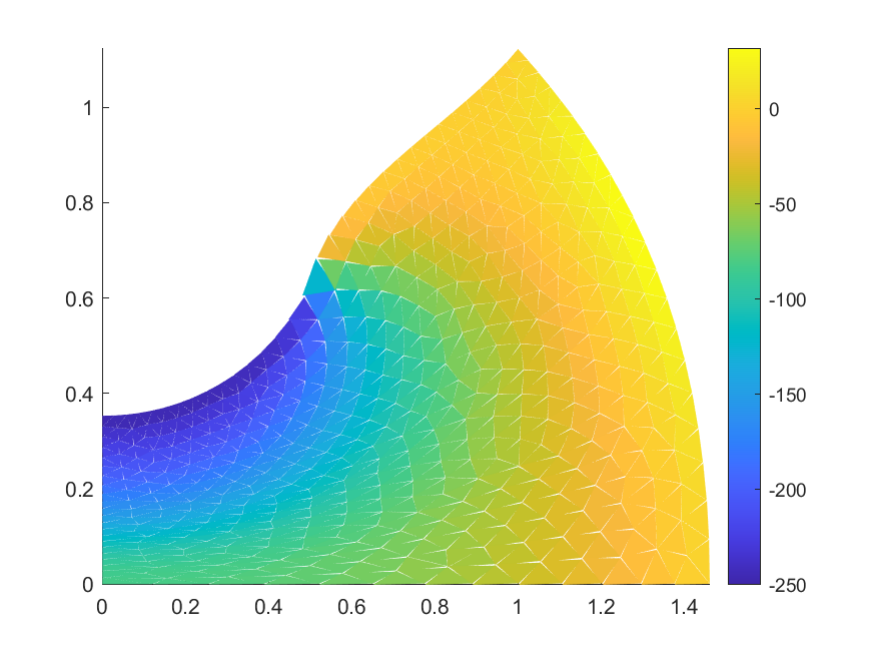}
         \caption{}
     \end{subfigure}
     \vfill
     \begin{subfigure}{0.48\textwidth}
         \centering
         \includegraphics[width=\textwidth]{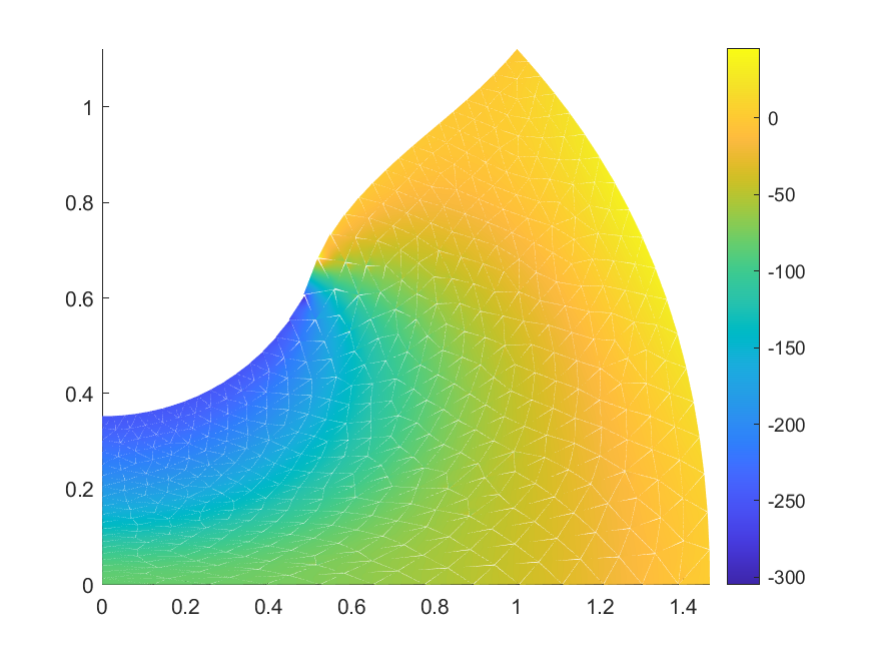}
         \caption{}
     \end{subfigure}
     \hfill
     \begin{subfigure}{0.48\textwidth}
         \centering
         \includegraphics[width=\textwidth]{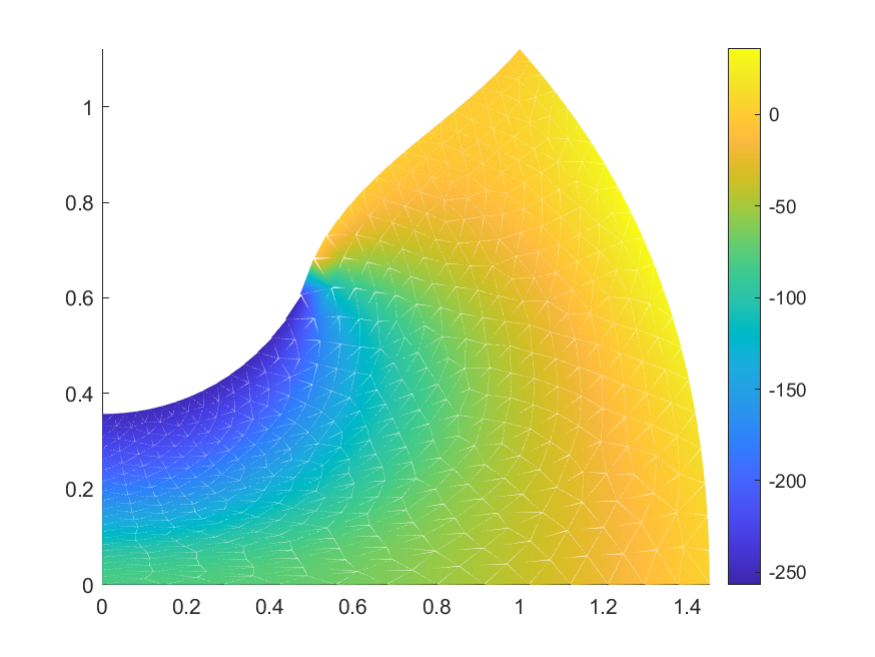}
         \caption{}
     \end{subfigure}
        \caption{(a) Representative unstructured mesh for the punch problem. Contour plots of the hydrostatic stress for the punch problem plotted on the deformed configuration using unstructured meshes. (b) B-bar VEM, (c) SH-VEM, and (d) PSH-VEM. }
        \label{fig:punch_contour_unstructured}
\end{figure}

\begin{figure}[!htb]
     \centering
     \begin{subfigure}{0.48\textwidth}
         \centering
         \includegraphics[width=\textwidth]{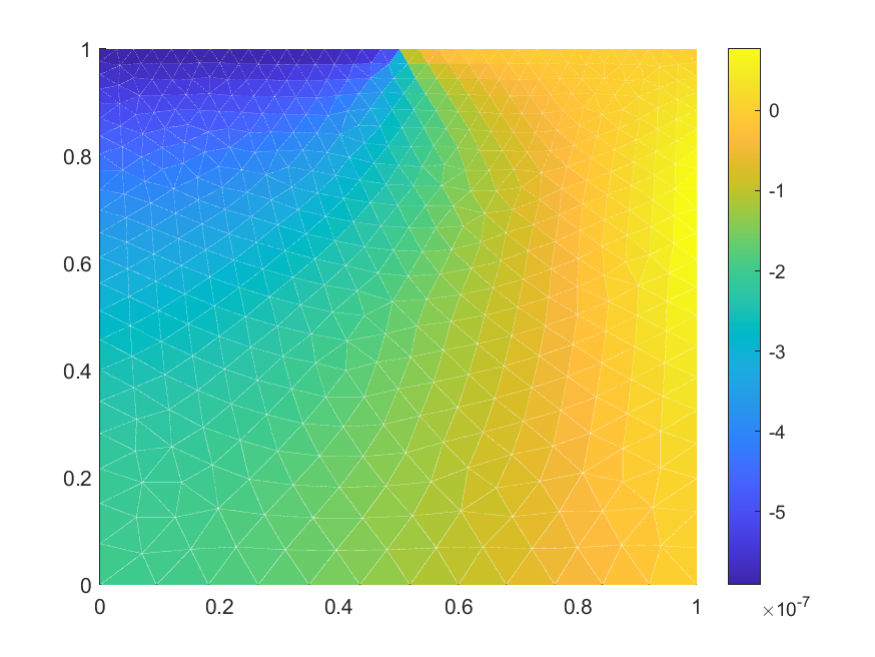}
         \caption{}
     \end{subfigure}
     \hfill
     \begin{subfigure}{0.48\textwidth}
         \centering
         \includegraphics[width=\textwidth]{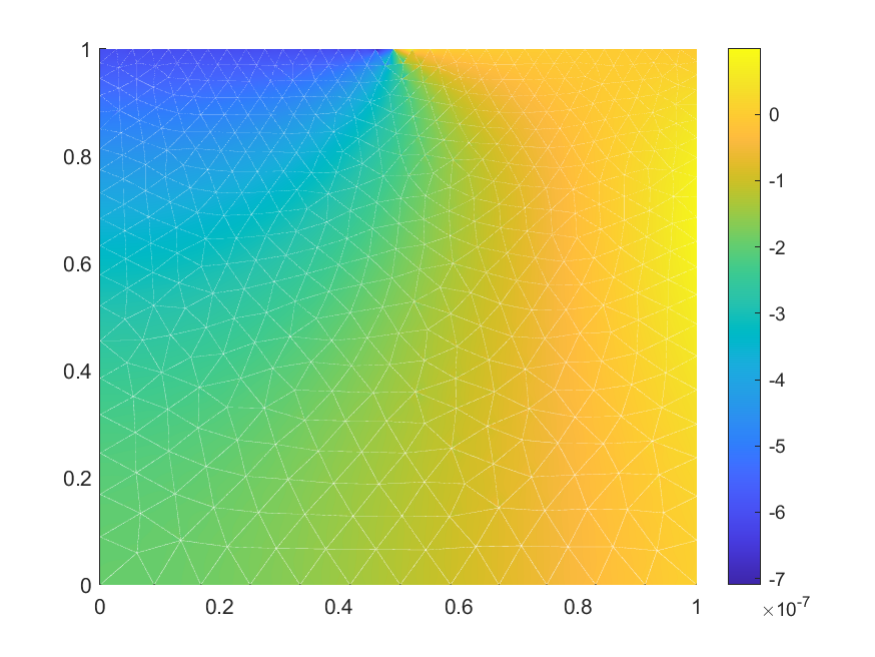}
         \caption{}
     \end{subfigure}
      \hfill
     \begin{subfigure}{0.48\textwidth}
         \centering
         \includegraphics[width=\textwidth]{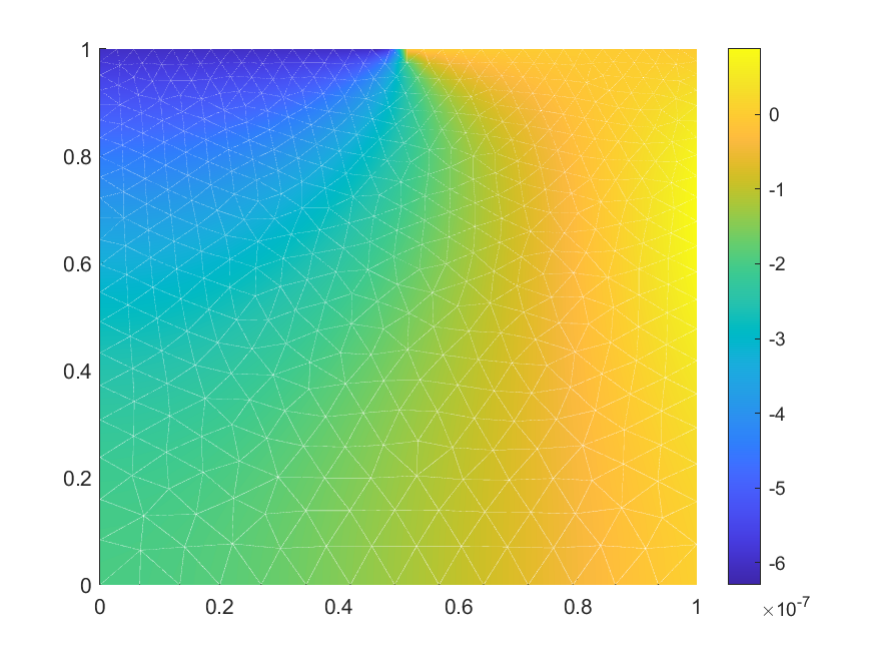}
         \caption{}
     \end{subfigure}
        \caption{Contour plots of the trace of the strain field for the punch problem plotted on the undeformed configuration using unstructured meshes (\ac{see~\fref{fig:punch_contour_unstructured_a}}). (a) B-bar VEM, (b) SH-VEM, and (c) PSH-VEM. }
        \label{fig:punch_trace_contour_unstructured}
\end{figure}

Next, we adapt the problem of a punch being driven into a solid as described in~\cite{Wriggers:2017:EVE} for nearly incompressible hyperelastic materials. This problem is used in~\cite{Wriggers:2017:EVE} to test the robustness of mixed virtual element methods for large deformations and to compare to standard mixed finite element formulations. In our tests, we assume a linearly elastic material and solve the problem on a unit square domain using $E_Y = 250$ psi and $\nu=0.4999999$.   
The left and top edges are horizontally constrained, while the bottom edge is vertically constrained. Along half of the top edge, a uniform load of $F=-250$ lbf per unit length is applied. For this problem, we examine the three methods B-bar VEM, SH-VEM, and PSH-VEM. We first use an unstructured triangular mesh (see~\fref{fig:punch_contour_unstructured_a}) and plot the resulting contours of the hydrostatic stress in~\fref{fig:punch_contour_unstructured}. The contours are plotted on the deformed configuration. From the plots, we find that the three methods produced relatively smooth hydrostatic stress fields, with B-bar and PSH-VEM having a similar range. In~\fref{fig:punch_trace_contour_unstructured}, plots of the trace of the strain field are shown on the undeformed configuration. The three methods yielded nearly traceless strain fields, which is consistent for a nearly incompressible material.  

\begin{figure}[!htb]
     \centering
     \begin{subfigure}{0.45\textwidth}
         \centering
         \includegraphics[width=\textwidth]{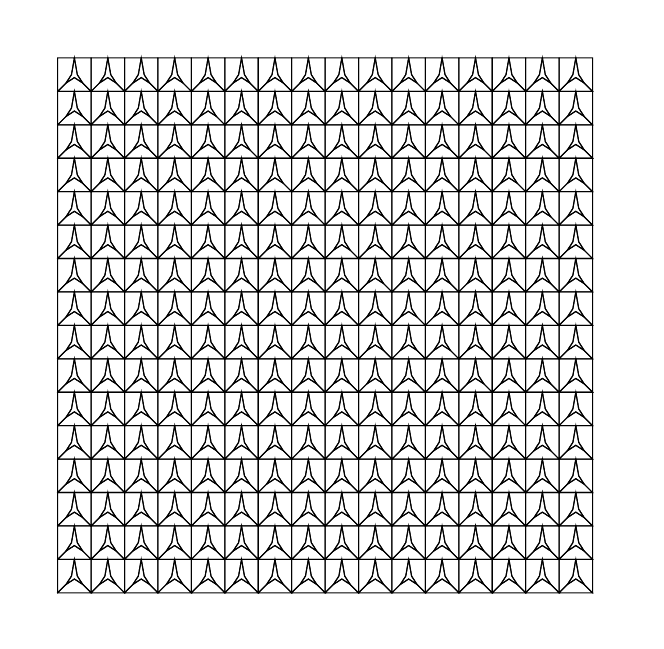}
         \caption{}\label{fig:punch_contour_nonconvex_a}
     \end{subfigure}
     \hfill
     \begin{subfigure}{0.48\textwidth}
         \centering
         \includegraphics[width=\textwidth]{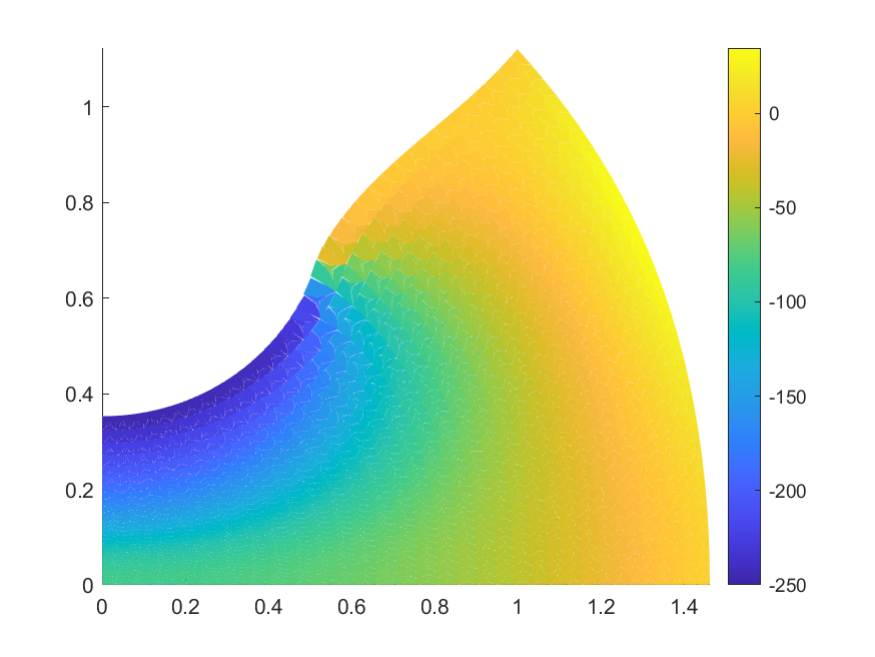}
         \caption{}
     \end{subfigure}
     \vfill
     \begin{subfigure}{0.48\textwidth}
         \centering
         \includegraphics[width=\textwidth]{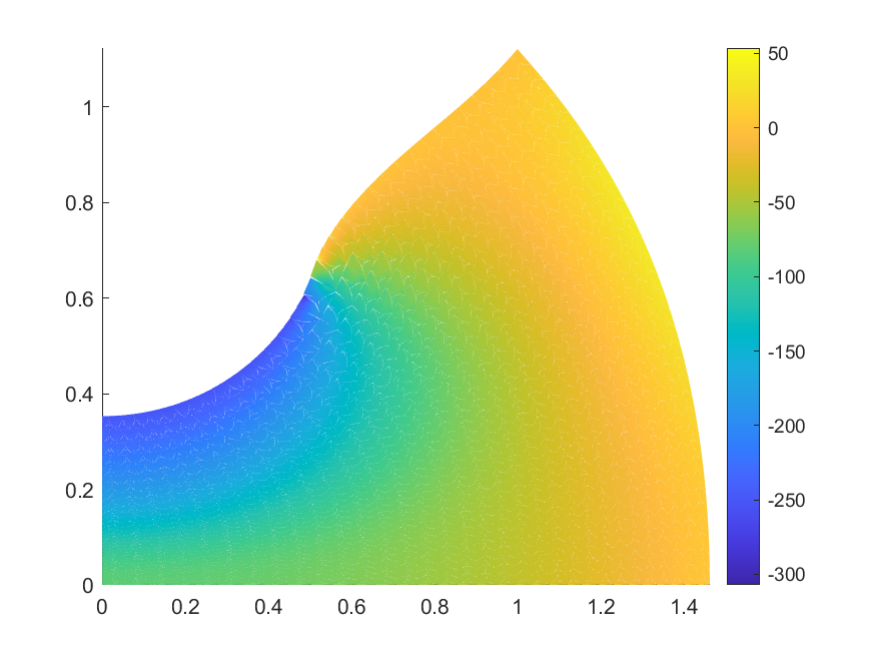}
         \caption{}
     \end{subfigure}
     \hfill
     \begin{subfigure}{0.48\textwidth}
         \centering
         \includegraphics[width=\textwidth]{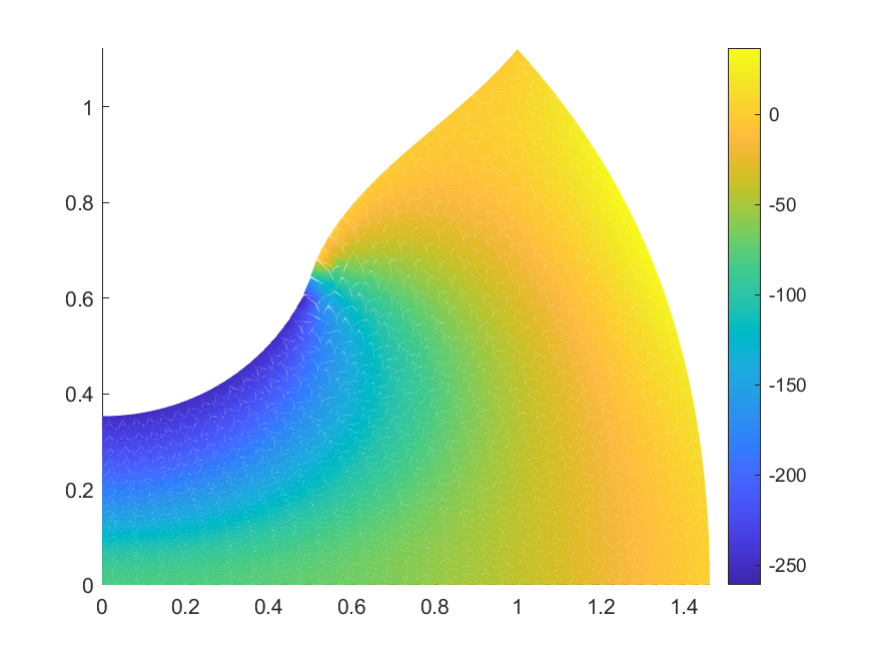}
         \caption{}
     \end{subfigure}
        \caption{(a) Representative nonconvex mesh for the punch problem. Contour plots of the hydrostatic stress of PSH-VEM for the punch problem plotted on the deformed configuration using nonconvex meshes. (b) B-bar VEM, (c) SH-VEM, and (d) PSH-VEM. }
        \label{fig:punch_contour_nonconvex}
\end{figure}

\begin{figure}[!htb]
     \centering
     \begin{subfigure}{0.48\textwidth}
         \centering
         \includegraphics[width=\textwidth]{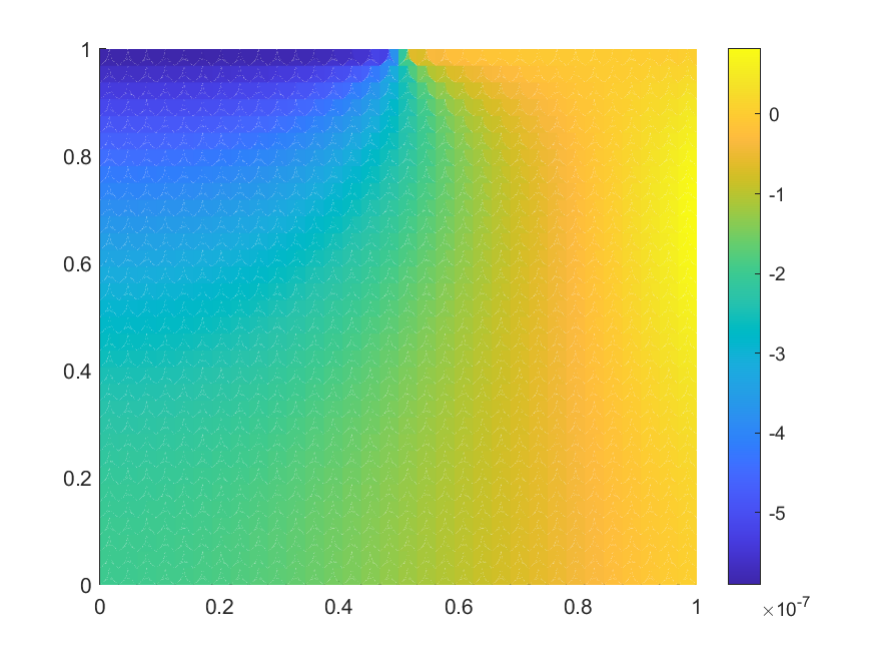}
         \caption{}
     \end{subfigure}
     \hfill
     \begin{subfigure}{0.48\textwidth}
         \centering
         \includegraphics[width=\textwidth]{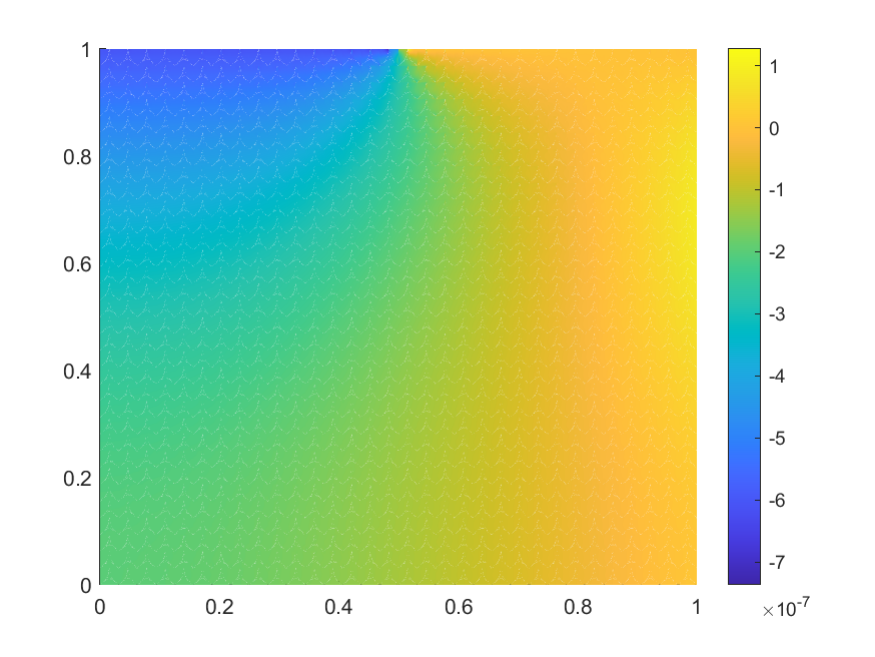}
         \caption{}
     \end{subfigure}
      \hfill
     \begin{subfigure}{0.48\textwidth}
         \centering
         \includegraphics[width=\textwidth]{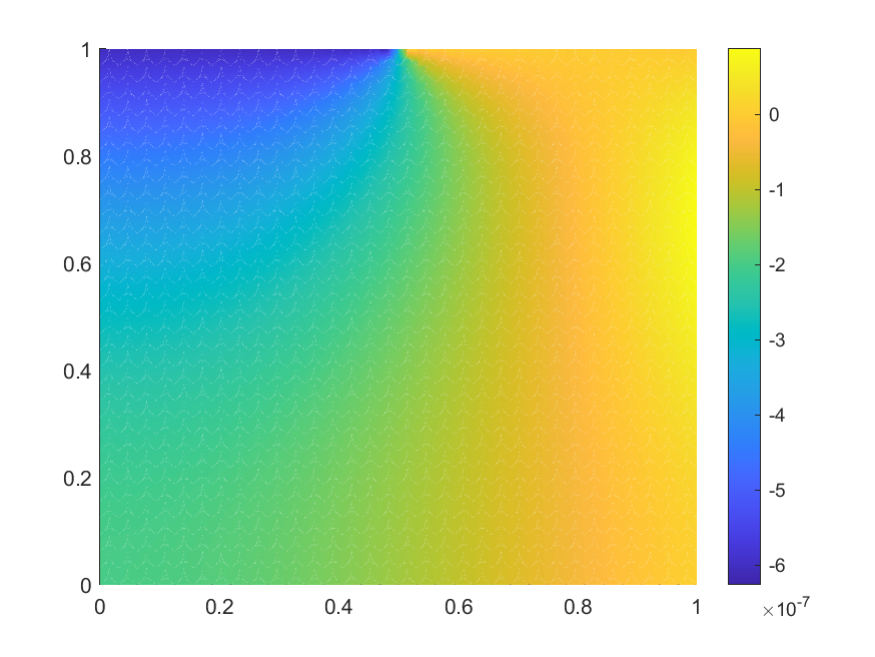}
         \caption{}
     \end{subfigure}
        \caption{Contour plots of the trace of the strain field for the punch problem plotted on the undeformed configuration using nonconvex meshes (\ac{see~\fref{fig:punch_contour_nonconvex_a}}). (a) B-bar VEM, (b) SH-VEM, and (c) PSH-VEM. }
        \label{fig:punch_trace_contour_nonconvex}
\end{figure}

Next, we use a mesh with convex and nonconvex elements with an example mesh shown in~\fref{fig:punch_contour_nonconvex_a}. The contour plots of the hydrostatic stress are given in~\fref{fig:punch_contour_nonconvex}. The plots show that even for nonconvex elements, the three methods retain relatively smooth hydrostatic stress fields. In~\fref{fig:punch_trace_contour_nonconvex}, the contours of the trace of the strain field is presented. Similar to the unstructured case, the strain field of the three methods are nearly traceless.  

\subsection{Stabilized stress-hybrid methods}
\ac{We also examined a stabilized formulation of $9\beta$ and $11\beta$ SH-VEM (both are rank deficient) using a $9$-term and $11$-term divergence-free basis, respectively. We found that the two stabilized stress-hybrid formulations did not show signs of volumetric locking, did not have spurious eigenvalues and attained higher accuracy than $15\beta$ SH-VEM. However, both stabilized methods were less accurate than the PSH method and they did not exhibit superconvergence
in the hydrostatic stress.
Further details can be found in Appendix A.}

\section{Conclusion}\label{sec:conclusion}  
In this paper, we studied an extension of the stress-hybrid virtual element method~\cite{chen:2023:shv} to first-order six-noded triangular virtual elements for linear elastic problems. 
Since complex geometries are amenable to Delaunay mesh generation, 
this work is a first attempt to study two-dimensional virtual element formulations for incompressible elasticity on general triangular meshes (including Delaunay meshes).
In this approach, we used the Hellinger--Reissner variational principle to construct a weak equilibrium and strain-displacement conditions. The weak strain-displacement condition is used along with the virtual element formulation to define a computable projection operator for the stress field. In the initial approach, we used a divergence-free polynomial tensor basis and by combining it with the divergence theorem, we were able to compute the stress projection with only the displacements along the boundary. When using a nondivergence-free basis, we relied on the modification of the VEM space introduced in~\cite{berrone:2023:los} and a secondary energy projection to recover the stress projection. However, we found that this requirement is rather restrictive and leads to elements being overly stiff. Therefore, we relaxed the condition by introducing a penalty term to weakly satisfy the element equilibrium condition and alleviate locking. The SH and PSH-formulations were then tested for stability using an eigenvalue analysis. For the divergence-free basis, a $15$-term expansion was used to ensure no zero-energy modes appeared for highly distorted elements, while a complete bilinear $12$-term expansion was used for the penalty element. For the thin cantilever beam and Cook's membrane problems we found both SH-VEM and PSH-VEM were not sensitive to shear locking. The PSH-VEM was able to reproduce nearly exact bending solutions even on coarse meshes. For the plate with a circular hole and pressurized cylinder the stress-hybrid method produced optimal convergence rates in the $L^2$ norm of the displacement, energy seminorm, and $L^2$ norm of the hydrostatic stress; while the penalty formulation produced superconvergent rates in energy seminorm and $L^2$ norm of the hydrostatic stress. The plate with a hole had relatively smooth hydrostatic stress fields, while for the pressurized cylinder it was observed that large errors concentrated around the interior boundary but the stress-hybrid methods were more accurate. In the manufactured problem, SH-VEM and PSH-VEM showed no signs of locking and had relatively smooth contours of hydrostatic stress. The PSH-VEM showed better accuracy for the manufactured problem; however, for our choice of penalty parameter, the PSH-VEM did not achieve superconvergent rates for the energy seminorm or $L^2$ norm of the hydrostatic stress. For the punch problem, both SH-VEM and PSH-VEM did not show instabilities in the hydrostatic stress field and produced nearly traceless strain fields. In the Appendix, we examined a stabilized version of the $9\beta$ SH-VEM and $11\beta$ SH-VEM. The $9\beta$ and $11\beta$ SH-VEM uses a $9$-term and $11$-term divergence-free basis, respectively. The stress-hybrid methods were then tested on additional benchmark problems. An eigenvalue analysis was used to test stability and volumetric locking of the two stress-hybrid methods. The thick cantilever beam problem showed that the stabilized methods with fewer basis functions were less stiff in bending. The plate with a hole problem showed that the stress-hybrid approaches converged optimally in all cases. The two stabilized methods attained better accuracy than $15\beta$ SH-VEM, with $9\beta$ SH-VEM having performance that is more comparable to the penalized approach. From these tests, we found that the two stabilized stress-hybrid formulations did not show signs of volumetric locking and did not have spurious eigenvalues. 
As part of future work, finding optimal stress basis functions, incorporating additional techniques to mitigate volumetric locking and extending the stress-hybrid and penalty stress-hybrid virtual element method to nonlinear and three dimensional problems are of interest.

\section*{Acknowledgements}
The research support of Sandia National Laboratories to the University of California, Davis, is gratefully acknowledged.
Sandia National Laboratories is a multi-mission laboratory managed and operated by National Technology $\&$ Engineering Solutions of Sandia, LLC (NTESS), a wholly owned subsidiary of Honeywell International Inc., for the U.S. Department of Energy’s National Nuclear Security Administration (DOE/NNSA) under contract DE-NA0003525.  This written work is authored by an employee of NTESS. The employee, not NTESS, owns the right, title, and interest in and to the written work and is responsible for its contents.  Any subjective views or opinions that might be expressed in the written work do not necessarily represent the views of the U.S. Government.  The publisher acknowledges that the U.S. Government retains a non-exclusive, paid-up, irrevocable, world-wide license to publish or reproduce the published form of this written work or allow others to do so, for U.S. Government purposes.  The DOE will provide public access to results of federally sponsored research in accordance with the DOE Public Access Plan.

\appendix
\setcounter{figure}{0}
\setcounter{table}{0}

\section{Stabilized stress-hybrid methods}
In the standard virtual element method, a stabilization term is necessary to ensure that the element stiffness matrix has correct rank. For completeness, we examine effects of a stabilization term on the two methods $9\beta$ SH-VEM and $11\beta$ SH-VEM, which were shown earlier to have zero-energy modes. The $9\beta$ SH-VEM is constructed from the first $9$ terms of the basis given in~\eqref{eq:airy_basis}, while the $11\beta$ SH-VEM uses the first $11$ terms. For the stability term, we follow the construction given in~\cite{Artoli:2017:cmech}:
\begin{subequations}
\begin{align}
    \vm{K}_S &= \tau\left[\vm{I} - \vm{D}(\vm{D^T\vm{D}})^{-1}\vm{D}^T\right],\\
\intertext{where $\tau$ is a scaling factor, $\vm{I}$ is the identity matrix and $\vm{D}\in \mathbb{R}^{12\times 6}$ is the matrix containing the degrees of freedom of the polynomials $\vm{m}_\alpha \in \widehat{\vm{M}}(E)$ given by}
    \vm{D} &= \begin{bmatrix}
        \vm{m}_1(\vx_1) & \vm{m}_2(\vx_1) &\dots & \vm{m}_6(\vx_1) \\ 
        \vm{m}_1(\vx_2) & \vm{m}_2(\vx_2) &\dots & \vm{m}_6(\vx_2) \\
        \dots & \dots & \dots & \dots \\
        \vm{m}_1(\vx_6) & \vm{m}_2(\vx_6) &\dots & \vm{m}_6(\vx_6)
    \end{bmatrix}.\label{eq:dof_matrix}
    \end{align}
\end{subequations}
It is common to choose the scaling factor to be proportional to the trace of the element stiffness matrix $\vm{K}_E$; however, this choice of scaling will lead to an overly stiff solution for nearly incompressible materials~\cite{Park:2020:meccanica,Reddy:2019:AVE}. For simplicity, we set $\tau=\frac{1}{2}$ in the following examples. The stabilized element stiffness matrix is then given by
\begin{equation}
    \vm{K} = \vm{K}_E + \vm{K}_S. 
\end{equation}

\subsection{Eigenvalue analysis}
\begin{figure}[!bht]
    \centering
    \begin{subfigure}{0.48\textwidth}
    \centering
    \includegraphics[width=\textwidth]{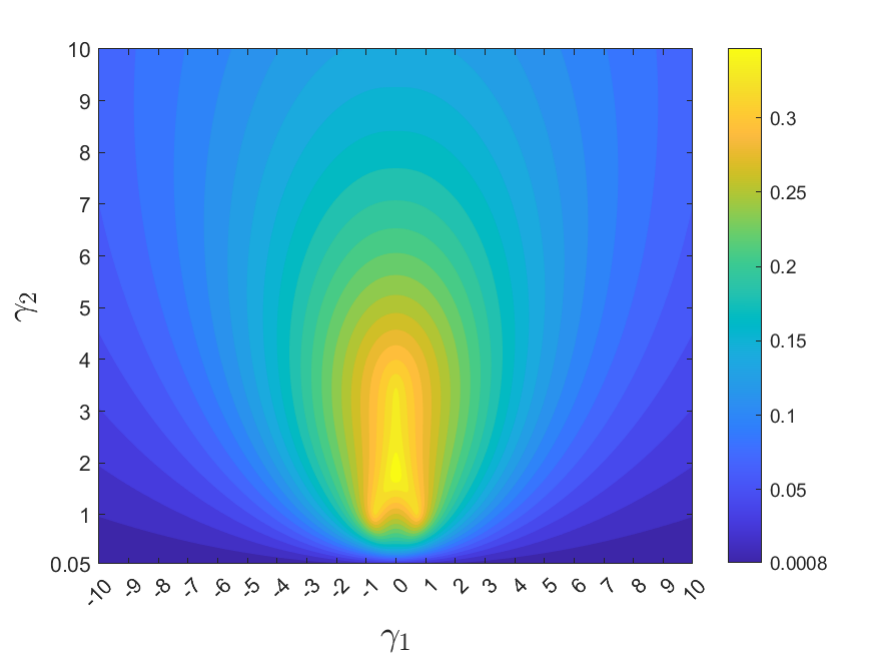}
    \caption{}
    \end{subfigure}
     \hfill
     \begin{subfigure}{0.48\textwidth}
         \centering
         \includegraphics[width=\textwidth]{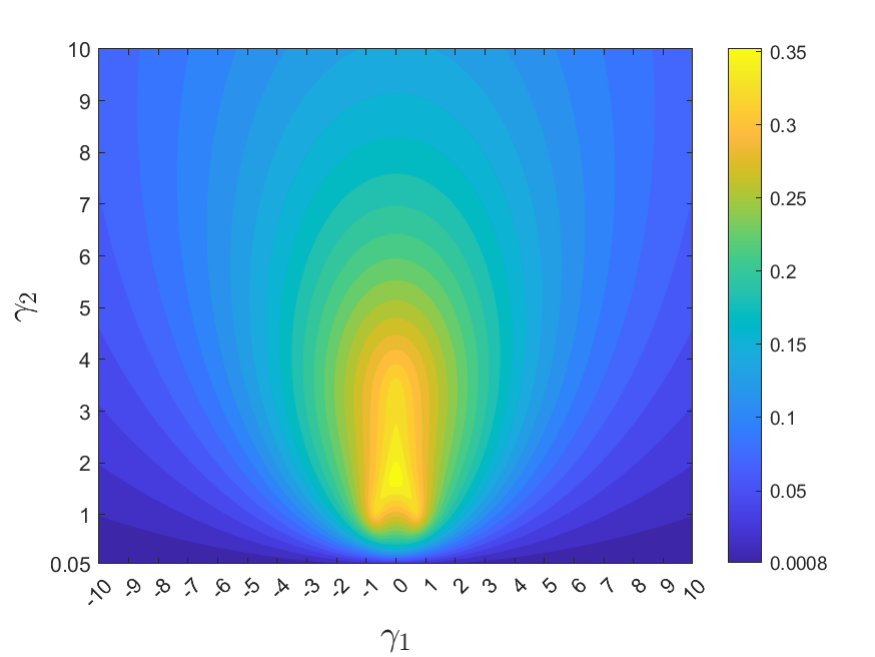}
         \caption{}
     \end{subfigure}
    \caption{Contour plots of the fourth smallest eigenvalue of (a) stabilized $9\beta$ SH-VEM and (b) stabilized $11\beta$ SH-VEM.}
    \label{fig:eig_contour_stabilized}
\end{figure}

We first repeat the eigenvalue analysis presented in~\sref{subsec:eigenanalysis}. Contour 
plots for the fourth smallest eigenvalue are presented in~\fref{fig:eig_contour_stabilized}. The plots reveal that the stabilization term has eliminated the spurious eigenvalue for both $9\beta$ SH-VEM and $11\beta$ SH-VEM. 

Next, we repeat the eigenvalue analysis in~\sref{subsec:eigenanalysis_incompressible} to test for locking in the near-incompressible limit. The five largest eigenvalues of the element stiffness matrix for $9\beta$ and $11\beta$ SH-VEM are given in~\tref{tab:eig_largest_regular_appendix} for a regular six-noded triangular element and in~\tref{tab:eig_largest_nonconvex_appendix} for a nonconvex element. The tables show that both methods only have a single eigenvalue that tends toward infinity. This suggests that the two stabilized elements are not prone to volumetric locking. 
\begin{table}[!bht]
    \centering
    \resizebox{.38\textwidth}{!}{
    \begin{tabular}{|c|c|c|}
    \hline
     Eigenvalue & $9\beta$ SH-VEM & $11\beta$ SH-VEM \\
    \hline
       1&  $7.8 \times 10^{-1}$ & $9.0 \times 10^{-1}$ \\
       2&  $1.2 \times 10^0$& $1.3\times 10^0$\\ 
       3&  $1.3 \times 10^0$& $1.5\times 10^0$\\
       4&  $4.6 \times 10^0$& $4.6\times 10^0$\\
       5&  $4.2 \times 10^{6}$& $4.2\times 10^{6}$\\
         \hline
    \end{tabular}}
    \caption{Comparison of the five largest eigenvalues of the element stiffness matrix on a six-noded triangular element.}
    \label{tab:eig_largest_regular_appendix}
\end{table}
\begin{table}[!bht]
    \centering
    \resizebox{.38\textwidth}{!}{
    \begin{tabular}{|c|c|c|}
    \hline
     Eigenvalue & $9\beta$ SH-VEM & $11\beta$ SH-VEM \\
    \hline
         1& $5.9 \times 10^{-1}$ & $7.9 \times 10^{-1}$ \\
         2& $1.8 \times 10^0$& $1.8\times 10^0$\\ 
         3& $2.1 \times 10^0$& $4.9\times 10^0$\\
         4& $9.9 \times 10^0$& $9.9\times 10^0$ \\
         5& $6.7 \times 10^{6}$& $6.7\times 10^{6}$  \\
         \hline
    \end{tabular}}
    \caption{Comparison of the five largest eigenvalues of the element stiffness matrix on a six-noded nonconvex element.}
    \label{tab:eig_largest_nonconvex_appendix}
\end{table}

\subsection{Bending of a thick cantilever}
\begin{figure}[!htb]
\begin{minipage}{.48\textwidth}
     \centering
     \begin{subfigure}{\textwidth}
         \centering
         \includegraphics[width=\textwidth]{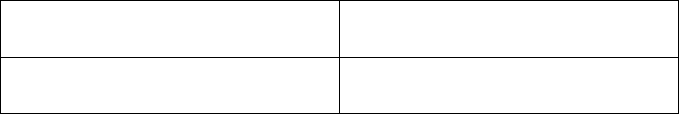}
         \caption{}
     \end{subfigure}
     \vfill
     \begin{subfigure}{\textwidth}
         \centering
         \includegraphics[width=\textwidth]{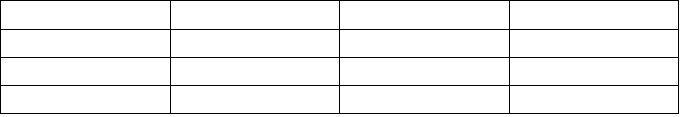}
         \caption{}
     \end{subfigure}
     \vfill
     \begin{subfigure}{\textwidth}
         \centering
         \includegraphics[width=\textwidth]{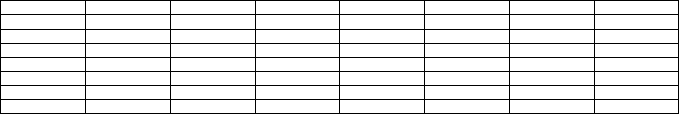}
         \caption{}
     \end{subfigure}
     \end{minipage}
     \hfill
     \begin{minipage}{.48\textwidth}
     \centering
     \begin{subfigure}{\textwidth}
         \centering
         \includegraphics[width=\textwidth]{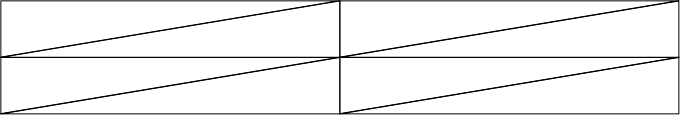}
         \caption{}
     \end{subfigure}
     \vfill
     \begin{subfigure}{\textwidth}
         \centering
         \includegraphics[width=\textwidth]{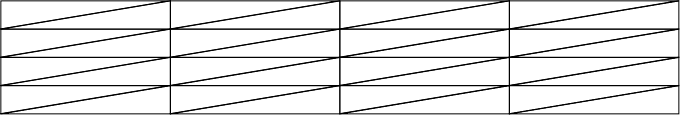}
         \caption{}
     \end{subfigure}
     \vfill
     \begin{subfigure}{\textwidth}
         \centering
         \includegraphics[width=\textwidth]{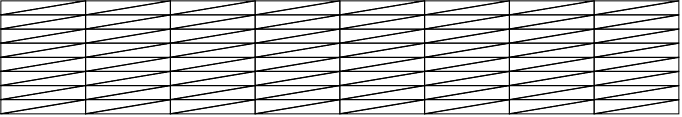}
         \caption{}
     \end{subfigure}
     \end{minipage}
        \caption{Structured meshes for the thick cantilever beam problem. Meshes in (a)-(c) consists of quadrilaterals with an aspect ratio of $4:1$ and (d)-(f) consist of triangular meshes where each corresponding quadrilateral element is cut along the diagonal.}
        \label{fig:beammesh_structured_aspect_ratio_4}
\end{figure}

We consider the bending of a cantilever beam under plane stress conditions with a shear end load as shown in~\cite{Eom:2009:AMP,Felippa:2003:ASO} for different triangular elements. The material has a Young's modulus $E_Y=30000$ psi and Poisson's ratio $\nu=0.25$. The beam has a length of $L=48$ inch, a height of $D=12$ inch, and unit thickness. A shear load of $P=40$ lbf is applied on the right boundary, while the left boundary is fixed. For this problem, we compare the tip displacement of the stress-hybrid methods: $15\beta$ SH-VEM, $12\beta$ PSH-VEM, $9\beta$ SH-VEM, $11\beta$ SH-VEM, and the $5\beta$ SH-VEM for quadrilaterals given in~\cite{chen:2023:shv}. The mesh for the $5\beta$ SH-VEM consists of $N\times N$ $(N=1,2,4,8,16)$ structured quadrilateral elements with an aspect ratio of $4:1$, while the corresponding triangular mesh is constructed by splitting each quadrilateral element along a diagonal (see~\fref{fig:beammesh_structured_aspect_ratio_4}). In~\tref{tab:beam_bending_aspect_4_structured}, the normalized tip displacement for the different methods on each mesh is shown. The table shows that all the SH-VEM and PSH-VEM methods converge with mesh refinement, but the $15\beta$ SH-VEM, which uses higher order terms, is slightly stiffer. The other three methods on six-noded triangles have comparable performance to the $5\beta$ SH-VEM on quadrilaterals.     

\begin{table}[!bht]
    \centering
    \resizebox{\textwidth}{!}{
    \begin{tabular}{|c|c|c|c|c|c|c|}
    \hline
      & $5\beta$ SH-VEM & $15\beta$ SH-VEM & $12\beta$ PSH-VEM & $9\beta$ SH-VEM & $11\beta$ SH-VEM \\
    \hline
        $1\times 1$ & $0.7637$ & $0.4536$ & $1.6185$ & $0.4798$ & $0.4778$  \\
         $2\times 2$ & $0.9413$ & $0.8236$ & $1.0665$ & $0.8668$ & $0.8612$ \\ 
         $4\times 4$ & $0.9856$ & $0.9610$ & $1.0133$ & $0.9787$ & $0.9770$ \\
         $8\times 8$ & $0.9965$ & $0.9917$ & $1.0030$ & $0.9971$ & $0.9966$ \\
         $16\times 16$ & $0.9992$ & $0.9982$ & $1.0007$ & $0.9997$ & $0.9996$\\
         \hline
    \end{tabular}}
    \caption{Comparison of the normalized tip displacements for the cantilever beam problem on structured meshes. The meshes are constructed from $N\times N$ $(N=1,2,4,8,16)$ quadrilaterals with aspect ratio of $4:1$. }
    \label{tab:beam_bending_aspect_4_structured}
\end{table}

\subsection{Plate with a circular hole}
\begin{figure}[!bht]
    \centering
    \begin{subfigure}{0.48\textwidth}
        \centering
        \includegraphics[width=\textwidth]{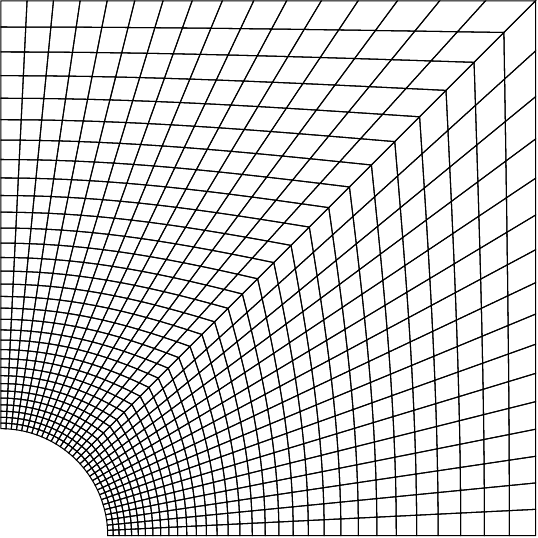}
        \caption{}
    \end{subfigure}
    \hfill
    \begin{subfigure}{0.48\textwidth}
        \centering
        \includegraphics[width=\textwidth]{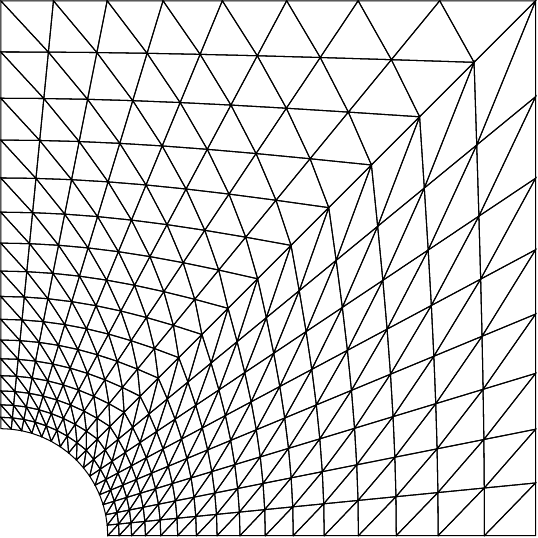}
        \caption{}
    \end{subfigure}
    \caption{(a) A structured quadrilateral mesh and (b) a structured triangular mesh with the same number of degrees of freedom for the plate with a hole problem.  }
    \label{fig:plate_hole_quad_tri_mesh}
\end{figure}

\begin{figure}[!bht]
     \centering
     \begin{subfigure}{0.32\textwidth}
         \centering
         \includegraphics[width=\textwidth]{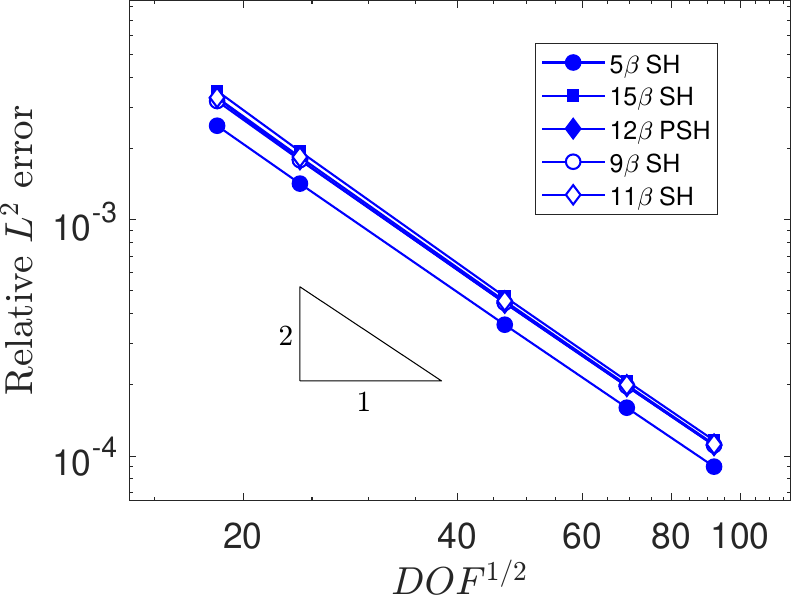}
         \caption{}
     \end{subfigure}
     \hfill
     \begin{subfigure}{0.32\textwidth}
         \centering
         \includegraphics[width=\textwidth]{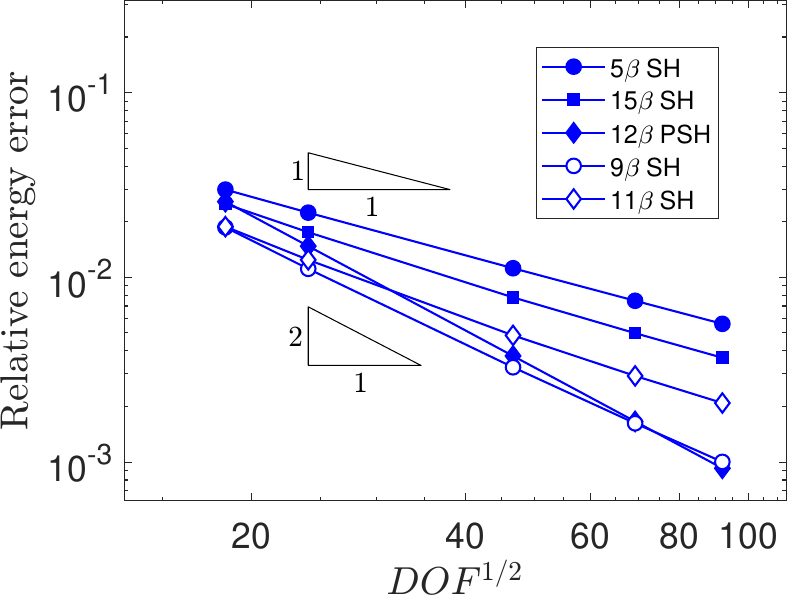}
         \caption{}
     \end{subfigure}
     \hfill
     \begin{subfigure}{0.32\textwidth}
         \centering
         \includegraphics[width=\textwidth]{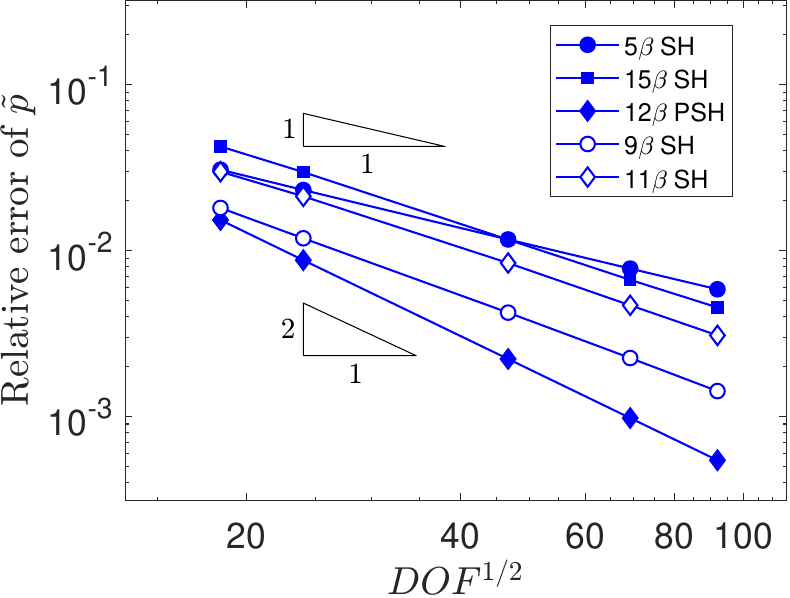}
         \caption{}
     \end{subfigure}
        \caption{Comparison of $5\beta$ SH-VEM, $15\beta$ SH-VEM, $12\beta$ PSH-VEM, $9\beta$ SH-VEM, and $11\beta$ SH-VEM for the plate with a hole problem on structured meshes (\ac{see~\fref{fig:plate_hole_quad_tri_mesh}}). (a) $L^2$ error of displacement, (b) energy error, and (c) $L^2$ error of hydrostatic stress. }
        \label{fig:plate_convergence_sh_compare_structured}
\end{figure}

Finally, we revisit the plate with a circular hole problem that is presented in~\sref{sec:results}. For this problem, we compare the convergence of the stress-hybrid methods: $15\beta$ SH-VEM, $12\beta$ PSH-VEM, $9\beta$ SH-VEM, $11\beta$ SH-VEM, and the $5\beta$ SH-VEM, in the error norms given in~\eqref{eq:error_norms}. We first test this problem on structured quadrilateral and triangular meshes with the same number of global degrees of freedom; a few sample meshes are shown in~\fref{fig:plate_hole_quad_tri_mesh}. In~\fref{fig:plate_convergence_sh_compare_structured}, the convergence results are given for the stress-hybrid methods and show that all the methods converge optimally in the displacement $L^2$ norm, energy seminorm and $L^2$ norm of the hydrostatic stress (the penalty approach attains superconvergence in energy and hydrostatic stress). The $5\beta$ SH-VEM has the smallest error in displacement, but the triangular SH-VEM have better accuracy in energy and hydrostatic stress, with the $12\beta$ PSH-VEM and $9\beta$ SH-VEM having the smallest errors.      

\begin{figure}[!htb]
    \centering
    \begin{subfigure}{0.48\textwidth}
        \centering
        \includegraphics[width=\textwidth]{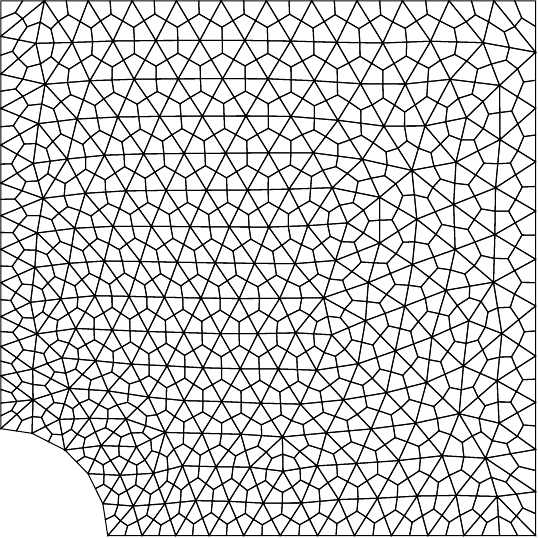}
        \caption{}
    \end{subfigure}
    \hfill
    \begin{subfigure}{0.48\textwidth}
        \centering
        \includegraphics[width=\textwidth]{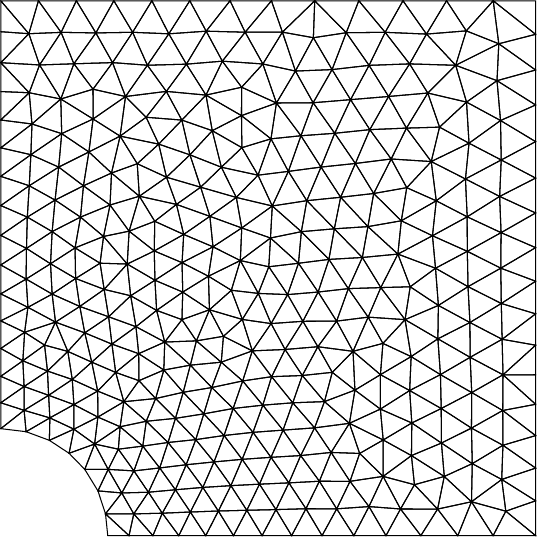}
        \caption{}
    \end{subfigure}
    \caption{(a) An unstructured quadrilateral mesh and (b) an unstructured triangular mesh with a similar number of degrees of freedom for the plate with a hole problem.}
    \label{fig:plate_hole_quad_tri_mesh_unstructured}
\end{figure}

\begin{figure}[!htb]
     \centering
     \begin{subfigure}{0.32\textwidth}
         \centering
         \includegraphics[width=\textwidth]{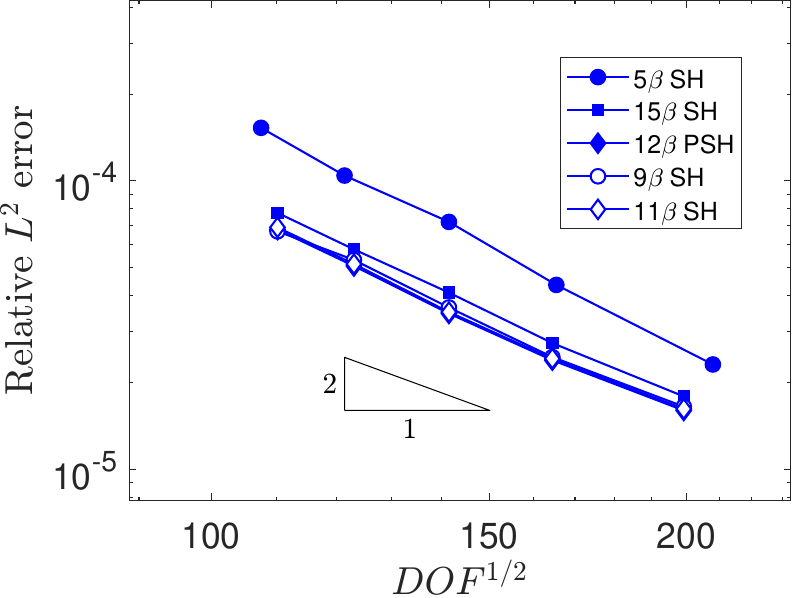}
         \caption{}
     \end{subfigure}
     \hfill
     \begin{subfigure}{0.32\textwidth}
         \centering
         \includegraphics[width=\textwidth]{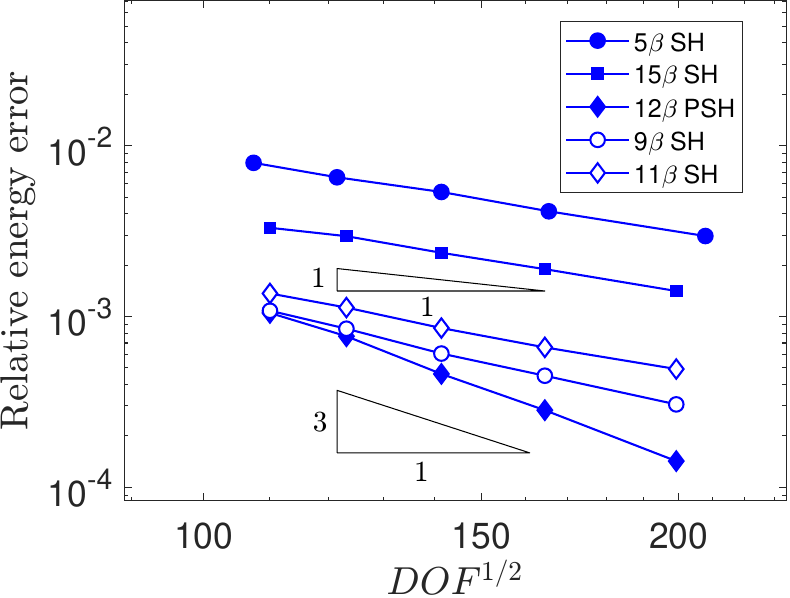}
         \caption{}
     \end{subfigure}
     \hfill
     \begin{subfigure}{0.32\textwidth}
         \centering
         \includegraphics[width=\textwidth]{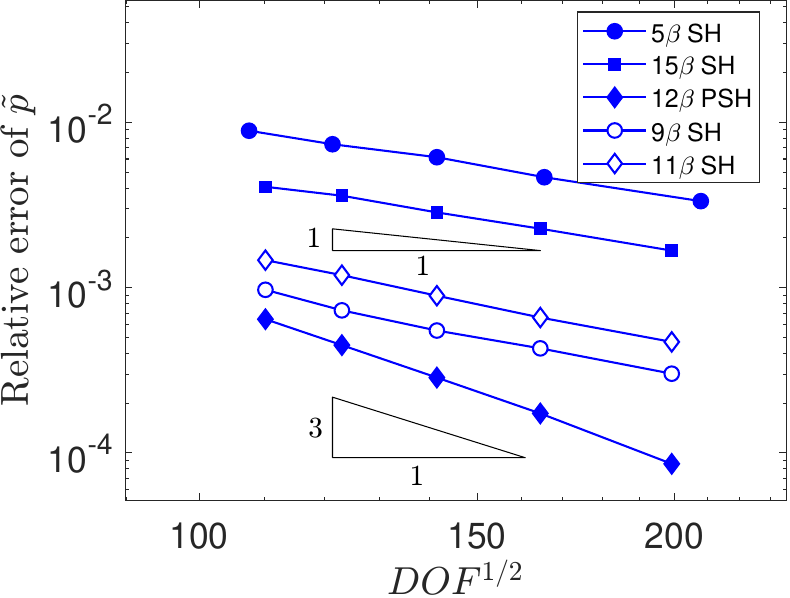}
         \caption{}
     \end{subfigure}
        \caption{Comparison of $5\beta$ SH-VEM, $15\beta$ SH-VEM, $12\beta$ PSH-VEM, $9\beta$ SH-VEM, and $11\beta$ SH-VEM for the plate with a hole problem on unstructured meshes (\ac{see~\fref{fig:plate_hole_quad_tri_mesh_unstructured}}). (a) $L^2$ error of displacement, (b) energy error, and (c) $L^2$ error of hydrostatic stress. }
        \label{fig:plate_convergence_sh_compare_unstructured}
\end{figure}

We now test the plate with a circular hole problem on unstructured meshes. The quadrilateral meshes for $5\beta$ SH-VEM are constructed by taking unstructured triangular meshes and splitting each triangle into three elements. Representative meshes with similar global degrees of freedom are shown in~\fref{fig:plate_hole_quad_tri_mesh_unstructured}. In~\fref{fig:plate_convergence_sh_compare_unstructured}, we present the convergence rates of the five stress-hybrid methods. Like the case of structured meshes, all five methods deliver optimal convergence rates on
unstructured meshes, with $12\beta$ PSH-VEM followed by $9\beta$ SH-VEM delivering the lowest errors in the energy seminorm and $L^2$ norm of hydrostatic stress.

\end{document}